\documentclass[12pt]{article}
\pdfoutput=1
\usepackage{amsbsy,amsmath,amsthm,amssymb,graphicx,verbatim}
\usepackage{float, subfig}
\usepackage[longnamesfirst,sort]{natbib}
\newcommand{\pcite}[1]{\citeauthor{#1}'s \citeyearpar{#1}}
\usepackage{multirow}

\usepackage[letterpaper]{geometry}
\usepackage{amsmath}
\usepackage{amsthm}
\usepackage{amssymb}
\usepackage{amsfonts}
\usepackage{times}
\usepackage{graphicx}
\usepackage{authblk}
\usepackage{enumitem}
\usepackage{color}

\usepackage{ifthen}
\usepackage{xcomment}
\usepackage{versions}
\usepackage[longnamesfirst]{natbib}
\shortcites{JonesEtal:2006, KongEtal:2003}

\newtheoremstyle{mytheorem}{}{}{\slshape}{}{\bfseries}{}{2mm}{}
\theoremstyle{mytheorem}

\newtheorem{theorem}{Theorem}

\theoremstyle{remark}
\newtheorem{remark}{Remark}

\newcommand{\ind}{\buildrel \text{ind} \over \sim}
\newcommand{\Real}{\mathbb{R}}
\newcommand{\cas}{\buildrel \text{a.s.} \over \longrightarrow}

\newcommand{\cd}{\buildrel d \over \rightarrow}
\newcommand{\sX}{\textsf{X}}
\newcommand{\pib}{\bar{\pi}}
\newcommand{\xn}{\{X_n\}_{n\ge 0}}
\newcommand{\hu}{\hat{u}}
\newcommand{\ba}{\boldsymbol{a}}
\newcommand{\baone}{\boldsymbol{a^{[1]}}}
\newcommand{\batwo}{\boldsymbol{a^{[2]}}}

\newcommand{\bd}{\boldsymbol{d}}
\newcommand{\bw}{\boldsymbol{w}}
\newcommand{\hatbd}{\boldsymbol{\hat{d}}}
\newcommand{\bm}{\boldsymbol{m}}

\newcommand{\bzeta}{\boldsymbol{\zeta}}
\newcommand{\hatbzeta}{\boldsymbol{\hat{\zeta}}}
\newcommand{\bt}{\boldsymbol{t}}
\newcommand{\bz}{\boldsymbol{z}}
\newcommand{\by}{\boldsymbol{y}}
\newcommand{\E}{E}
\newcommand{\baoneopt}{\boldsymbol{\hat{a}}^{[1, \text{opt}]}}

\newcommand{\eb}{\boldsymbol{b}}
\renewcommand{\mid}{\,|\,}
\newcommand{\m}{M}
\newcommand{\ap}{\alpha}

\graphicspath{{./bvs-cookie/}}
\graphicspath{{./t-figures/}}
\newcount\dotcnt\newdimen\deltay
\def\Ddot#1#2(#3,#4,#5,#6){\deltay=#6\setbox1=\hbox to0pt{\smash{\dotcnt=1
\kern#3\loop\raise\dotcnt\deltay\hbox to0pt{\hss#2}\kern#5\ifnum\dotcnt<#1
\advance\dotcnt 1\repeat}\hss}\setbox2=\vtop{\box1}\ht2=#4\box2}

\DeclareMathOperator{\Cov}{Cov}

\numberwithin{equation}{section}
\begin{document}

\title{Estimating standard errors for importance sampling estimators with
multiple Markov chains}
\author[1]{Vivekananda Roy}
\author[2]{Aixin Tan}
\author[3]{James M. Flegal}
\affil[1]{Department of Statistics, Iowa State University}
\affil[2]{Department of Statistics and Actuarial Science, University of Iowa}
\affil[3]{Department of Statistics, University of California, Riverside}
\date{Aug 10, 2016}

\maketitle
\renewcommand{\baselinestretch}{1.2}\large\normalsize

%
%\renewcommand{\baselinestretch}{1.2}
%
%\markright{ \hbox{\footnotesize\rm Statistica Sinica
%%{\footnotesize\bf 24} (201?), 000-000
%}\hfill\\[-13pt]
%\hbox{\footnotesize\rm
%%\href{http://dx.doi.org/10.5705/ss.20??.???}{doi:http://dx.doi.org/10.5705/ss.20??.???}
%}\hfill }
%
%\markboth{\hfill{\footnotesize\rm VIVEKANANDA ROY,  AIXIN TAN, AND JAMES M. FLEGAL} \hfill}
%{\hfill {\footnotesize\rm IMPORTANCE SAMPLING WITH MULTIPLE MARKOV CHAINS} \hfill}
%
%\renewcommand{\thefootnote}{}
%$\ $\par

%%%%%%%%%%%%%%%%%%%%%%%%%%%%%%%%%%%%%%%%%%%%%%%%%%%%%%%%%%%%%%%%%%%%%%%%%%%%%%%%%%%%%%%%%%%%%%%%%%%%%%%%%%%%%%%%%%%%%%%%%%%%

%\fontsize{10.5}{14pt plus.8pt minus .6pt}\selectfont
%\vspace{0.8pc}
%\centerline{\large\bf ESTIMATING STANDARD ERRORS FOR IMPORTANCE SAMPLING }
%\vspace{2pt}
%\centerline{\large\bf ESTIMATORS  WITH
%MULTIPLE MARKOV CHAINS}
%\vspace{.4cm}
%\centerline{Vivekananda Roy, Aixin Tan, and James M. Flegal}
%\vspace{.4cm}
%\centerline{\it Iowa State University, University of Iowa, and University of California, Riverside}
%\vspace{.55cm}
%\fontsize{9}{11.5pt plus.8pt minus .6pt}\selectfont

%%%%%%%%%%%%%%%%%%%%%%%%%%%%%%%%%%%%%%%%%%%%%%%%%%%%%%%%%%%%%%%%%%%%%%%%%%%%%%%%%%%%%%%%%%%%%%%%%%%%%%%%%%%%%%%%%%%%%%%%%%%%

%\begin{quotation}
%\noindent {\it Abstract:}
\begin{abstract}
The naive importance sampling estimator, based on samples from a single importance density, can be numerically unstable. Instead, we consider generalized importance sampling estimators where samples from more than one probability distribution are combined.  We study this problem in the Markov chain Monte Carlo context, where independent samples are replaced with Markov chain samples. If the chains converge to their respective target distributions at a polynomial rate, then under two finite moment conditions, we show a central limit theorem holds for the generalized estimators. Further, we develop an easy to implement method to calculate valid asymptotic standard errors based on batch means. We also provide a batch means estimator for calculating asymptotically valid standard errors of \pcite{Geyer:1994} reverse logistic estimator. We illustrate the method via three examples.  In particular, the generalized importance sampling estimator is used for Bayesian spatial modeling of binary data and to perform empirical Bayes variable selection where the batch means estimator enables standard error calculations in high-dimensional settings.
\end{abstract}

\thanks{
%\noindent %
}

\section{Introduction}
Let $\pi(x) = \nu(x)/m$ be a probability density function (pdf) on
$\sX$ with respect to a measure $\mu(\cdot)$. Suppose $f: \sX
\rightarrow \mathbb{R}$ is a $\pi$ integrable function and we want to
estimate $E_\pi f := \int_{\sX} f(x) \pi(x) \mu (dx)$. Let $\pi_1(x) =
\nu_1(x)/m_1$ be another pdf on $\sX$ such that $\{x: \pi_1(x) = 0\}
\subset \{x: \pi(x) = 0\}$. The importance sampling (IS) estimator of
$E_\pi f$ based on independent and identically distributed (iid) samples $X_1, \dots, X_n$ from the importance
density $\pi_1$ is % (see e.g.\ \citet{robe:case:2004})
\begin{equation}
  \label{eq:is0}
\frac{\sum_{i=1}^n f(X_i)\nu(X_i)/ \nu_1(X_i)} { \sum_{i=1}^n \nu(X_i)/ \nu_1(X_i)}  \cas  \int_{\sX} \frac{f(x) \nu(x) / m} {\nu_1(x) / m_1}
  \pi_1(x) \, \mu(dx) \bigg/ \int_{\sX} \frac{\nu(x) / m} {\nu_1(x)
  / m_1} \pi_1(x) \, \mu(dx) = E_\pi f,
\end{equation}
as $n \to \infty$.  This estimator can also be used in the Markov chain Monte Carlo (MCMC) context when $X_1, \dots, X_n$ are realizations from a suitably irreducible Markov chain with stationary density $\pi_1$ (\citet{hast:1970}). Note that \eqref{eq:is0} requires the functions $\nu$ and $\nu_1$ to be known. On the other hand, it does not depend on normalizing constants $m$ and $m_1$, which are generally unknown.

In this article, we consider situations where one wants to estimate
$E_\pi f$ for all $\pi$ belonging to a large collection, say
$\Pi$. This situation arises in both frequentist and Bayesian
statistics. Although \eqref{eq:is0} provides consistent estimators of
$E_\pi f$ for all $\pi \in \Pi$ based on a single Markov chain $\xn$
with stationary density $\pi_1$, it does not work well when $\pi$
differs greatly from $\pi_1$.  In that case the ratios
$\nu(x)/\nu_1(x)$ can be arbitrarily large for some sample values
making the estimator at \eqref{eq:is0} unstable.  In general, there is
not a single good importance density $\pi_1$ which is close to all
$\pi \in \Pi$ (see e.g.\ \citet{Geyer:1994}). 
Hence a natural
modification is to replace $\pi_1$ in \eqref{eq:is0} with a mixture of
densities where each density in $\Pi$ is close to a subset of the $k$ reference densities. 
%Hence a natural
%modification is to replace $\pi_1$ in \eqref{eq:is0} with a mixture of
%densities where each density is ``close'' to a subset of densities in $\Pi$.
To this end, denote $\overline{\pi} \equiv \sum_{i=1}^k (a_i/ |\ba|)
\pi_i$, where $\ba= (a_1, \dots, a_k)$ are $k$ positive constants,
$|\ba| = \sum_{i=1}^k a_i$, and $\pi_i(x) = \nu_i(x)/m_i$ for $i=1,
\dots, k$ are $k$ densities known up to their normalizing
constants. Suppose further that $n_1, \dots, n_k$ are positive
integers and $d_i := m_i/m_1$ for $i= 2, \dots, k,$ with $d_1 \equiv
1$. Then define the $(k-1)$ dimensional vector
\begin{equation}
\label{eq:ddef}
  \bd = (m_2 / m_1, \ldots, m_k / m_1).
\end{equation}
Finally for $l=1, \dots, k$, let $\{X_i^{(l)}\}_{i=1}^{n_l}$ be an iid sample from $\pi_l$ or realizations from a positive Harris Markov chain with invariant density $\pi_l$ (for definitions see \citet{MeynTweedie:1993}). Then as $n_l
\rightarrow \infty,$ for all $l=1, \dots, k$, we have
\begin{align}
  \label{eq:mainest}
  \hat{\eta} &\equiv  \Biggl( \sum_{l=1}^k
  \frac{a_l} {n_l} \sum_{i=1}^{n_l} \frac{f(X_i^{(l)})
  \nu(X_i^{(l)})} {\sum_{s=1}^k a_s \nu_s(X_i^{(l)}) / d_s} \Biggr)
  \Bigg/ \Biggl( \sum_{l=1}^k \frac{a_l} {n_l} \sum_{i=1}^{n_l}
  \frac{\nu(X_i^{(l)})} {\sum_{s=1}^k a_s \nu_s(X_i^{(l)}) / d_s}
  \Biggr)\\ & \cas \Biggl( \sum_{l=1}^k a_l \int_{\sX} f(x) \frac{\nu(x)}
           {\sum_{s=1}^k a_s \nu_s(x) / d_s} \pi_l(x) \, \mu(dx)
           \Biggr) \Bigg/ \Biggl( \sum_{l=1}^k a_l \int_{\sX}
           \frac{\nu(x)} {\sum_{s=1}^k a_s \nu_s(x) / d_s} \pi_l(x)
           \, \mu(dx) \Biggr) \nonumber \\ &= \int_{\sX} f(x)
           \frac{\nu(x)} {\pib(x)} \pib(x) \, \mu(dx) \bigg/
           \int_{\sX} \frac{\nu(x)} {\pib(x)} \pib(x) \, \mu(dx) = E_\pi f \nonumber.
\end{align}

The generalized IS estimator \eqref{eq:mainest} has been discussed widely in the literature, e.g.\ applications include Monte Carlo maximum likelihood estimation and Bayesian sensitivity analysis. \citet{GillVardiWellner:1988}, \citet{KongEtal:2003}, \citet{MengWong:1996}, \citet{Tan:2004}, and \citet{Vardi:1985} consider estimation using \eqref{eq:mainest} based on iid samples. The estimator is applicable to a much larger class of problems if Markov chain samples are allowed, see e.g.\ \citet{buta:doss:2011}, \citet{Geyer:1994}, and \citet{tan:doss:hobe:2015}, which is the setting of this paper.  

Alternative importance weights have also been proposed.  In the case when the normalizing constants $m_i$'s are known, the estimator \eqref{eq:mainest} resembles the {\it balance heuristic} estimator of \cite{veac:guib:1995}, which is revisited in \cite{owen:zhou:2000} as the {\it deterministic mixture}. The standard population Monte Carlo algorithm of \cite{capp:guil:mari:robe:2004} uses a weighted ratio of the target $\pi$ and the proposal $\pi_j$ it was drawn from (evaluated at the sample itself). However, when iid samples are available from $\pi_j, j=1, 2, \dots, k$, \cite{elvi:mart:luen:buga:2015} have shown that the normalized estimator ($m_i$'s known) version of \eqref{eq:mainest} always has a smaller variance than that of the population Monte Carlo algorithm. Further, it may be difficult in practice to find fully known importance densities that approximate the target densities. Indeed, applications such as in empirical Bayes analysis and Bayesian sensitivity analysis routinely select representatives from the large number of target posterior densities to serve as proposal densities, and they are known only up to normalizing constants. See \citet{buta:doss:2011}, \citet{doss:2010}, as well as section~\ref{sec:spatial} for examples.
%\red{Further, it may be difficult in practice to find fully known importance densities that approximate the target densities.  Specifically, applications such as in empirical Bayes analysis and Bayesian sensitivity analysis routinely use target densities as proposal densities that are known only up to normalizing constants \citep{buta:doss:2011, doss:2010}.}  
Although there is no known proof for the self normalized estimator \cite[][p.\ 16]{elvi:mart:luen:buga:2015}, it is reasonable to assume the superiority of \eqref{eq:mainest} over estimators corresponding to other weighting schemes.

As noted in \eqref{eq:mainest}, the estimator $\hat{\eta}$ converges to $E_\pi f$ as the sample sizes increase to infinity, for iid samples as well as Markov chain samples satisfying the usual regularity conditions. Now given samples of finite size in practice, it is of fundamental importance to provide some measure of uncertainty, such as the standard errors (SEs) associated with this consistent estimator. Take estimators that are sample averages based on iid Monte Carlo samples for example, it is a basic requirement to report their SEs. But the very same issue is often overlooked in practice when the estimators have more complicated structure, and when they are based on MCMC samples, largely due to the difficulty of doing so. See, for e.g. \citet{fleg:hara:jone:2008} on the issue concerning MCMC experiments and \citet{koeh:etal:2009} for more general simulation studies. 
%We believe the relatively small literature in the computing community on reporting SE for MCMC estimators is largely due to the complexity of the problem. 
For calculating SEs of $\hat{\eta}$ based on MCMC samples, \citet{tan:doss:hobe:2015} provide a solution using the method of regenerative simulation (RS).  However, this method crucially depends on the construction of a practical minorization condition, i.e.\ one where sufficient regenerations are observed in finite simulations (for definitions and a description of RS see \citet{mykl:tier:yu:1995}).  Further, the usual method of identifying regeneration times by splitting becomes impractical for high-dimensional problems \citep{gilks1998adaptive}. 
% Further, RS is problematic for multi-dimensional updates with Metropolis steps since a regeneration can only occur if all the components are accepted \citep{mykl:tier:yu:1995}.  
Hence, successful applications of RS involve significant trial and error and are usually limited to low-dimensional Gibbs samplers (see e.g.\ \citet{tan:hobe:2009, roy:hobe:2007}).  In this paper we avoid RS and provide SE estimators of $\hat{\eta}$ using the batch means (BM) method, which is straightforward to implement and can be routinely applied in practice. In obtaining this estimator, we also establish a central limit theorem (CLT) for $\hat{\eta}$ that generalizes some results in \citet{buta:doss:2011}. % See the remarks following Theorem~\ref{thm:elex} for more discussions on how the results in this paper improve on those given in \citet{tan:doss:hobe:2015}.

The estimator $\hat{\eta}$ in \eqref{eq:mainest} depends on the ratios of normalizing constants, $\bd$, which are unknown in practical applications. We consider the two-stage scheme studied in \citet{buta:doss:2011} where first an estimate $\hatbd$ is obtained using \pcite{Geyer:1994} ``reverse logistic regression'' method based on samples from $\pi_l$, and then independently, new samples are used to estimate $E_\pi f$ for $\pi \in \Pi$ using the estimator $\hat{\eta} (\hatbd)$ in \eqref{eq:mainest}. \citet{buta:doss:2011} showed that the asymptotic variance of $\hat{\eta} (\hatbd)$ depends on the asymptotic variance of $\hatbd$. Thus we study the CLT of $\hatbd$ and provide a BM estimator of the asymptotic covariance matrix of $\hatbd$. Since $\hatbd$ involves multiple Markov chain samples, we utilize a multivariate BM estimator. Although, the form of the asymptotic covariance matrix of $\hatbd$ is complicated, our consistent BM estimator is straightforward to code.

The problem of estimating $\bd$, the ratios of normalizing constants of unnormalized densities is important in its own right and has many applications in frequentist and Bayesian inference.  For example, when the samples are iid sequences this is the biased sampling problem studied in \citet{Vardi:1985}.  In addition, the problem arises naturally in the calculations of likelihood ratios in missing data (or latent variable) models, mixture densities for use in IS, and Bayes factors.  

Our work considers the problem of estimating $\bd$ using \pcite{Geyer:1994} reverse logistic regression method. Specifically, we study the general quasi-likelihood function proposed in \citet{DossTan:2014}. Unlike \pcite{Geyer:1994} method, this extended quasi-likelihood function has the advantage of using user defined weights which are appropriate in situations where the multiple Markov chains have different mixing rates. We establish the CLT for the resulting estimators of $\bd$ and develop the BM estimators of their asymptotic covariance matrix.

Thus we consider two related problems in this paper--firstly, estimating (ratios of) normalizing constants given samples from $k$ densities, and secondly, estimating expectations with respect to a large number of (other) target distributions using these samples. In both cases, we establish CLTs for our estimators and provide easy to calculate SEs using BM methods. 

Prior results of \citet{buta:doss:2011}, \citet{DossTan:2014}, \citet{Geyer:1994}, and \citet{tan:doss:hobe:2015} all assume that the
underlying Markov chains are geometrically ergodic.  We weaken this condition substantially in that we only require the chains to be {\it polynomial ergodic}.  To this end, let $K_l(x, \cdot)$ be the Markov transition function for the Markov chain $\Phi_l = \{X_t^{(l)}\}_{t
  \ge 1}$, so that for any measurable set $A$, and $s, t \in \{1, 2, \ldots\}$ we have $P\bigl( X_{s+t}^{(l)} \in A \mid X_s^{(l)} = x
\bigr)= K_l^{t} (x, A)$. Let $\| \cdot \|$ denote the total variation norm and $\Pi_l$ be the probability measure corresponding to the density $\pi_l$. The Markov chain $\Phi_l$ is {\it polynomially ergodic of order $m$} where $m >0$ if there exists $W : \sX \rightarrow \mathbb{R}^+$ with $E_{\pi_l} W < \infty$ such that
\begin{equation*}
  \|K_l^{t} (x, \cdot) - \Pi_l(\cdot)\| \le W(x) t^{-m} .
\end{equation*}
There is substantial MCMC literature establishing that Markov chains are at least polynomially ergodic (see \citet{vats:fleg:jone:2015:spectral} and the references therein).  

We illustrate the generalized IS method and importance of obtaining SEs through three examples.  First, we consider a toy example to demonstrate that BM and RS estimators are consistent and investigate the benefit of allowing general weights to be used in generalized IS.  Second, we consider a Bayesian spatial model for a root rot disease dataset where we illustrate the importance of calculating SEs by considering different designs and performing samples size calculations.  Finally, we consider a standard linear regression model with a large number of variables and use the BM estimator developed here for empirical Bayes variable selection.   

% In addition, a crucial assumption in \pcite{DossTan:2014} and \citet{tan:doss:hobe:2015} for variance estimation via RS is the existence of a minorization condition for each of the $k$ Markov chains.  In practice, each chain needs to regenerate frequently, which is extremely challenging for high-dimensional state spaces. Our work does not require a minorization condition and hence we provide a practical procedure of estimating the standard error of IS estimators.
%weaken the necessary conditions on the underlying Markov chains.

The rest of the paper is organized as follows. Section~\ref{sec:ernc-mcmc} is devoted to the important problem of estimating the ratios of normalizing constants of unnormalized densities, that is estimating $\bd$.  Section~\ref{sec:imp} contains the construction of a CLT for $\hat{\eta}$ and describes how valid SEs of $\hat{\eta}$ can be obtained using BM.  Sections~\ref{sec:toy} and Appendix~\ref{sec:toyapp} contain toy examples illustrating the benefits of different weight functions. Section~\ref{sec:spatial} considers a Bayesian spatial models for binary responses.  The empirical Bayes variable selection example is contained in Appendix~\ref{sec:realdata}.  We conclude with a discussion in Section~\ref{sec:disc}.  All proofs are relegated to the appendices.

\section{Estimating ratios of normalizing constants}
\label{sec:ernc-mcmc}
Consider $k$ densities $\pi_l =
\nu_l/m_l, l=1,\ldots,k$ with respect to the measure $\mu$, where the
$\nu_l$'s are known functions and the $m_l$'s are unknown constants.
For each $l$ we have a positive Harris Markov chain $\Phi_l = \{
X_1^{(l)}, \ldots, X_{n_l}^{(l)} \}$ with invariant density
$\pi_l$. Our objective is to estimate all possible ratios $m_i / m_j,
\, i \neq j$ or, equivalently, the vector $\bd$ defined in \eqref{eq:ddef}.

\pcite{Geyer:1994} reverse logistic regression is described as follows.  Let $n=\sum n_l$ and set $a_l=n_l/n$ for now. For $l = 1, \ldots, k$ define the vector $\bzeta$ by
\begin{equation*}
  \zeta_l = -\log(m_l) + \log(a_l)
\end{equation*}
and let
\begin{equation}
  \label{eq:pl}
  p_l(x, \bzeta) = \frac{ \nu_l(x) e^{\zeta_l} } { \sum_{s=1}^k
  \nu_s(x) e^{\zeta_s} }.
\end{equation}
Given the value $x$ belongs to the pooled sample $\big\{
X_i^{(l)}, \, i = 1, \ldots, n_l, \, l = 1, \ldots, k \big\}$, $p_l(x,
\bzeta)$ is the probability that $x$ came from the $l^{\text{th}}$
distribution. Of course, we know which distribution the sample $x$ came
from, but here we pretend that the only thing we know about $x$ is its
value and estimate $\bzeta$ by
maximizing the log quasi-likelihood function
\begin{equation}
  \label{eq:lql}
  l_n(\bzeta) = \sum_{l=1}^k \sum_{i=1}^{n_l} \log \bigl(
  p_l(X_i^{(l)}, \bzeta) \bigr)
\end{equation}
with respect to $\bzeta$.  Since $\bzeta$ has a one-to-one
correspondence with $\bm = (m_1, \ldots, m_k)$, by estimating $\bzeta$
we can estimate $\bm$. 

As \citet{Geyer:1994} mentioned, there is a
non-identifiability issue regarding $l_n(\bzeta)$: for any constant $c
\in \Real$, $l_n(\bzeta)$ is same as $l_n(\bzeta + c 1_k)$ where
$1_k$ is the vector of $k$ $1$'s. So we can estimate the true $\bzeta$
only up to an additive constant. Thus, we can estimate $\bm$ only up
to an overall multiplicative constant, that is, we can estimate only
$\bd$. Let $\bzeta_0 \in \Real^k$ be defined by $[\bzeta_0]_l =
[\zeta]_l - \bigl( \sum_{s=1}^k [\zeta]_s \bigr)/k$, the true $\bzeta$ normalized to add to
zero. \citet{Geyer:1994} proposed to estimate $\bzeta_0$ by
$\hatbzeta$, the maximizer of $l_n$ subject to the linear constraint
$\bzeta^{\top} 1_k = 0$, and thus obtain an estimate of $\bd$.  The
estimator $\hatbd$ (written explicitly in Section~\ref{sec:batch}),
was introduced by \citet{Vardi:1985}, and studied further by
\citet{GillVardiWellner:1988}, who proved that in the iid setting,
$\hatbd$ is consistent and asymptotically normal, and established its
efficiency. %optimality properties.
 \citet{Geyer:1994} proved the consistency
and asymptotic normality of $\hatbd$ when $\Phi_l, \ldots, \Phi_k$ are $k$ Markov chains satisfying
certain mixing conditions. In the iid setting, \citet{MengWong:1996},
\citet{KongEtal:2003}, and \citet{Tan:2004} rederived the estimate 
under different computational schemes. 

However, none of these articles discuss how to consistently estimate the
covariance matrix of $\hatbd$, even in the iid setting. Recently
\citet{DossTan:2014} address this important issue and obtain a RS estimator of the
covariance matrix of $\hatbd$ in the Markov chain
setting. \citet{DossTan:2014} also mention the optimality results of
\citet{GillVardiWellner:1988} do not hold in the Markov chain
case. In particular, when using Markov chain samples, the choice of
the weights $a_j =n_j/n$ to the probability density $\nu_j/m_j$ in the
denominator of \eqref{eq:pl} is no more optimal and should instead
incorporate the effective sample size of different chains as they
might have quite different rates of mixing. They 
introduce the following more general log quasi-likelihood function
\begin{equation}
  \label{eq:new-lql}
  \ell_n(\bzeta) = \sum_{l=1}^k w_l \sum_{i=1}^{n_l} \log \bigl(
  p_l(X_i^{(l)}, \bzeta) \bigr),
\end{equation}
where the vector $w \in \Real^k$ is defined by $w_l = a_l n / n_l$ for $l = 1, \ldots, k$ for an arbitrary probability vector $\ba$. (Note the change of notation from $l$ to $\ell$.) Clearly if $a_l = n_l/n$, then $w_l = 1$ and~\eqref{eq:new-lql} becomes \eqref{eq:lql}. 

When RS can be used, \citet{DossTan:2014} proved the consistency (to the true value $\bzeta_0$) and asymptotic normality of the constrained maximizer $\hatbzeta$ (subject to the constraint $\bzeta^{\top} 1_k = 0$) of \eqref{eq:new-lql} under geometric ergodicity. They also obtain a RS estimator of the asymptotic covariance matrix and describe an empirical method for choosing the optimal $\ba$ based on minimizing the trace of the estimated covariance matrix of $\hatbd$. However, their procedure requires a practical minorization condition for each of the $k$ Markov chains, which can be extremely difficult. Without a minorization condition, we show $\hatbd$ is a consistent estimator of $\bd$, show $\hatbd$ satisfies a CLT under significantly weaker mixing conditions, and provide a strongly consistent BM estimator of the covariance matrix of $\hatbd$.

\subsection{Central limit theorem and asymptotic covariance estimator}
\label{sec:batch}
Within each Markov chain $l = 1, \ldots, k$, assume $n_l \rightarrow \infty$ in such a way that $n_l /n \rightarrow s_l \in (0, 1)$. In order to obtain the CLT result for $\hatbd$, we first establish a CLT for $\hatbzeta$. Note that the function $g \colon {\Real}^{k} \rightarrow {\Real}^{k-1}$ that maps $\bzeta_0$ into $\bd$ is given by
\begin{equation}
  \label{eq:g}
  g(\bzeta) =
  \begin{pmatrix}
    e^{\zeta_1 - \zeta_2} a_2/a_1 \\
    e^{\zeta_1 - \zeta_3} a_3/a_1 \\
    \vdots                            \\
    e^{\zeta_1 - \zeta_k} a_k/a_1
  \end{pmatrix},
\end{equation}
and its gradient at $\bzeta_0$ (in terms of $\bd$) is
\begin{equation}
  \label{eq:D}
  D =
  \begin{pmatrix}
    d_2    & d_3    & \ldots & d_k    \\
    -d_2   & 0      & \ldots & 0      \\
    0      & -d_3   & \ldots & 0      \\
    \vdots & \vdots & \ddots & \vdots \\
    0      & 0      & \ldots & -d_k
  \end{pmatrix}.
\end{equation}
Since $\bd = g(\bzeta_0)$, and by definition $\hatbd = g(\hatbzeta)$,
we can use the CLT result of $\hatbzeta$ to get a CLT for $\hatbd$.

%In order to state the CLT of $\hatbzeta$, 
First, we introduce the following
notations.  For $r = 1, \dots, k$, let
\begin{equation}
\label{eq:Ys}
Y_i^{(r,l)} = p_r(X_i^{(l)}, \bzeta_0) - E_{\pi_l} \bigl( p_r(X, \bzeta_0) \bigr), \qquad i = 1, \ldots, n_l.
\end{equation}
The asymptotic covariance matrix in the CLT of $\hatbzeta$, involves two
$k \times k$ matrices $B$ and $\Omega$, which we now define. The matrix $B$ is given by
\begin{equation}
  \label{eq:B}
  \begin{split}
    B_{rr} & = \sum_{j=1}^k a_j E_{\pi_j} \bigl( p_r(X,
               \bzeta) [1 - p_r(X, \bzeta)]
               \bigr) \text{ and}\\
    B_{rs} & = - \sum_{j=1}^k a_j E_{\pi_j} \bigl( p_r(X,
               \bzeta) p_s(X, \bzeta) \bigr)
               \text{ for } r \neq s.
  \end{split}
\end{equation}
Let $\Omega$ be the $k \times k$ matrix defined (for $r, s = 1, \ldots, k$) by
\begin{equation}
  \label{eq:Omega}
 \Omega_{rs} =
  \sum_{l=1}^k \frac{a_l^2}{s_l} \Big[  E_{\pi_l}\{Y_1^{(r,l)} Y_1^{(s,l)}\} +  \sum_{i=1}^{\infty} E_{\pi_l}\{Y_1^{(r,l)} Y_{1+i}^{(s,l)}\}+  \sum_{i=1}^{\infty} E_{\pi_l}\{Y_{1+i}^{(r,l)} Y_1^{(s,l)}\}\Big] .
\end{equation}
\begin{remark}
  The right hand side of \eqref{eq:Omega} involves terms
  of the form $E_{\pi_l}\{Y_1^{(r,l)} Y_{1+i}^{(s,l)}\}$ and
  $E_{\pi_l}\{Y_{1+i}^{(r,l)} Y_1^{(s,l)}\}$. For any fixed $l, r, s$
  and $i$, the two expectations are the same if $X_{1}^{(l)}$ and
  $X_{1+i}^{(l)}$ are exchangeable, e.g.\ if the chain $\Phi_l$
  is reversible.  In general, the two
  expectations are not equal.%necessarily the same if $\Phi_l$ is not reversible for some $l$. 
\end{remark}
The matrix $B$ will be estimated by its natural estimate $\widehat{B}$ defined by
\begin{equation}
  \label{eq:Bhat}
  \begin{split}
    \widehat{B}_{rr} & = \sum_{l=1}^k a_l \biggl( \frac{1}{n_l}
                         \sum_{i=1}^{n_l} p_r(X_i^{(l)}, \hatbzeta)
                         \bigl[ 1 - p_r(X_i^{(l)}, \hatbzeta) \bigr]
                         \biggr) \text{ and}\\
    \widehat{B}_{rs} & = - \sum_{l=1}^k a_l \biggl( \frac{1}{n_l}
                         \sum_{i=1}^{n_l} p_r(X_i^{(l)}, \hatbzeta)
                         p_s(X_i^{(l)}, \hatbzeta) \biggr) \text{ for } r \neq s.
  \end{split}
\end{equation}

To obtain a BM estimate $\widehat{\Omega}$, suppose we
simulate the Markov chain $\Phi_l$ for $n_l=e_l b_l$ iterations (hence
$e_l=e_{n_l}$ and $b_l=b_{n_l}$ are functions of $n_l$) and define for
$r, l=1,\ldots,k$
\[
\bar{Z}^{(r,l)}_{m} := \frac{1}{b_l} \sum_{j=m b_l + 1}^{(m+1) b_l} p_r(X_j^{(l)}, \hatbzeta) \hspace*{5mm} \text{ for } m=0,\dots, e_l - 1 \;.
\]
%where $\hat{Y}_{i}^{(r,l)}$ is $Y_{i}^{(r,l)}$ defined in \eqref{eq:Ys} with $\hatbzeta$ substituted for $\bzeta_0$.

Now set $\bar{Z}^{(l)}_{m} = \left( \bar{Z}^{(1,l)}_{m}, \ldots,
  \bar{Z}^{(k,l)}_{m} \right) ^{\top}$ for $m=0,\dots, e_l - 1$.  For $l =1,\ldots, k$, denote
$\bar{\bar{Z}}^{(l)} = \left( \bar{\bar{Z}}^{(1,l)}, \dots, \bar{\bar{Z}}^{(k,l)} \right)^{\top}$ where
$\bar{\bar{Z}}^{(r,l)} =\sum_{i=1}^{n_l} p_r(X_i^{(l)}, \hatbzeta)/n_l$. Let
\begin{equation}
\label{eq:BM}
\widehat{\Sigma}^{(l)} = \frac{b_l}{e_l - 1} \sum_{m=0}^{e_l - 1} \left[ \bar{Z}^{(l)}_{m} - \bar{\bar{Z}}^{(l)} \right] \left[ \bar{Z}^{(l)}_{m} - \bar{\bar{Z}}^{(l)} \right]^T \; \mbox{ for} \;\; l =1,\ldots, k.
\end{equation}
%Note that the above expression of $\widehat{\Sigma}^{(l)}$ does not involve the expectation terms in $\hat{Y}_{i}^{(r,l)}$ and
%$\hat{W}^{(r,l)}$ as they cancel out each other.
Let
\begin{equation}
  \label{eq:Sigmahat}
\widehat{\Sigma} = \begin{pmatrix}
\widehat{\Sigma}^{(1)}\Ddot{11}.(6pt,-2pt,6pt,-1.5pt)&\quad&\quad&\quad 0&\\
&&&&\\
&0&&&\widehat{\Sigma}^{(k)}\\
\end{pmatrix} \;
\end{equation}
and define the following $k \times k^2$ matrix
\begin{equation}
  \label{eq:defA}
A_n = \left( - \sqrt{\frac{n}{n_1}} a_1 I_k \quad - \sqrt{\frac{n}{n_2}} a_2 I_k \quad \dots \quad - \sqrt{\frac{n}{n_k}} a_k I_k \right) \;,
\end{equation}
where $I_k$ denotes the $k \times k$ identity matrix. Finally, define
\begin{equation}
  \label{eq:Omegahat}
\widehat{\Omega} = A_n \widehat{\Sigma} A_n^{\top}.
\end{equation}

We are now ready to describe conditions that ensure strong consistency and asymptotic normality of $\hatbd$. The following theorem
also provides consistent estimate of the asymptotic covariance matrix
of $\hatbd$ using BM method. Consistency of $\hatbd$ holds under minimal assumptions, i.e.\ if $\Phi_1, \ldots,
\Phi_k$ are positive Harris chains. On the other hand, CLTs and consistency of BM
estimator of asymptotic covariance require some mixing conditions on
the Markov chains. For a square matrix $C$, let $C^{\dagger}$ denote the Moore-Penrose inverse
of $C$.
\begin{theorem}
  \label{thm:CLT}
  Suppose that for each $l = 1, \ldots, k$, the Markov chain $\{
  X_1^{(l)}, X_2^{(l)}, \ldots \}$ has invariant distribution
  $\pi_l$.
  \begin{enumerate}
  \item If the Markov chains $\Phi_1, \ldots, \Phi_k$ are positive
    Harris, the log quasi-likelihood function~\eqref{eq:new-lql} has a
    unique maximizer subject to the constraint $\bzeta^{\top} 1_k =
    0$.  Let $\hatbzeta$ denote this maximizer, and let $\hatbd =
    g(\hatbzeta)$.  Then $\hatbd \cas
    \bd$ as $n_1, \ldots,n_k \rightarrow \infty$.
  \item If the Markov chains $\Phi_1, \ldots, \Phi_k$ are
    polynomially ergodic of order $m > 1$, as
    $n_1, \ldots,n_k \rightarrow \infty$, $\sqrt{n} (\hatbd - \bd) \cd {\cal N} (0, V)$ where 
    $V = D^{\top} B^{\dagger} \Omega B^{\dagger}D$.
  \item Assume that the Markov chains $\Phi_1, \ldots, \Phi_k$
    are polynomially ergodic of order $m > 1$ and for all $l=
    1,\ldots,k$, $b_l = \lfloor n_l^{\nu} \rfloor$ where $1 > \nu >
    0$.  Let $\widehat{D}$ be the matrix $D$ in~\eqref{eq:D} with
    $\hatbd$ in place of $\bd$, and let $\widehat{B}$ and
    $\widehat{\Omega}$ be defined by~\eqref{eq:Bhat}
    and~\eqref{eq:Omegahat}, respectively.  Then, $\widehat{V} :=
    \widehat{D}^{\top} \widehat{B}^{\dagger} \widehat{\Omega}
    \widehat{B}^{\dagger} \widehat{D}$ is a strongly consistent
    estimator of $V$.
  \end{enumerate}
\end{theorem}

% Theorem~\ref{thm:CLT} does not require the chains to be stationary so there is no need for burn-in.  Further, spectral methods can also be used to obtain a strongly consistent estimator of $V$ using results in \citet{vats:fleg:jone:2015}.

\section{IS with multiple Markov chains}
\label{sec:imp}
This section considers a CLT and SEs for
the generalized IS estimator $\hat{\eta}$.  From \eqref{eq:mainest}
we see that $\hat{\eta} \equiv \hat{\eta}^{[f]} (\pi; \ba, \bd) =
\hat{v}^{[f]} (\pi, \pi_1; \ba, \bd) / \hat{u}(\pi, \pi_1; \ba, \bd)$, where
\begin{equation}
  \label{eq:uvhat}
  \begin{split}
\hat{v} &\equiv \hat{v}^{[f]} (\pi, \pi_1; \ba, \bd) := \sum_{l=1}^k \frac{a_l} {n_l} \sum_{i=1}^{n_l} v^{[f]}(X_i^{(l)}; \ba, \bd) \mbox{ and} \\
  \hat{u} & \equiv \hat{u}(\pi, \pi_1; \ba, \bd) := \sum_{l=1}^k \frac{a_l} {n_l} \sum_{i=1}^{n_l}
  u(X_i^{(l)}; \ba, \bd)
\end{split}
\end{equation}
with
\begin{equation}
  \label{eq:uv}
  v^{[f]}(x; \ba, \bd) := f(x) u(x; \ba, \bd)
  \quad \text{and} \quad u(x; \ba, \bd) := \frac{\nu(x)} {\sum_{s=1}^k a_s \nu_s(x) / d_s}.
\end{equation}
% $d_i = m_{h_i}/m_{h_1}, i=2, \dots, k, \bd = (d_2, \dots, d_k)$ for
% $k$ skeleton points $h_1, \ldot, h_k$.
Note that $\hat{u}$ converges almost surely to
\begin{equation}
\label{eq:nucon}
\sum_{l=1}^k a_l E_{\pi_l} u(X; \ba, \bd) =  \int_{\sX} \frac{\sum_{l=1}^k a_l
\nu_l(x) / m_l} {\sum_{s=1}^k a_s \nu_s(x) / (m_s /
m_1)} \nu(x) \, \mu(dx) = \frac{m}{m_1},
\end{equation}
as $n_1, \ldots, n_k \rightarrow \infty$. Thus $\hat{u}$ itself is a useful quantity as it consistently estimates the ratios of normalizing constants $\{u(\pi, \pi_1) \equiv m/ m_1 | \pi \in \Pi\}$. Unlike the estimator $\hatbd$ in Section~\ref{sec:ernc-mcmc}, $\hat{u}$ does not require a sample from each density $\pi \in \Pi$. Thus $\hat{u}$ is well suited for situations where one wants to estimate the ratios $u(\pi, \pi_1)$ for a very large number of $\pi$'s based on samples from a small number of skeleton densities, say $k$.  Thus this method is particularly efficient when obtaining samples from the target distributions is computationally demanding and the distributions within $\Pi$ are similar.

In the context of Bayesian analysis, let $\pi(x) = \text{lik}(x)
p(x)/m$ be the posterior density corresponding to the likelihood
function $\text{lik}(x)$ and prior $p(x)$ with normalizing constant
$m$. In this case, $u(\pi, \pi_1)$ is the so-called Bayes factor
between the two models, which is commonly used in model
selection.

The estimators $\hu$ and $\hat{v}$ in \eqref{eq:uvhat} depend on
$\bd$, which is generally unknown in practice. Here we consider a two-stage procedure for evaluating
$\hu$. In the 1st stage, $\bd$ is estimated by its reverse logistic
regression estimator $\hatbd$ described in Section~\ref{sec:ernc-mcmc}
using Markov chains $\tilde{\Phi}_l \equiv
\{\tilde{X_i^l}\}_{i=1}^{N_l}$ with stationary densities $\pi_l$, for
$l=1,\ldots,k$. Note the change of notation from
Section~\ref{sec:ernc-mcmc} where we used $n_l$'s to denote the length
of the Markov chains. In order to avoid more notations, we use
$\tilde{\Phi}_l$'s and $N_l$'s to denote the stage 1 chains and their
length respectively. Once $\hatbd$ is formed, new MCMC samples
$\Phi_l \equiv \{X_i^l\}_{i=1}^{n_l}, l=1\ldots,k$ are obtained and
$u(\pi, \pi_1) (E_\pi f)$ is estimated using $\hat{u}(\pi, \pi_1; \ba,
\hatbd)$ ($\hat{\eta}^{[f]} (\pi; \ba, \hatbd)$) based on these 2nd
stage samples. \citet{buta:doss:2011} propose this two-stage method and quantify its benefits over the method where the same MCMC samples are used to estimate both $\bd$ and $u(\pi, \pi_1)$. 

\subsection{Estimating ratios of normalizing constants}
\label{sec:elrnc}
Before we state a CLT for $\hat{u}(\pi, \pi_1; \ba, \hatbd)$, we require some notation.  Let 
\begin{equation}
  \label{eq:taul}
  \tau^2_l(\pi ; \ba, \bd) = \mbox{Var}_{\pi_l} (u(X_1^{(l)}; \ba, \bd)) + 2 \sum_{g=1}^{\infty}
\mbox{Cov}_{\pi_l} (u(X_1^{(l)}; \ba, \bd), u(X_{1+g}^{(l)}; \ba, \bd))
\end{equation}
and $\tau^2 (\pi ; \ba, \bd) = \sum_{l =1}^k (a_l^2/s_l) \tau^2_l(\pi ;
\ba, \bd)$.  Further, define $c(\pi; \ba, \bd)$ as a vector of length $k-1$ with $(j-1)$th coordinate as
\begin{equation}
  \label{eq:cdef}
   [c(\pi; \ba, \bd)]_{j-1} = \frac{u(\pi, \pi_1)}{ d_j^2} \int_{\sX} \frac{a_j
  \nu_{j} (x)}{\sum_{s = 1}^k a_s \nu_{s} (x)/d_s} \pi (x) dx \mbox{ for } j=2,\dots,k
\end{equation}
and $\hat{c}(\pi; \ba, \bd)$ as a vector of length $k-1$ with $(j-1)$th coordinate as
\begin{equation}
  \label{eq:chatdef}
   [\hat{c}(\pi; \ba, \bd)]_{j-1} = \sum_{l=1}^k \frac{1}{n_l}
\sum_{i=1}^{n_l} \frac{a_j a_l \nu (X_i^{(l)}) \nu_{j}
  (X_i^{(l)})}{ (\sum_{s = 1}^k a_s \nu_{s} (X_i^{(l)})/d_s)^2 d_j^2} \mbox{ for } j=2,\dots,k.
\end{equation}
Assuming $n_l=e_lb_l$, then let
\begin{equation}
  \label{eq:tauldef}
\hat{\tau}^2_l(\pi ; \ba, \bd)=\frac{b_l}{e_l - 1} \sum_{m=0}^{e_l
  - 1} \left[ \bar{u}_{m}(\ba, \bd) - \bar{\bar{u}}(\ba,
  \bd)\right]^2,
\end{equation}
 where $\bar{u}_{m}(\ba, \bd)$ is the average
of the $(m+1)$st block $\{u(X_{mb_l+1}^{(l)}; \ba, \bd),\cdots,
u(X_{(m+1)b_l}^{(l)}; \ba, \bd)\}$, and $\bar{\bar{u}}(\ba,
\bd)$ is the overall average of $\{u(X_{1}^{(l)}; \ba, \bd),
\cdots, u(X_{n_l}^{(l)}; \ba, \bd)\}$.  Here, $b_l$ and $e_l$ are
the block sizes and the number of blocks respectively. Finally let
$\hat{\tau}^2 (\pi ; \ba, \bd) = \sum_{l =1}^k (a_l^2/s_l)
\hat{\tau}^2_l(\pi ; \ba, \bd)$.

\begin{theorem}
  \label{thm:elnc}
  Suppose that for the stage 1 chains, conditions of
  Theorem~\ref{thm:CLT} holds such that $N^{1/2} (\hatbd - \bd) \cd
  {\cal N} (0, V)$ as $N\equiv \sum_{l=1}^k N_l \rightarrow
  \infty$. Suppose there exists $q \in [0, \infty)$
such that $n/N \rightarrow q$ where $n=\sum_{l=1}^k n_l$ is the total
sample size for stage 2. In addition, let $n_l/n \rightarrow s_l$ for
$l = 1,\cdots,k$.
\begin{enumerate}
  \item Assume that the stage 2 Markov chains $\Phi_1, \ldots,
  \Phi_k$ are polynomially ergodic of order $m$, and for some $\delta >0$ $E_{\pi_l} |u(X; \ba, \bd)|^{2+\delta} < \infty$ for each $l = 1,\cdots,k$ where $m > 1+2/\delta$. Then as $n_1, \ldots,n_k \rightarrow \infty$,
  \begin{equation}
    \label{eq:elncclt}
    \sqrt{n} (\hat{u}(\pi, \pi_1; \ba, \hatbd) - u(\pi, \pi_1)) \cd N(0, q
c(\pi; \ba, \bd)^{\top} V c(\pi; \ba, \bd) + \tau^2 (\pi ; \ba, \bd)) .
  \end{equation}
\item Let $\widehat{V}$ be the consistent estimator of $V$ given
  in Theorem~\ref{thm:CLT} (3). Assume that the Markov chains $\Phi_1,
  \ldots, \Phi_k$ are polynomially ergodic of order $m \ge
    (1+\epsilon)(1 + 2/\delta)$ for some $\epsilon, \delta >0$ such
    that $E_{\pi_l} |u(X; \ba, \bd)|^{4+\delta} < \infty$, and for all $l=
    1,\ldots,k$, $b_l = \lfloor n_l^{\nu} \rfloor$ where $1 > \nu >
    0$. Then $q \hat{c}(\pi; \ba, \hatbd)^{\top}
    \widehat{V} \hat{c}(\pi; \ba, \hatbd) + \hat{\tau}^2 (\pi ; \ba,
    \hatbd))$ is a strongly consistent estimator of the asymptotic
    variance in \eqref{eq:elncclt}.
  \end{enumerate}
\end{theorem}

% The proof of Theorem~\ref{thm:elnc} is given in Appendix~\ref{sec:appthm2}.
Note that the asymptotic variance in \eqref{eq:elncclt} has two
components. The second term is the variance of $\hu$ when $\bd$ is
known. The first term is the increase in the variance of $\hu$ resulting
from using $\hatbd$ instead of $\bd$. Since we are interested in
estimating $u(\pi, \pi_1)$ for a large number of $\pi$'s and for every
$\pi$ the computational time needed to calculate $\hu$ in
\eqref{eq:uvhat} is linear in the total sample size $n$, this can not
be very large. If generating MCMC samples is not computationally
demanding, then long chains can be used in the 1st stage (that is,
large $N_l$'s can be used) to obtain a precise estimate of $\bd$, and thus greatly
reducing the first term in the variance expression \eqref{eq:elncclt}.

\subsection{Estimation of expectations using generalized IS}
\label{sec:estexp}
This section discusses estimating SEs of the generalized IS estimator $\hat{\eta}$ given in \eqref{eq:mainest}.  In order to state a CLT for $\hat{\eta}$ we define the following notations:
\[
\gamma^{11}_l \equiv \gamma^{11}_l(\pi; \ba, \bd) = \mbox{Var}_{\pi_l} (v^{[f]}(X_1^{(l)}; \ba, \bd)) + 2 \sum_{g=1}^{\infty}
\mbox{Cov}_{\pi_l} (v^{[f]}(X_1^{(l)}; \ba, \bd), v^{[f]}(X_{1+g}^{(l)}; \ba, \bd)),
\]
\[
\begin{split}
  \gamma^{12}_l & \equiv \gamma^{12}_l(\pi; \ba, \bd) = \gamma^{21}_l
  \equiv \gamma^{21}_l(\pi; \ba, \bd)= \mbox{Cov}_{\pi_l}
  (v^{[f]}(X_1^{(l)}; \ba, \bd), u(X_1^{(l)}; \ba, \bd)) \\ &+
  \sum_{g=1}^{\infty} [\mbox{Cov}_{\pi_l} (v^{[f]}(X_1^{(l)}; \ba,
  \bd), u(X_{1+g}^{(l)}; \ba, \bd)) +\mbox{Cov}_{\pi_l}
  (v^{[f]}(X_{1+g}^{(l)}; \ba, \bd), u(X_1^{(l)}; \ba, \bd))],
\end{split}
\]
\[
\gamma^{22}_l \equiv\gamma^{22}_l(\pi; \ba, \bd) = \mbox{Var}_{\pi_l} (u(X_1^{(l)}; \ba, \bd)) + 2 \sum_{g=1}^{\infty}
\mbox{Cov}_{\pi_l} (u(X_1^{(l)}; \ba, \bd), u(X_{1+g}^{(l)}; \ba, \bd)),
\]
(note $\gamma^{22}_l$ is the same as $\tau^2_l(\pi; \ba, \bd)$ defined in \eqref{eq:taul}) and
\begin{equation}
  \label{eq:defgam}
  \Gamma_l(\pi; \ba, \bd) = \left(
  \begin{array}{lr}
    \gamma^{11} & \gamma^{12}\\
    \gamma^{21} & \gamma^{22}\\
  \end{array}
\right);
\Gamma (\pi; \ba, \bd) = \sum_{l =1}^k \frac{a_l^2}{s_l} \Gamma_l(\pi; \ba, \bd).
\end{equation}
Since $\hat{\eta}$ has the form of a ratio, to establish a CLT for it, we apply the
Delta method on the function $h(x, y) = x/y, \;\mbox{with}\; \nabla h(x, y) = (1/y, -x/y^2)' $.
Let
\begin{equation}
  \label{eq:defrho}
  \rho(\pi; \ba, \bd) = \nabla h(E_\pi f u(\pi, \pi_1), u(\pi, \pi_1))' \Gamma(\pi; \ba, \bd) \nabla h( E_\pi f u(\pi, \pi_1), u(\pi, \pi_1)),
\end{equation}
$e(\pi; \ba, \bd)$ is a vector of length $k-1$ with $(j-1)$th coordinate as
\begin{equation}
  \label{eq:defe}
 [e(\pi; \ba, \bd)]_{j-1} = \frac{a_j}{d_j^2} \int_{\sX} \frac{[f(x) - E_\pi f]
  \nu_{j} (x)}{\sum_{s = 1}^k a_s \nu_{s} (x)/d_s} \pi(x) dx ,\;\; j=2,\dots,k,
\end{equation}
and $\hat{e}(\pi; \ba, \bd)$ is a vector of length $k-1$ with $(j-1)$th coordinate as
\begin{align}
  \label{eq:defhate}
  [\hat{e}(\pi; \ba, \bd)]_{j-1}  \equiv \frac{\sum_{l=1}^k \frac{a_l}{n_l}
    \sum_{i=1}^{n_l} \frac{a_j f(X_i^{(l)}) \nu
      (X_i^{(l)}) \nu_{j} (X_i^{(l)})}{d_j^2 (\sum_{s =
        1}^k a_s \nu_{s} (X_i^{(l)})/d_s)^2}}{\hat{u} (\pi, \pi_1; \ba, \bd)} - \frac{[c(\pi; \ba, \bd)]_{j-1} \hat{\eta}^{[f]} (\pi; \ba, \bd)}{\hat{u} (\pi, \pi_1; \ba, \bd)},
\end{align}
where $[c(\pi; \ba, \bd)]_{j-1}$ is defined in \eqref{eq:chatdef}.
% We use the same two-stage procedure, as in Section~\ref{sec:elrnc}, for evaluating $\hat{\eta}$. Again
Assuming $n_l=e_lb_l$, let
\begin{align*}
  \widehat{\Gamma}_{l} (\pi ; \ba, \bd) &= \frac{b_l}{e_l - 1} \sum_{m=0}^{e_l - 1} \Bigg[ \left(
  \begin{array}{c}
    \bar{v}_{m}^{[f]}\\
    \bar{u}_{m}\\
  \end{array}
\right) -  \left(
  \begin{array}{c}
    \bar{\bar{v}}^{[f]}\\
    \bar{\bar{u}}\\
  \end{array}
\right) \Bigg] \Bigg[ \left(
  \begin{array}{c}
    \bar{v}_{m}\\
    \bar{u}_{m}\\
  \end{array}
\right) -  \left(
  \begin{array}{c}
    \bar{\bar{v}}^{[f]}\\
    \bar{\bar{u}}\\
  \end{array}
\right) \Bigg]^{\top} \\
&=
\frac{b_l}{e_l - 1} \left(
  \begin{array}{cc}
    \sum_{m=0}^{e_l - 1} \left[ \bar{v}^{[f]}_{m} - \bar{\bar{v}}^{[f]}\right]^2& \sum_{m=0}^{e_l - 1} \left[ \bar{v}^{[f]}_{m} - \bar{\bar{v}}^{[f]}\right]\left[ \bar{u}_{m} - \bar{\bar{u}}\right] \\
\sum_{m=0}^{e_l - 1} \left[ \bar{v}^{[f]}_{m} - \bar{\bar{v}}^{[f]}\right]\left[ \bar{u}_{m} - \bar{\bar{u}}\right]  &\sum_{m=0}^{e_l - 1}\left[ \bar{u}_{m} - \bar{\bar{u}}\right]^2\\
  \end{array}
\right)\\
&=
\left(\begin{array}{cc}
    \hat{\gamma}^{11}(\pi ; \ba, \bd) & \hat{\gamma}^{12}(\pi ; \ba, \bd)\\
    \hat{\gamma}^{21}(\pi ; \ba, \bd) & \hat{\gamma}^{22}(\pi ; \ba, \bd)\\
\end{array}
\right),
\end{align*}
 where $\bar{v}^{[f]}_{m}$ is the average
of the $(m+1)$st block $\{v^{[f]}(X_{mb_l+1}^{(l)}; \ba, \bd),\cdots,
v^{[f]}(X_{(m+1)b_l}^{(l)}; \ba, \bd)\}$, $\bar{\bar{v}}^{[f]}$ is the overall average of $\{v^{[f]}(X_{1}^{(l)}; \ba, \bd),
\cdots, v^{[f]}(X_{n_l}^{(l)}; \ba, \bd)\}$ and $\bar{u}_{m} \equiv \bar{u}_{m} (\pi, \ba, \bd), \bar{\bar{u}} \equiv \bar{\bar{u}} (\pi, \ba, \bd)$ defined in Section~\ref{sec:elrnc}.
Finally let
$\widehat{\Gamma} (\pi ; \ba, \bd) = \sum_{l =1}^k (a_l^2/s_l)
\widehat{\Gamma}_l(\pi ; \ba, \bd)$, and
\[\hat{\rho}(\pi; \ba, \hatbd) = \nabla h(\hat{v}^{[f]}(\hatbd), \hat{u}(\hatbd))' \widehat{\Gamma}(\pi; \ba, \hatbd) \nabla h(\hat{v}^{[f]}(\hatbd), \hat{u}(\hatbd)).
\]
\begin{theorem}
  \label{thm:elex}
  Suppose that for the stage 1 chains, conditions of
  Theorem~\ref{thm:CLT} hold such that $N^{1/2} (\hatbd - \bd) \cd
  {\cal N} (0, V)$ as $N\equiv \sum_{l=1}^k N_l \rightarrow
  \infty$. Suppose there exists $q \in [0, \infty)$
such that $n/N \rightarrow q$ where $n=\sum_{l=1}^k n_l$ is the total
sample size for stage 2. In addition, let $n_l/n \rightarrow s_l$ for
$l = 1,\cdots,k$.
\begin{enumerate}
  \item Assume that the stage 2 Markov chains $\Phi_1, \ldots,
  \Phi_k$ are polynomially ergodic of order $m$, and for some $\delta >0$ $E_{\pi_l} |u(X; \ba, \bd)|^{2+ \delta} < \infty$ and $E_{\pi_l} |v^{[f]}(X; \ba, \bd)|^{2+ \delta} < \infty$, for each $l = 1,\cdots,k$ where $m > 1+2/\delta$. Then as $n_1, \ldots,n_k \rightarrow \infty$,
  \begin{equation}
    \label{eq:elexclt}
    \sqrt{n} (\hat{\eta}^{[f]}(\pi; \ba, \hatbd) - E_\pi f) \cd N(0, q
e(\pi; \ba, \bd)^{\top} V e(\pi; \ba, \bd) + \rho (\pi ; \ba, \bd)) .
  \end{equation}
\item Let $\widehat{V}$ be the consistent estimator of $V$ given
  in Theorem~\ref{thm:CLT} (3). Assume that the Markov chains
    $\Phi_1, \ldots, \Phi_k$ are polynomially ergodic of order $m \ge
    (1+\epsilon)(1 + 2/\delta)$ for some $\epsilon, \delta >0$ such
    that $E_{\pi_l} |u(X; \ba, \bd)|^{4+\delta} < \infty, E_{\pi_l}
    |v^{[f]}(X; \ba, \bd)|^{4+ \delta} < \infty$, and for each $l =
    1,\cdots,k$, $b_l = \lfloor n_l^{\nu} \rfloor$ where
  $1 > \nu > 0$. Then $q \hat{e}(\pi; \ba, \hatbd)^{\top} \widehat{V}
  \hat{e}(\pi; \ba, \hatbd) + \hat{\rho} (\pi ; \ba, \hatbd))$ is a
  strongly consistent estimator of the asymptotic variance in
  \eqref{eq:elexclt}.
  \end{enumerate}
\end{theorem}

%The proof of Theorem~\ref{thm:elex} is given in Appendix~\ref{sec:appthm3}.
\begin{remark}
  \label{rem:genBD}
Part (1) of Theorems~\ref{thm:elnc} and~\ref{thm:elex} extend \pcite{buta:doss:2011} Theorems~1 and~3, respectively.  Specifically, they require $a_l = n_l/n$ which is a non-optimal choice for $\ba$ \citep{tan:doss:hobe:2015}.  Our results also substantially weaken the Markov chain mixing conditions.
\end{remark}

\begin{remark}
  \label{rem:tdhdiff}
Theorems~\ref{thm:elnc} and~\ref{thm:elex} prove consistency of the BM estimators of the variances of $\hu$ and $\hat{\eta}$ for a general $\ba$.  This extends results in \citet{tan:doss:hobe:2015}, which provides RS based estimators of the asymptotic variance of $\hu$ and $\hat{\eta}$ in the special case when $\ba = (1, \hatbd)$. With this particular choice, $u(x; \ba, \hatbd)$ and $v^{[f]}(x; \ba, \hatbd)$ in \eqref{eq:uv} become free of $\hatbd$ leading to independence among certain quantities.  However, one can set $\ba = w*(1, \hatbd)$ for any user specified fixed vector $w$, which allows the expression in (2.18) of \citet{tan:doss:hobe:2015} to be free of $\hatbd$ and thus the necessary independence.  Hence, their RS estimator can also be applied to an arbitrary vector $\ba$ (details are given in Appendix~\ref{sec:regen}).%the supplement).
\end{remark}

\begin{remark}
A sufficient condition for the moment assumptions for $u$ in Theorems~\ref{thm:elnc} and \ref{thm:elex} is that, for any $\pi \in \Pi$, $\sup_x \left\{\pi(x)\big/\sum_{s=1}^k a_s\pi_s(x)\right\}<\infty$. That is, in any given direction, the tail of at least one of $\{\pi_s, s =1, \dots, k\}$ is heavier than that of $\pi$. This is not hard to achieve in practice by properly choosing $\{\pi_s\}$ with regard to $\Pi$ \cite[see e.g.][]{roy:2014}. Further, if $\E_\pi |f|^{4+\delta}<\infty$ then the moment assumptions for $v^{[f]}$ are satisfied. %The $(4+\delta)$th moment condition is a typical requirement for using the BM method in single Markov chain settings.}
 \end{remark}
% \begin{remark}
% Spectral methods can also be used in Part (2) of Theorems~\ref{thm:elnc} and~\ref{thm:elex} to obtain conditions that ensure strong consistency using results in \citet{vats:fleg:jone:2015}.
% \end{remark}

\section{Toy example}
\label{sec:toy}

In this section, we employee a toy example to confirm that both the BM and
the RS estimators are consistent, 
%, and to leave the readers the impression that despite the BM method is much easier to use, it may experience more variability than the RS method, provided the Markov chains regenerate often enough. A second purpose is to demonstrate
as well as demonstrate the benefit of allowing general weights to be used in the generalized IS estimator. 
Let $t_{r,\mu}$ denote the t-distribution with degree of freedom $r$ and central parameter $\mu$. We set $\pi_1(\cdot)=\nu_1(\cdot)$ and $\pi_2(\cdot)=\nu_2(\cdot)$, which are the density functions for a $t_{5, \mu_1=1}$ and $t_{5, \mu_2=0}$, respectively. Pretending that we do not know the value of the ratio between the two normalizing constants, $d=m_2/m_1=1/1$, we estimate it by the stage~1 estimator 
%general reverse logistic regression estimator 
$\hat{d}$ from section~\ref{sec:ernc-mcmc}, 
%, which maximizes the log quasi-likelihood function in (2.3) for a given weight $\ba$. 
and compare the BM and the RS method in estimating the asymptotic variance. As for the stage~2 estimators from section~\ref{sec:imp}, the choice of weight and performance of the BM and the RS methods in assessing estimators' uncertainty are studied in Appendix~\ref{sec:toyapp}. %the supplement.
%We will also study a sea of $t$-distributions $\Pi=\{t_{5,\mu}: \mu \in M\}$ where $M$ is a fine grid over $[0,1]$, say $M=\{0,.01,\cdots,.99, 1\}$. For each $\mu \in M$, the goal is to estimate $\E_{t_{5,\mu}}X$ and the ratio between its normalizing constant and $m_1$.  Full details can be found in %Appendix~E.%
%the supplement.

We draw iid samples from $\pi_1$ and Markov chain samples from $\pi_2$ using the independent Metropolis Hastings algorithm with proposal density $t_{5,1}$. It is simple to show $\inf_x\frac{t_{5,\mu}(x)}{t_{5,0}(x)}> 0$, which implies the algorithm is uniformly ergodic (\citet{mengersen1996rates} Theorem 2.1) and hence polynomially ergodic and geometrically ergodic. For RS, our carefully tuned minorization condition enables the Markov chain for $\pi_2$ to regenerate about every 3 iterations. In contrast, the BM method proposed here requires no such theoretical development. 

\begin{figure}[htb]
  \begin{center}
    \includegraphics[width=.45\linewidth]{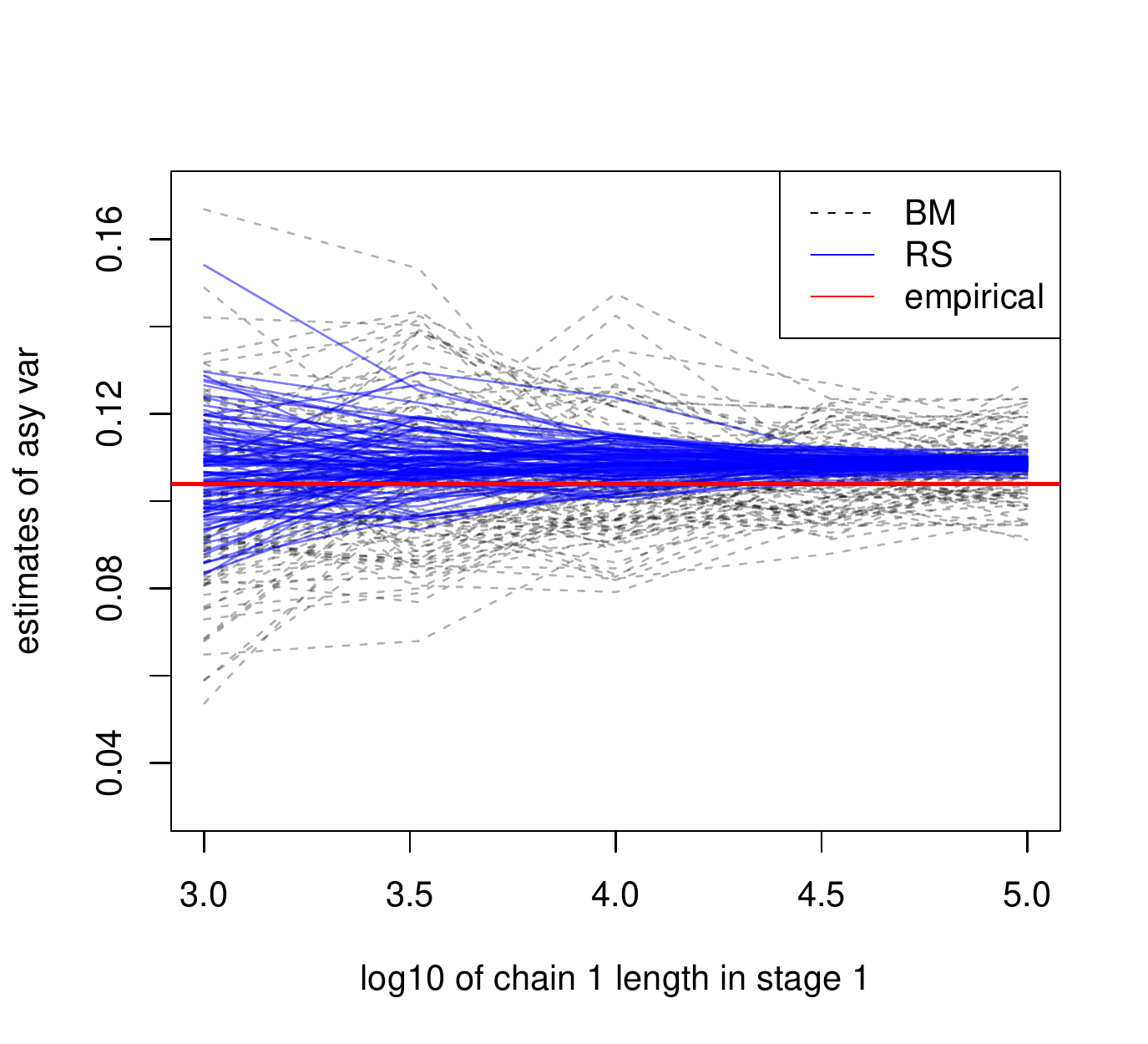}
        \includegraphics[width=.45\linewidth]{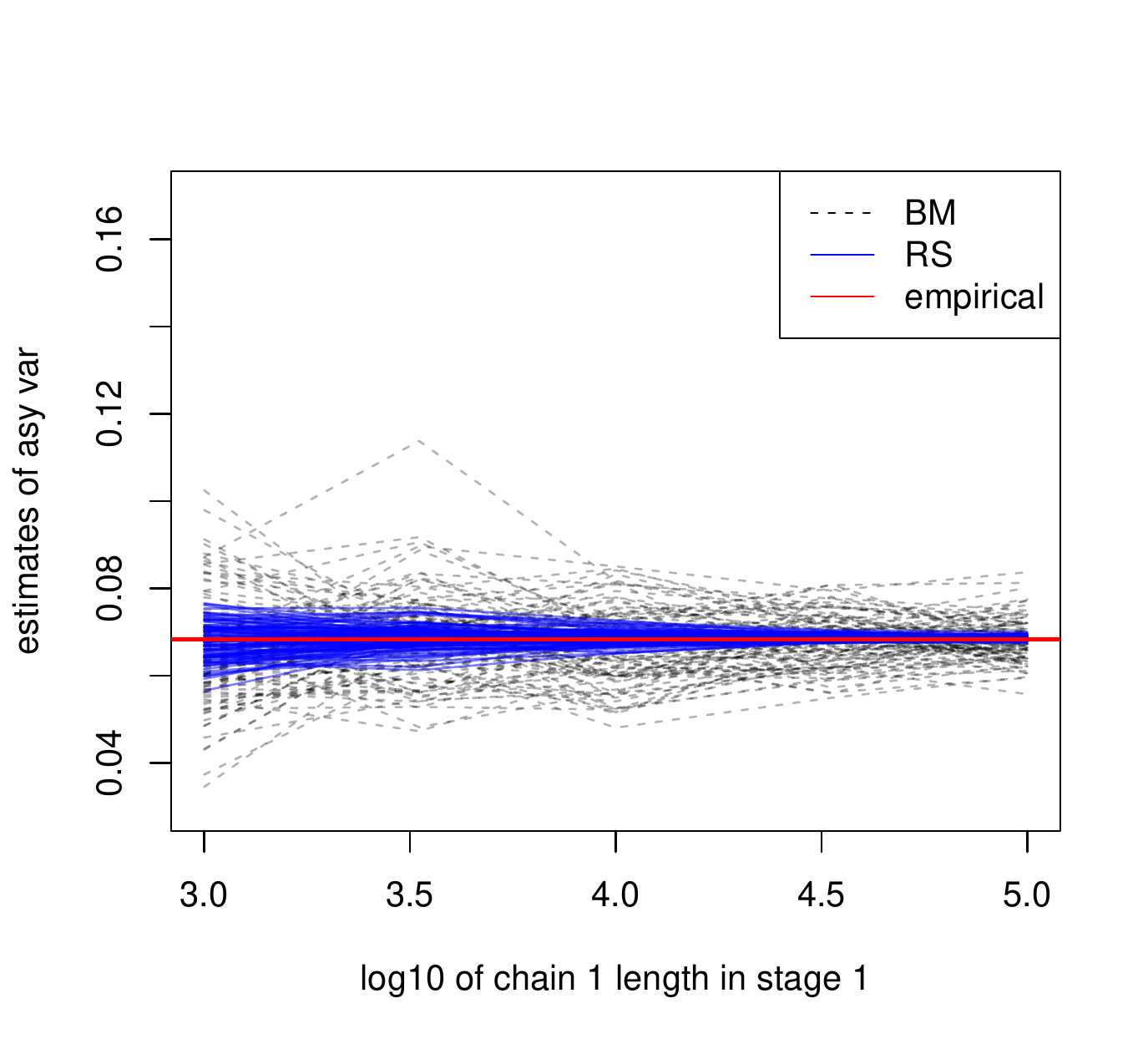}
  \end{center}
       \caption{Plots of BM and RS estimates of the asymptotic variance of $\hat{d}$ in stage~1 for 100 randomly chosen replications. The left panel is based on the naive weight, $\baone=(0.5, 0.5)$ and the right panel is based on a close-to-optimal weight, $\baone=(0.82, 0.18)$. Horizontal lines represent the empirical asymptotic variance of $\hat{d}$ obtained over all replications.}
\label{fig:a1}
\end{figure}

We evaluated the variance estimators at various sample sizes with different choices of weight. Figure~\ref{fig:a1} displays traces of the BM and the RS estimates of the asymptotic variance of $\hat{d}$, %in stage~1, 
in dashed and solid lines, respectively. Overall, the BM and the RS estimates approach the empirical asymptotic variance as the sample size increases, suggesting their consistency. Due to the frequency of regenerations, BM estimates are generally more variable than RS estimates. Further, the left panel of Figure~\ref{fig:a1} is for estimators based on the naive weight, $\ba=(0.5, 0.5)$, that is proportional to the sample sizes; and the right panel is for estimators based on $\ba=(0.82, 0.18)$, that emphasizes the iid sample more than the Markov chain sample. Indeed, the latter weight is a close-to-optimal weight obtained with a small pilot study (see Appendix~\ref{sec:toyapp} for details). Using such a method to choose weight can lead to big improvement in the efficiency of $\hat{d}$ if the mixing rate of the multiple samples differ a lot.

%Similarly for stage~2 we observe convergence of the BM and the RS estimates to the corresponding empirical asymptotic variances of $\hat{d}_\mu$ and $\hat{\E}_\mu(X)$.

\section{Bayesian spatial models for binary responses}
\label{sec:spatial}

%\blue{[Aixin's comment: There is no mention of the purpose of including this example in this paper of multiple IS  until two paragraphs later.  That is, there is no mention of "Bayesian sensitivity analysis". This might lead the reader to lose track of the connection between this section and the previous ones, and on a separate matter, make them to look forward to a disease map more eagerly. Vivek, if you think it's appropriate to bring back some of the commented off descriptions of Bayesian sensitivity analysis in the tex file here, please do so.]}

In this section, we analyze a root rot disease dataset collected from a 90-acre farm in the state of Washington \citep{zhan:2002}. All computations are done in R, using the package \text{geoBayes} \citep{evan:roy:manual:2015}. Recorded at $M=100$ chosen sites are the longitude and the latitude $s_i$, the root counts $\ell(s_i)$, and the number of infected roots $y(s_i), i=1,\ldots, \m$.  Of interest is a map of the disease rate over the entire area for precision farming. We consider the following spatial generalized linear mixed model (SGLMM), similar to that used by \citet{zhan:2002} and \citet{roy:evan:zhu:2016}. %To account for smooth changing of the disease rate and other spatial correlations, w
Taking $\ell(s_i)$ and $s_i$ as fixed, let
\[y(s_i) | z(s_i) \ind \text{Binomial}\left(\ell(s_i), \Phi\left(z(s_i)\right)\right), i=1,\ldots, \m.\]
Here $\bz = (z(s_1), \dots, z(s_{\m}))$ is a vector of latent variables, which is assumed to be a subvector of
a Gaussian random field (GRF) $\{z_s, s\in S\}$, that has a constant mean $\mu$, and a covariance function 
\[\Cov\left(z(s), z(s')\right) = \sigma^2 \rho_\phi (\|s-s'\|) + \omega\, \sigma^2 I_s(s')\,.\] 
Here, $\sigma^2$ is the partial sill, $\|\cdot\|$ denotes the Euclidean distance, and
$\rho_\phi$ is a correlation function from the spherical family with range parameter $\phi$. That is, $\rho_\phi(u) = 1 - \frac{3}{2}\frac{u}{\phi} + \frac{1}{2} \left(\frac{u}{\phi}\right)^3$ for $u\in (0, \phi)$. Next, $I_s(s')$ is an indicator that equals $1$ if $s=s'$, and equals $0$ otherwise. Finally, $\omega \sigma^2$ is the nugget effect, accounting for any remaining variability at site $s$ such as measurement error, while $\omega \in \mathcal{R}^+$ is the relative size of the nugget to the partial sill.  
%\[L(\beta, \sigma^2, \theta, \nu | \by) = \int_{\mathcal{R}^n} \bigg[\prod_{i=1}^n p(y_i |z_i, \nu)\bigg] f(\bz | \beta, \sigma^2, \theta) d \bz .\]
Following \cite{roy:evan:zhu:2016} we assign a non-informative Normal-inverse-Gamma prior to $(\mu, \sigma^2)$ which is (conditionally) conjugate for the model. Specifically, 
\[ \mu|\sigma^2 \sim \text{N}(0, 100\sigma^2), \;\;\text{and} \;\;
f(\sigma^2) \propto \left(\sigma^2\right)^{-\frac{1}{2}-1} \exp\left(-\frac{1}{2\sigma^2}\right)\,.\]
%\red{Are there reference to this? Aixin got these numbers from the R code examples of roy:evan:zhu:2016, that set $V=100, n_\sigma=1, a_\sigma=1$ in
%\[ \mu|\sigma^2 \sim \text{N}(0, V\sigma^2), \;\;\text{and} \;\;
%f(\sigma^2) \propto \left(\sigma^2\right)^{-\frac{n_\sigma}{2}-1} \exp\left(-\frac{n_\sigma a_\sigma}{2\sigma^2}\right)\,.\]
%}
%(Christensen and Waagepetersen 2002, Biometrics). 
Assigning priors for $h = (\phi, \omega)$ in the correlation function of the Gaussian random field is usually
difficult, and the choice of prior may influence the inference 
\citep{chri:2004}. %p. 716). 
Hence we perform a sensitivity analysis focusing on obtaining the Bayes factor (BF) of the model at $h$ relative to a baseline $h_0$ %$m_h(\by)/m_{h_0}(\by)$, 
for a range of values $h \in {\cal H}$. 
Note that for a fixed $h = (\phi, \omega)$, the Bayesian model described above has parameters $\psi = (\mu, \sigma^2)$. Conditioning on the observed data $\by = (y(s_1), \dots, y(s_M))$, inference is based on the posterior density
\begin{equation}\label{eq:pih}
\pi_{h}(\psi | \by) = \frac{L_h(\psi | \by) \pi(\psi)}{{ m_h}(\by)},
\end{equation}
where $L_h(\psi | \by) = \int_{\mathcal{R}^{\m}} f(\by | \bz) f_h(\bz | \psi) d \bz$ is the likelihood, $\pi(\psi)$ is the
prior on $\psi$, and $m_h(\by) = \int_{\mathcal{R} \times \mathcal{R}_+}  L_h(\psi | \by) \pi(\psi) d\psi$ is the normalizing constant, also called the marginal likelihood. The Bayes factor between any two models indexed by $h$ and $h_0$ is $m_h(\by)/m_{h_0}(\by)$, and
the empirical Bayes choice of $h$ is $\underset{h \in \mathcal{H}}{\arg\max}\, {m_h}(\by)=\underset{h \in \mathcal{H}}{\arg\max}\, \left[{m_h}(\by)/m_{h_0}(\by) \right]$.  Our plan is to get MCMC samples for a small reference set of $h$, to estimate the BF among them using the reverse logistic estimator, and then get new samples to estimate $\{m_h(\by)/m_{h'}(\by), h \in {\cal H} \}$ using the generalized IS method. Below, we describe the MCMC algorithms and the practical concern of how long to run them, which illustrates the importance of calculating a SE. 

The two high-dimensional integrals lend the posterior density in \eqref{eq:pih} intractable. But there exist MCMC algorithms to sample from the augmented posterior distribution, that is 
%$\pi_{h}(\psi, \bz | \by)$. 
%Specifically, 
%\[
%   f_h(\by, \bz | \psi) = f(\by | \bz) f_h(\bz | \psi) .
%\]
%%Note that
%%\[\int_{\mathcal{R}^{\m}} f_h(\by, \bz | \psi) d\bz = L_h(\psi |
%%\by). \]
%The so-called {\it complete} posterior density 
\begin{equation}\label{eq:pi-da}
   \pi_{h}( \psi, \bz| \by) %= \frac{f(\by, \bz | \psi, h) \pi(\psi)}{{ m_h}(\by)} 
   =  \frac{f(\by | \bz) f_h(\bz | \psi)\pi(\psi)}{{ m_h}(\by)}.
\end{equation}
Note that $\int_{\mathcal{R}}\pi_{h}(\psi, \bz | \by) d\bz = \pi_{h}(\psi | \by)$. 
Hence, a two-component Gibbs sampler that updates $\psi$ and $\bz$ in turn from their respective conditional distributions based on \eqref{eq:pi-da} yields a Markov chain
$\{\psi^{(i)}, \bz^{(i)}\}_{i \ge 1}$ with stationary distribution $\pi_h( \psi, \bz| \by)$. As a result, the marginal 
$\{\psi^{(i)}\}_{i \ge 1}$ is also a Markov chain with stationary distribution ${ \pi_h}(
\psi| \by)$ \citep{tann:wong:1987}.

As a starting point, we use a small pilot study to identify a range for $h = (\phi, \omega)$ that corresponds to reasonably large BF values. This step is carried out by obtaining the reverse logistic estimator of BF at a coarse grid of $h$ values over a wide area, based on short runs of Markov chains. Specifically, $(\phi, \omega) \in [80, 200] \times [0.2 , 2]$ and within this range the minimum BF is about 1\% the size of the maximum. To more carefully estimate BF over this range, we examine a fine grid ${\cal H}$ that consists of $130$ different $h$ values, with increments of size $10$ for the $\phi$ component, and that of $.2$ for the $w$ component. 

\begin{figure}[htb]
  \begin{center}
     \includegraphics[width=.48\linewidth]{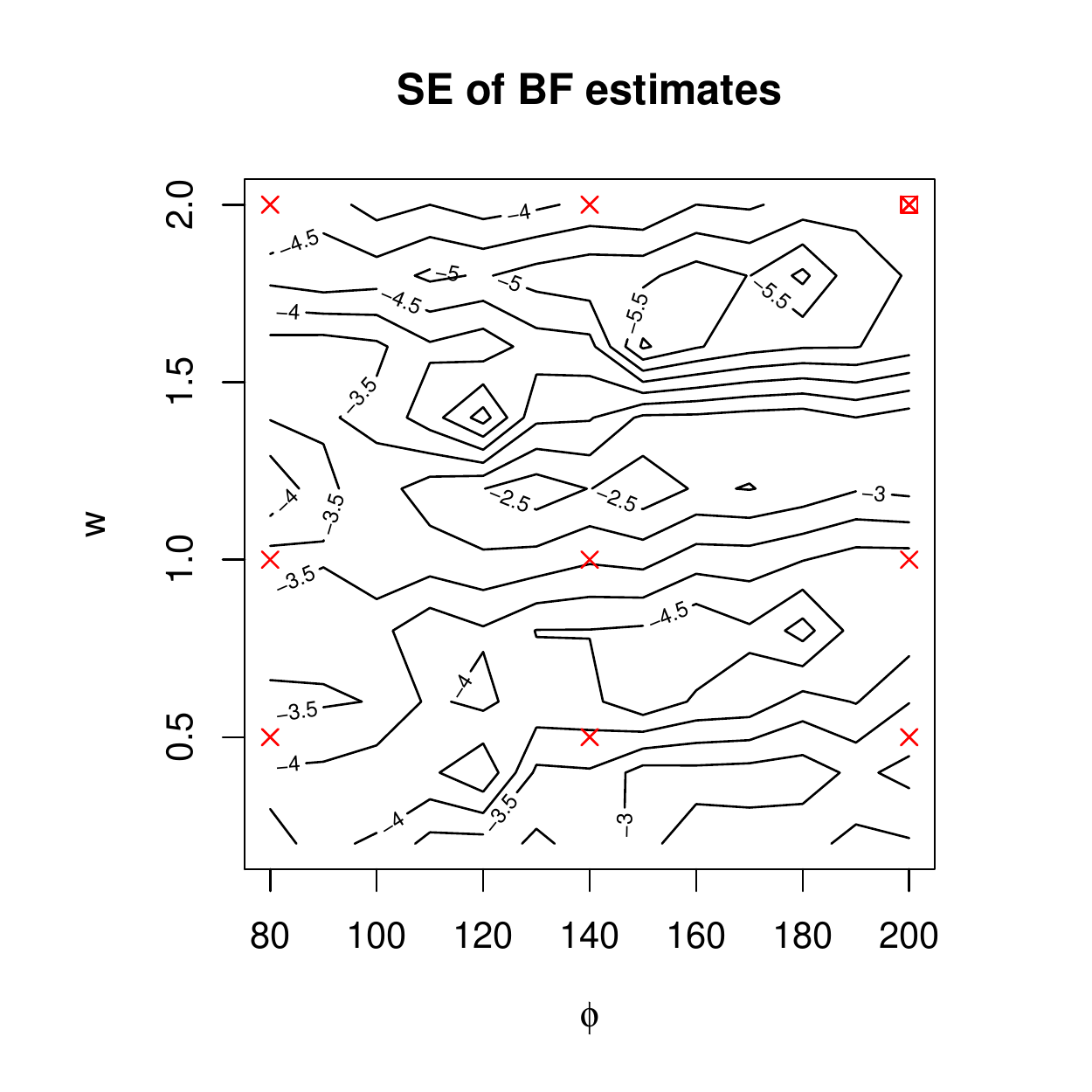}
     \includegraphics[width=.48\linewidth]{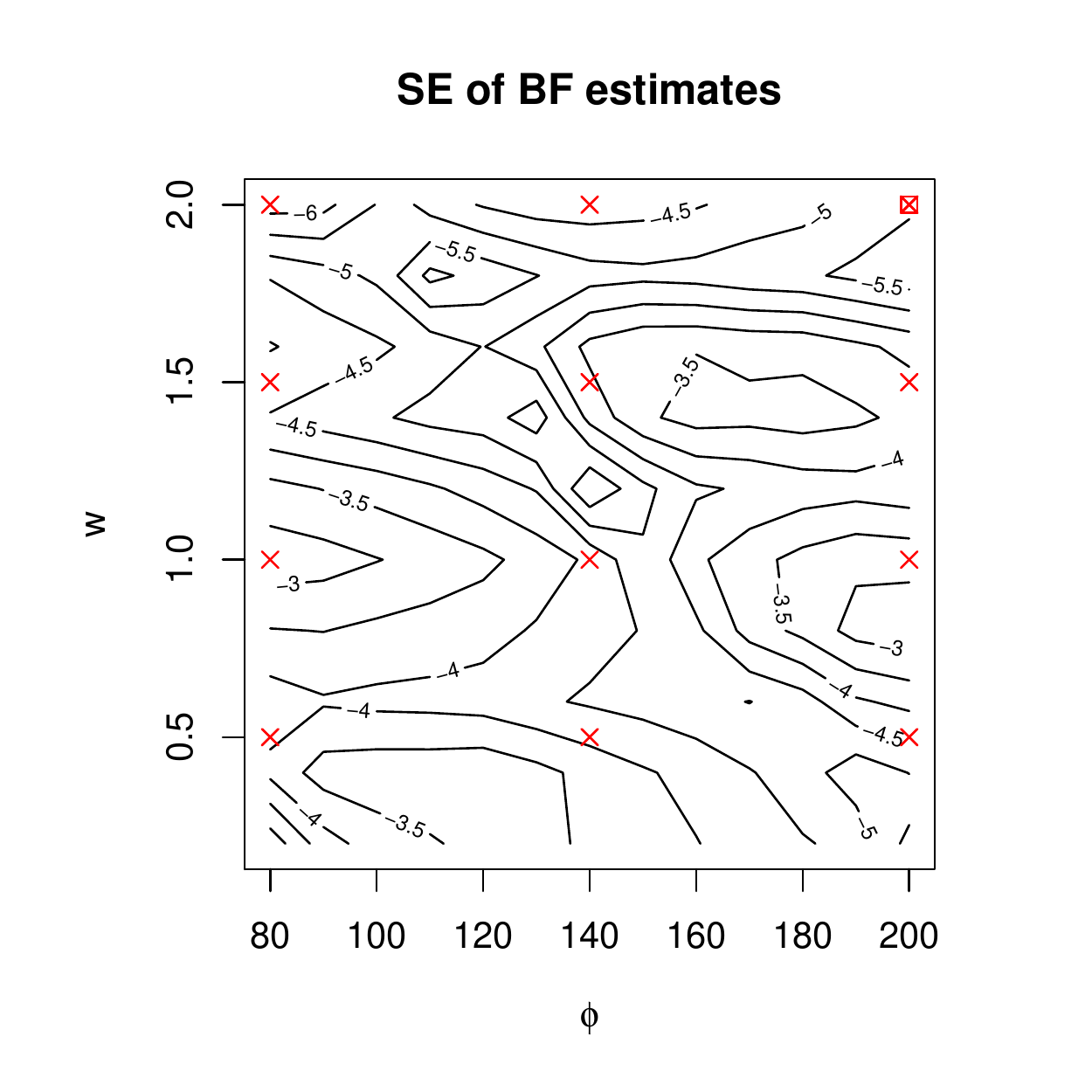}
       \end{center}
       \caption{
Contour plots of the SEs (in log scale) evaluated for the BF estimates. The two plots are based on the original and alternative designs. Skeleton points used in each design are marked by crosses, with the baselines marked by boxed crosses.}
\label{fig:root2se}
\end{figure}

A natural choice for the set of skeleton points is $\{80, 140, 200\} \times \{0.5, 1,  2\}$, with an arbitrarily chosen baseline at $(200, 2)$. We first experiment with samples of sizes $n_1=\cdots = n_9 = 500$ at the skeleton points (after a burn-in period of $500$ iterations and a thinning procedure that keeps one sample every $10$ iterations), of which the first $80\%$ are used in stage~1, and the remaining in stage~2 of the generalized IS procedure. BF estimates at all $h \in {\cal H}$ are obtained, though not shown.  Given the current Monte Carlo sample sizes, it is natural to consider how trustworthy these BF estimates are.  To this end, Figure~\ref{fig:root2se} shows the point-wise SEs over the entire plot obtained via the BM method.  In this setting for some $h$, the magnitude of the SE is over $63\%$ of the corresponding BF estimate.  Hence, inference based on these BF estimates could be questioned because of the high computational uncertainty.  

\begin{figure}[htb]
  \begin{center}
       \includegraphics[width=\linewidth]{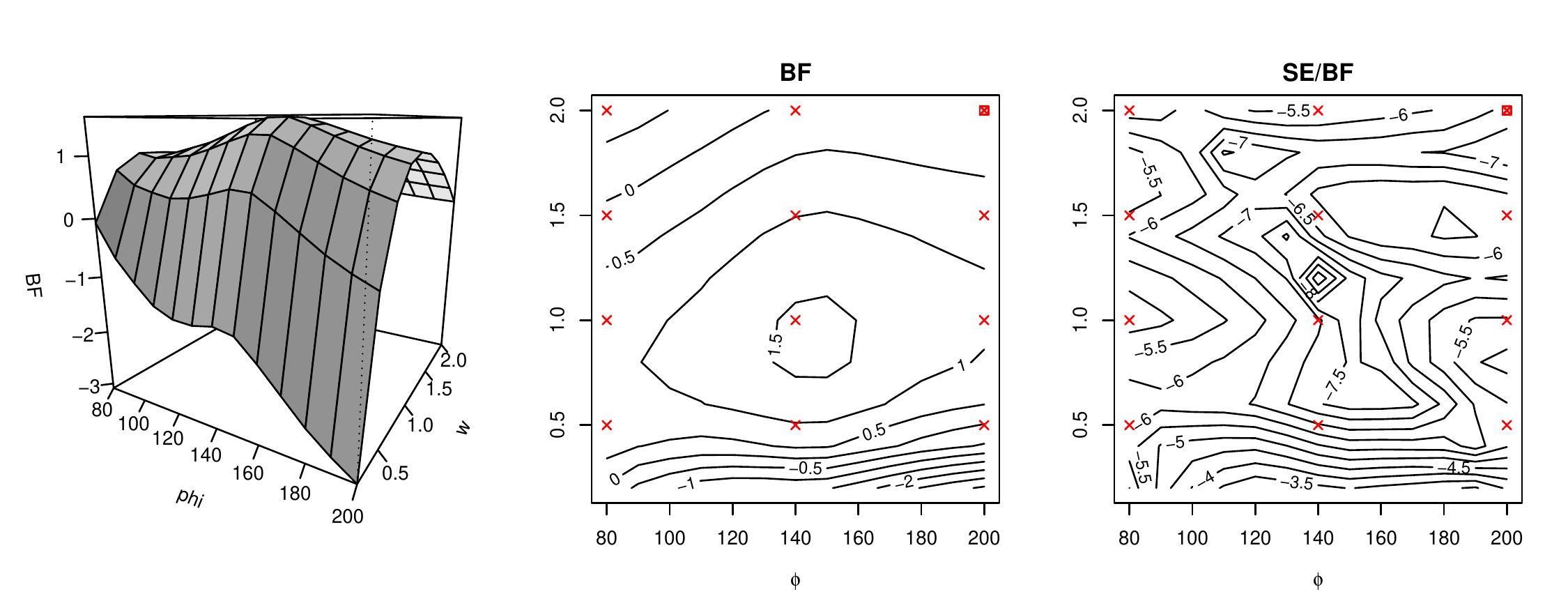}
  \end{center}
  \caption{
  Based on samples of size $6000$ of each of the $16$ Markov chains, the left and the middle panels display a surface plot and a contour plot for the BF estimates (in log scale), respectively. The right panel shows the ratio of the SE to the BF estimates (in log scale), where the SEs are evaluated using the BM method.}
\label{fig:rootbf}
\end{figure}

Suppose we wish to reduce the relative SE to $5\%$ or less for all $h\in {\cal H}$, then roughly $(60\%/5\%)^2=144$ times as many samples would be required under the current design.  Instead, we consider an alternative design using a bigger set of skeleton points, $\{80, 140, 200\} \times \{0.5, 1, 1.5, 2\}$ keeping the baseline unchanged at $(200, 0.2)$. In this alternative design, with sample sizes $n_1= \cdots = n_{12} = 500$, the largest relative SE is $17\%$.  This is almost one fourth that of the previous design, but the computing time only increased by $40\%$.  Accordingly, running the alternative design would achieve the computing goal much faster. In this case, we increase the sample sizes to $n_1= \cdots =n_{12} = 6000$, which is approximately $(17\%/5\%)^2$ times $500$. Overall, the new process takes $2$ minutes to run on a 3.4GHz Intel Core i7 running linux.  The resulting BF estimates are shown in Figure~\ref{fig:rootbf}, with maximum relative SE reduced to $5.0\%$. For the sake of comparison, running the original design for the same amount of time allows a common sample size of $n_i=8000$ resulting in a maximum relative SE of $14.4\%$. In short, easily obtainable SE estimates allows us to experiment, choose among different designs, and perform samples size calculations.

The simplicity of the method matters when it comes to estimating SEs in practice. Using the BM method to obtain SE requires no extra input beyond what is needed for obtaining the generalized IS estimates. Indeed, as long as one can run existing software to obtain the Markov chain samples, there is no need to know the Markov transition kernels utilized in the background. Unlike the BM method, the RS method depends on identifying regeneration times, typically through constructing minorization conditions for the Markov transition kernels (see \citet{mykl:tier:yu:1995} for details).  Despite the fact that minorization conditions can be established for any Markov transition kernel, we demonstrate for the current example the amount of effort needed to obtain a regeneration can be prohibitively high. Recall the MCMC scheme involves sampling from $\pi_h( \psi| \bz, \by)$ and $\pi_h( \bz| \psi, \by)$ in turn. The former is a standard distribution hence easy to sample from. The latter is not, and we followed \citet{digg:tawn:moye:1998} that updates $z_j, j=1,\cdots, {\m}$ in turn, each using a one-dimensional Metropolis-Hastings step that keeps invariant the conditional posterior distribution of $z_j$ given all other components. Denote the transition density of these MH steps as $f_1,\cdots, f_M$, and suppressing the notations of their dependence on $\by$, the transition kernel of the Markov chain can be represented as
\[
p(\bz', \psi' | \bz, \psi) = 
f_1(z'_1| z_2,\cdots, z_M, \psi) f_2(z'_2| z_1', z_3,\cdots, z_M, \psi) \cdots f_n(z'_M | z_1', \cdots, z'_{{\m}-1}, \psi) \pi_h(\psi' | \bz')\,.
\]
According to a common method described in \citet{jone:hobe:2004}, one can build a minorization condition by finding $D \subset \Real^{{\m}}\times \Real \times \Real^+$, $\epsilon>0$, and $k(\cdot)$ such that, %for all $(\bz, \psi)\in D$, 
\[p(\bz', \psi' | \bz, \psi) \geq \epsilon I_D(\bz, \psi)\, k(\bz', \psi' )  \;\;\text{for all $(\bz', \psi') \in \Real^{{\m}}\times \Real \times \Real^+$\,.}
  \]
Further, the above condition can be established if %for all $(\bz, \psi)\in D$, 
%\[f_1(z'_1| z_2,\cdots, z_n, \psi) f_2(z'_2| z_1', z_3,\cdots, z_n, \psi) \cdots f_n(z'_n | z_1', \cdots, z'_{{\m}-1}, \psi)
%\geq \,\epsilon_1 k_1(z'_1) \, \epsilon_2 k_2(z'_1, z'_2) \, \cdots \, \epsilon_n k_n(z'_1, \cdots, z'_n)\] %f(\mu', \]
\[\begin{split}
 & f_1(z'_1| z_2,\cdots, z_M, \psi) f_2(z'_2| z_1', z_3,\cdots, z_M, \psi) \cdots f_M(z'_M | z_1', \cdots, z'_{{\m}-1}, \psi)\pi_h(\psi' | \bz')\\ %f(\mu', {\sigma^2}' | bz') \\
\geq & I_D(\bz, \psi)\,\epsilon_1 k_1(z'_1) \, \epsilon_2 k_2(z'_1, z'_2) \, \cdots \, \epsilon_M k_M(z'_1, \cdots, z'_M)\pi_h(\psi' | \bz') %f(\mu', {\sigma^2}' | bz') 
  \;\;\text{for all $(\bz', \psi') \in \Real^{{\m}}\times \Real \times \Real^+$,}
\end{split}\]
where the common term $\pi_h$ on both sides of the inequality will cancel, and hence the work is in finding  $\epsilon_1,\cdots, \epsilon_M,$ and $k_1(\cdot),\cdots, k_M(\cdot)$. It's easy to see that the smaller the set $D$, the larger $\epsilon= \Pi_{i=1}^{{\m}} \epsilon_j$ can possibly be, where $\epsilon$ can be interpreted as the conditional regeneration rate given $D$ is visited. Suppose we take $D$ to be small enough such that $\epsilon_j$ takes on a very large value of $0.8$ for each $j$, then the probability of getting a regeneration given a visit to $D$ is $\epsilon= (0.8)^{100}  \approx 2\times 10^{-10}$. Being overoptimistic that the Markov chain visits $D$ with probability close to $1$, it would still take 100 billon iterations for the chain to regenerate about twenty times, barely enough for the RS method to be effective. 

Using the EB estimate $\hat{h}$ of $h$, estimation of the remaining 
parameters $\psi$ and prediction of the spatial random field can be done 
in the standard method using MCMC samples from $\pi_{\hat{h}} (\psi | 
\by)$ \cite[see e.g.][section 3.2]{roy:evan:zhu:2016}. Thus we can 
produce the root rot disease prediction map similar to that in \citet[][Web Fig. 
10]{roy:evan:zhu:2016}.

%proposal $z_j^*\sim N(-\sum_l Q_{kl} b_l (Q_{kk})^{-1})$, where $Q$ is the inverse of the variance-covariance matrix of $\bz|$

% \blue{Finally, based on the empirical Bayes choice of the hyperparameter $(\phi, \omega)=(150, 0.8)$, we can obtain a map of the predicted disease rate as well as its point-wise variance over the entire farm. The formula and the maps are not shown due to space concerns, and we point interesting readers to \citet{roy:evan:zhu:2016} for details.}
%\red{Do we need to do wrap up the example with a disease rate map? Using the empirical choice of $h=(150, .8)$. If so, details of the technique can be put in an appendix? or simply cite roy:evan:zhu. 
%\begin{scriptsize}
%Further, prediction about $Z_0$, the values of
%$Z(s)$ at any set of locations of interest, $(s_{01}, s_{02},
%\dots,s_{0k})$ can be made according to 
%\[
%f(\bz_0 | \by) = \int_{\mathcal{R} \times \mathcal{R}_+} \int_{\mathcal{R}^{\m}} f(\bz_0 | \bz, \psi, g) \pi_{g} (\psi, \bz | \by) d\bz d\psi\,,
%\]
%where $f(\bz_0 | \bz, \psi, g)$ is some multivariate normal density with easy-to-obtain mean and variance covariance matrix.
%And $\E(t(\bz_0)|\by)$ can be approximated by
%$\sum_{i=1}^N t^*(\bz^{(i)}, \psi^{(i)}, \hat{g})/N$, if $t^*(\bz, \psi, \hat{g}) \equiv \E(t(\bz_0)|\bz, \psi, \hat{g})$ is easy to evaluate. Otherwise, approximate $\E(t(\bz_0)|\by)$ with $\frac{1}{N} \sum_{i=1}^N t(\bz_0^{(i)}))$.
%\end{scriptsize}
%}

\section{Discussion}
\label{sec:disc}
In this paper we consider two separate but related problems. The first
problem is estimating the ratios of unknown normalizing constants
given Markov chain samples from each of the $k>1$ probability densities. 
The second problem is estimating expectations of a function with respect to a large number of
probability distributions. These problems are related in the sense
that generalized IS estimators used for the latter utilize estimates
derived when solving the first problem. The first situation also arises
in a variety of contexts other than the generalized IS estimators. 

For both problems, we derive estimators with
flexible weights and thus these estimators are appropriate
for Markov chains with different mixing behaviors. We establish CLTs
for these estimators and develop BM methods for consistently
estimating their SEs. These easy to calculate SEs
 are important for at least three reasons. First, SEs are needed to assess
the quality of the estimates. Second, our ability to calculate
SEs allows us to search for optimal weights $\ba$ for both
stage 1 and 2. And last but not least, SEs form the
basis for comparison of generalized IS with other available methods for
estimating large number of (ratios of) normalizing constants.

Although we compare BM and RS in this paper, spectral estimators can also be derived for variance estimation using the results in \cite{vats:fleg:jone:2015:spectral}. However, estimation by spectral methods is generally more expensive computationally.  Further, \cite{fleg:jone:2010} compare the performance of confidence intervals produced by BM, RS and spectral methods for the time average estimator, and they conclude that if tuning parameters are chosen appropriately, all these three methods perform equally well. Control variates can be used to further improve the accuracy of our generalized IS estimators \citep{owen:zhou:2000, doss:2010}. A direction of future research would be to establish a BM estimator of the SEs for control variate based methods.

\begin{appendix}
\section{Proof of Theorem~\ref{thm:CLT}}
\label{sec:appthm1}
  In the appendices, we denote \citet{DossTan:2014} by D\&T. 
  The proof of the consistency of $\hatbd$ follows from
  D\&T section~A.1 and is omitted. Establishing a CLT for
  $\hatbd$ is analogous to section~A.2 of D\&T, but there are
  significant differences. Below we establish the CLT for
  $\hatbd$ and finally show that $\widehat{V}$ is a consistent
  estimator of $V$.

  We begin by considering ${n}^{1/2} (\hatbzeta - \bzeta_0)$.  As
  before, let $\nabla$ represents the gradient operator. % As in the
  %classical proof of asymptotic normality of maximum likelihood
  %estimators, 
  We expand $\nabla \ell_{n}$ at $\hatbzeta$ around
  $\bzeta_0$, and using the appropriate scaling factor, we get
\begin{equation}
  \label{eq:TE-score}
  -{n}^{-1/2}\bigl( \nabla \ell_{n} (\hatbzeta) -
  \nabla \ell_{n} (\bzeta_0) \bigr) = -{n}^{-1} \nabla^2
  \ell_{n} (\bzeta_*)  {n}^{1/2}(\hatbzeta - \bzeta_0),
\end{equation}
where $\bzeta_*$ is between $\hatbzeta$ and $\bzeta_0$.  Consider
the left side of~\eqref{eq:TE-score}, which is just $n^{1/2}
n^{-1} \nabla \ell_{n} (\bzeta_0)$, since $\nabla \ell_{n}
(\hatbzeta) = 0$.  There are several nontrivial components to the
proof, so we first give an outline.
\begin{enumerate}
\item Following D\&T we show that each element of the vector $n^{-1} \nabla
  \ell_n (\bzeta_0)$ can be represented as a linear combination
  of mean $0$ averages of functions of the $k$ chains. % plus a vanishingly small term.
\item Based on Step~$1$, applying CLT for each of the $k$ Markov chain averages, we obtain a CLT for
  the scaled score vector. In particular, we show that
  $n^{1/2} n^{-1} \nabla \ell_n (\bzeta_0) \cd {\cal N}
  (0, \Omega)$, where $\Omega$ defined in \eqref{eq:Omega} involves infinite sums of auto-covariances of each chain.
\item Following \citet{Geyer:1994} it can be shown that $-n^{-1} \nabla^2 \ell_n (\bzeta_*) \cas B$
  and that $\bigl( -n^{-1} \nabla^2 \ell_n (\bzeta_*)
  \bigr)^{\dagger} \cas B^{\dagger}$, where $B$ is defined
  in \eqref{eq:B}.
\item We conclude that $n^{1/2} (\hatbzeta - \bzeta_0) \cd
  {\cal N} (0, B^{\dagger} \Omega B^{\dagger})$.
\item Since $\bd = g(\bzeta_0)$ and $\hatbd = g(\hatbzeta)$, where $g$
  is defined in \eqref{eq:g}, by the Delta method it follows that $
  n^{1/2} (\hatbd - \bd) \cd {\cal N} (0, V)$ where $V = D^{\top}
  B^{\dagger} \Omega B^{\dagger} D.$
\end{enumerate}
We now provide the details.
\begin{enumerate}
\item Start by considering $n^{-1} \nabla \ell_n(\bzeta_0)$.
  For $r = 1, \ldots, k$, from D\&T we have
  \begin{equation}
    \label{eq:grad-lin-comb}
    \begin{split}
      \frac{\partial \ell_n(\bzeta_0)} {\partial
      {\zeta_r}} & = w_r \sum_{i=1}^{n_r} \bigl( 1 - p_r(X_i^{(r)},
                     \bzeta_0) \bigr) - \sum_{\substack{l=1 \\l \neq
                     r}}^k w_l \sum_{i=1}^{n_l} p_r(X_i^{(l)},
                     \bzeta_0) \\
                \text{(can be shown to)}\;\; & = w_r \sum_{i=1}^{n_r} \Bigl( 1 - p_r(X_i^{(r)},
                     \bzeta_0) - \bigl[ 1 - E_{\pi_r} \bigl( p_r(X,
                     \bzeta_0) \bigr) \bigr] \Bigr) \\
                 & \hspace{8mm} - \sum_{\substack{l=1 \\l \neq r}}^k
                     w_l \sum_{i=1}^{n_l} \bigl[ p_r(X_i^{(l)},
                     \bzeta_0) - E_{\pi_l} \bigl( p_r(X, \bzeta_0)
                     \bigr) \bigr] .
    \end{split}
  \end{equation}
  That is, \eqref{eq:grad-lin-comb} can be
  used to view $n^{-1} \partial \ell_n(\bzeta_0) / \partial
  {\zeta_r}$ as a linear combination of mean $0$ averages of
  functions of the $k$ chains.

\item Next, we need a CLT for the vector $ \nabla \ell_n(\bzeta_0) =(\partial \ell_n(\bzeta_0) / \partial
  {\zeta_1},\cdots, \partial \ell_n(\bzeta_0) / \partial
  {\zeta_k})^T$, that is, to show that $n^{-1/2} \nabla \ell_n(\bzeta_0) \cd N(0,
    \Omega)$ as $n \rightarrow \infty$.  Note that,
    \begin{align*}
      \frac{1}{\sqrt{n}} \frac{\partial \ell_n(\bzeta_0)}
      {\partial
      {\zeta_r}} & = - \frac{1}{\sqrt{n}} \sum_{l=1}^k w_l
                     \sum_{i=1}^{n_l} \bigl[ p_r(X_i^{(l)},
                     \bzeta_0) - E_{\pi_l} \bigl( p_r(X, \bzeta_0)
                     \bigr) \bigr] \nonumber\\
                 & = - \sum_{l=1}^k \sqrt{\frac{n}{n_l}} a_l
                     \frac{1}{\sqrt{n_l}} \sum_{i=1}^{n_l} \bigl[ p_r(X_i^{(l)},
                     \bzeta_0) - E_{\pi_l} \bigl( p_r(X, \bzeta_0)
                     \bigr) \bigr] \nonumber\\
&=- \sum_{l=1}^k \sqrt{n} a_l \bar{\bar{Y}}^{(r,l)},
  \end{align*}
  where $\bar{\bar{Y}}^{(r,l)} := \frac{1}{n_l} \sum_{i=1}^{n_l}
  Y_i^{(r,l)}$ and $Y_i^{(r,l)}$ is as defined in \eqref{eq:Ys}.
  Since $p_r(x, \bzeta) \in (0, 1)$ for all $x$, $r$ and $\bzeta$,
  we have $E_{\pi_l} \bigl( | p_r(X, \bzeta_0) - E_{\pi_l} (
    p_r(X, \bzeta_0)) |^{2+\delta} \bigr) < \infty$ for any $\delta >
    0$. Then since $\Phi_l$ is polynomially ergodic of order $m > 1$, we have asymptotic normality for the univariate quantities
    $\sqrt{n_l}\bar{\bar{Y}}^{(r,l)}$ (see e.g.\ Corollary 2 of \citet{jone:2004}).  Since $n_l/n \rightarrow s_l$
  for $l= 1,\dots,k$ and $a_l$'s are known, by independence of the
  $k$ chains, we conclude that
  \begin{equation*}
    \frac{1}{\sqrt{n}} \frac{\partial \ell_{n}(\bzeta_0)}
    {\partial {\zeta_r}} \cd {\cal N} (0, \Omega_{rr})
    \text{ as } n \rightarrow \infty,
  \end{equation*}
  where $\Omega$ is defined in \eqref{eq:Omega}.  Next, we extend the component-wise CLT to a joint CLT. Consider any $\bt\in (t_1,\cdots,t_k) \in \Real^k$, we have
  \[\begin{split}  & t_1    \frac{1}{\sqrt{n}} \frac{\partial \ell_{n}(\bzeta_0)}
    {\partial {\zeta_1}}+\cdots  +t_k    \frac{1}{\sqrt{n}} \frac{\partial \ell_{n}(\bzeta_0)}
    {\partial {\zeta_k}}\\
    =& -\sum_{l=1}^k\left(t_1 \sqrt{n} a_l \frac{\sum_{i=1}^{n_l} Y_i^{(1,l)}}{n_l}+\cdots + t_k \sqrt{n} a_l \frac{\sum_{i=1}^{n_l} Y_i^{(k,l)}}{n_l}\right) \\
    =& - \sum_{l=1}^k\sqrt{\frac{n}{n_l}} a_l \frac{
    \sum_{i=1}^{n_l}\left(t_1   Y_i^{(1,l)}+\cdots + t_k   Y_i^{(k,l)}\right)
    }{\sqrt{n_l}}\cd  {\cal N} (0, \bt^T \Omega \bt)
    \text{ as } n \rightarrow \infty.
    \end{split}\]
  Hence, the Cram\'{e}r-Wold
  device implies the joint
  CLT, \begin{equation}\label{eq:jtCLT} n^{-1/2} \nabla \ell_n (\bzeta_0) \cd {\cal N}
  (0, \Omega) \qquad \text{as } n \rightarrow \infty.\end{equation}
  \end{enumerate}
Steps 3-5 are omitted since the derivations are basically the same as in D\&T.

  Next we provide a proof of the consistency of the estimate of the
  asymptotic covariance matrix $V$, that is, we show that
  $\widehat{V} \equiv \widehat{D}^{\top} \widehat{B}^{\dagger}
  \widehat{\Omega} \widehat{B}^{\dagger} \widehat{D} \cas V \equiv
  D^{\top} B^{\dagger} \Omega B^{\dagger} D$ as $n \rightarrow
  \infty$. Since $\hatbzeta \cas \bzeta_0$ and $\hatbd \cas \bd$, it
  implies that $\widehat{D} \cas D$. From D\&T, we know that
  $\widehat{B} \cas B$ and using the spectral representation of
  $\widehat{B}$ and of $B$, it follows that $\widehat{B}^{\dagger}
  \cas B^{\dagger}$.

  To complete the proof, we now show that $\widehat{\Omega} \cas
  \Omega$ where the BM estimator $\widehat{\Omega}$ is
  defined in \eqref{eq:Omegahat}. This will be proved in couple of
  steps. First, we consider a single chain $\Phi_l$ used to calculate
  $k$ quantities and establish a multivariate CLT. We use the results
  in \citet{vats:fleg:jone:2015:output} who obtain conditions for the
  nonoverlapping BM estimator to be strongly consistent in
  multivariate settings. Second, we combine results from the $k$
  independent chains.  Finally, we show that $\widehat{\Omega}$ is a
  strongly consistent estimator of $\Omega$.

Denote $\bar{\bar{Y}}^{(l)} = \left( \bar{\bar{Y}}^{(1,l)}, \bar{\bar{Y}}^{(2,l)}, \dots, \bar{\bar{Y}}^{(k,l)} \right)^{\top}$.  Similar to deriving \eqref{eq:jtCLT} via the Cram\'{e}r-Wold device, we have the following joint CLT for $W_l$:, $\sqrt{n_l}\bar{\bar{Y}}^{(l)} \cd {\cal N} (0, \Sigma^{(l)})$ as $n_l \rightarrow \infty$, where $\Sigma^{(l)}$ is a $k \times k$ covariance matrix with
\begin{equation}\label{eq:Sigmars}
 \Sigma^{(l)}_{rs} = E_{\pi_l}\{Y_1^{(r,l)} Y_1^{(s,l)}\} +  \sum_{i=1}^{\infty} E_{\pi_l}\{Y_1^{(r,l)} Y_{1+i}^{(s,l)}\}+  \sum_{i=1}^{\infty} E_{\pi_l}\{Y_{1+i}^{(r,l)} Y_1^{(s,l)}\} .
\end{equation}

The nonoverlapping BM estimator of $\Sigma^{(l)}$ is given in
\eqref{eq:BM}. We now prove the strong consistency of
$\widehat{\Sigma}^{(l)}$. Note that $\widehat{\Sigma}^{(l)}$ is
defined using the terms $\bar{Z}^{(r,l)}_{m}$'s which involve the
random quantity $\hatbzeta$. We define $\widehat{\Sigma}^{(l)}
(\bzeta_0)$ to be $\widehat{\Sigma}^{(l)}$ with $\bzeta_0$ substituted
for $\hatbzeta$, that is,
\[
\widehat{\Sigma}^{(l)} (\bzeta_0) = \frac{b_l}{e_l - 1} \sum_{m=0}^{e_l - 1} \left[ \bar{Y}^{(l)}_{m} - \bar{\bar{Y}}^{(l)} \right] \left[ \bar{Y}^{(l)}_{m} - \bar{\bar{Y}}^{(l)} \right]^{\top} \; \mbox{ for} \;\; l =1,\ldots, k,
\]
where $\bar{Y}^{(l)}_{m}= \left( \bar{Y}^{(1,l)}_{m}, \ldots,
  \bar{Y}^{(k,l)}_{m} \right) ^{\top}$ with $\bar{Y}^{(r,l)}_{m} :=
\sum_{j=m b_l + 1}^{(m+1) b_l} Y_j^{(r,l)}/ b_l$. We prove $\widehat{\Sigma}^{(l)} \cas
\Sigma^{(l)}$ in two steps: (1) $\widehat{\Sigma}^{(l)} (\bzeta_0) \cas\Sigma^{(l)}$ and (2) $\widehat{\Sigma}^{(l)} - \widehat{\Sigma}^{(l)} (\bzeta_0) \cas 0$.  Strong consistency of the multivariate BM estimator
$\widehat{\Sigma}^{(l)} (\bzeta_0)$ requires both $e_{l} \to \infty$
and $b_{l} \to \infty$.  Since for all $r$, $E_{\pi_l} \bigl( |
  p_r(X, \bzeta_0) - E_{\pi_l} ( p_r(X, \bzeta_0)) |^{4+\delta} \bigr)
  < \infty$ for any $\delta > 0$, $\Phi_l$ is polynomially ergodic of
  order $m > 1$, and $b_l = \lfloor n_l^{\nu}
  \rfloor$ where $1 > \nu > 0$, it follows from \citet{vats:fleg:jone:2015:output} that $\widehat{\Sigma}^{(l)}
  (\bzeta_0) \cas \Sigma^{(l)}$ as $n_l \rightarrow \infty$. We show
$\widehat{\Sigma}_{rs}^{(l)} - \widehat{\Sigma}_{rs}^{(l)} (\bzeta_0)
\cas 0$ where $\widehat{\Sigma}_{rs}^{(l)}$ and
$\widehat{\Sigma}_{rs}^{(l)} (\bzeta_0)$ are the $(r,s)$th elements of
the $k \times k$ matrices $\widehat{\Sigma}_{rs}^{(l)}$ and
$\widehat{\Sigma}_{rs}^{(l)} (\bzeta_0)$ respectively. By the mean
value theorem (in multiple variables), there exists
$\bzeta^*=t\hatbzeta +(1-t) \bzeta_0$ for some $t\in(0,1)$, such that
\begin{equation}
  \label{eq:Sigmv}
\widehat{\Sigma}_{rs}^{(l)} - \widehat{\Sigma}_{rs}^{(l)} (\bzeta_0) =
\nabla \widehat{\Sigma}_{rs}^{(l)} (\bzeta^*) \cdot (\hatbzeta -
\bzeta_0),
\end{equation}
 where
$\cdot$ represents the dot product. Note that
 \[
 \widehat{\Sigma}_{rs}^{(l)} (\bzeta) = \frac{b_l}{e_l - 1} \sum_{m=0}^{e_l - 1}  [\bar{Z}^{(r,l)}_{m} (\bzeta) - \bar{\bar{Z}}^{(r,l)} (\bzeta)] [\bar{Z}^{(s,l)}_{m} (\bzeta) - \bar{\bar{Z}}^{(s,l)} (\bzeta)],
 \]
 where $\bar{Z}^{(r,l)}_{m}(\bzeta) := \sum_{j=m b_l + 1}^{(m+1) b_l}
 p_r(X_j^{(l)}, \bzeta)/b_l$ and $\bar{\bar{Z}}^{(r,l)} (\bzeta) := \sum_{j=1}^{n_l}
 p_r(X_j^{(l)}, \bzeta)/n_l$. Some calculations show that for $t \neq r$
 \[
 \frac{\partial \bar{Z}^{(r,l)}_{m} (\bzeta)}{\partial \bzeta_t} = - \frac{1}{b_l} \sum_{j=m b_l + 1}^{(m+1) b_l}
 p_r(X_j^{(l)}, \bzeta) p_t(X_j^{(l)}, \bzeta)
 \]
 and
 \[
 \frac{\partial \bar{Z}^{(r,l)}_{m} (\bzeta)}{\partial \bzeta_r} = \frac{1}{b_l} \sum_{j=m b_l + 1}^{(m+1) b_l}
 p_r(X_j^{(l)}, \bzeta) (1 - p_r(X_j^{(l)}, \bzeta)).
 \]
 We denote $\bar{U}^r_m := \bar{Z}^{(r,l)}_{m} (\bzeta) - E_{\pi_l}
 [p_r(X, \bzeta)]$, $\bar{\bar{U}}^r := \bar{\bar{Z}}^{(r,l)} (\bzeta)
 - E_{\pi_l} [p_r(X, \bzeta)]$, and similarly the centered versions of
 $\partial \bar{Z}^{(r,l)}_{m} (\bzeta)/\partial \bzeta_t$ and
 $\partial \bar{\bar{Z}}^{(r,l)} (\bzeta)/\partial \bzeta_t$ by
 $\bar{V}^{(r,t)}_m$ and $\bar{\bar{V}}^{(r,t)}$ respectively. Since
 $p_r(X, \bzeta)$ is uniformly bounded by 1 and $\Phi_l$ is
 polynomially ergodic of order $m > 1$,
 there exist $\sigma^2_r, \tau^2_{r,t} < \infty$ such that $\sqrt{b_l}
 \bar{U}^r_m \cd N(0, \sigma^2_r), \sqrt{n_l} \bar{\bar{U}}^r \cd N(0,
 \sigma^2_r)$, $\sqrt{b_l} \bar{V}^{(r,t)}_m \cd N(0, \tau^2_{r, t})$,
 and $\sqrt{n_l} \bar{\bar{V}}^{(r,t)} \cd N(0, \tau^2_{r, t})$.  We
 have
\begin{align*}
  \label{eq:delSig}
 & \frac{\partial  \widehat{\Sigma}_{rs}^{(l)} (\bzeta) }{\partial \bzeta_t} = \frac{1}{e_l - 1} \sum_{m=0}^{e_l - 1} [\sqrt{b_l}(\bar{U}^r_m - \bar{\bar{U}}^r) \sqrt{b_l}(\bar{V}^{(s,t)}_m - \bar{\bar{V}}^{(s,t)}) + \sqrt{b_l}(\bar{V}^{(r,t)}_m - \bar{\bar{V}}^{(r,t)}) \sqrt{b_l}(\bar{U}^s_m - \bar{\bar{U}}^s)]\\
=& \frac{1}{e_l - 1} \sum_{m=0}^{e_l - 1}
 \left[ \sqrt{b_l}\bar{U}^r_{m}\sqrt{b_l} \bar{V}^{(s,t)}_{m} + \sqrt{b_l}\bar{V}^{(r,t)}_{m}\sqrt{b_l} \bar{\bar{U}}^s_{m}\right]
 - \frac{1}{e_l - 1} \left[ \sqrt{n_l} \bar{\bar{U}}^r \sqrt{n_l} \bar{\bar{V}}^{(s,t)} + \sqrt{n_l} \bar{\bar{V}}^{(r,t)} \sqrt{n_l} \bar{\bar{U}}^s\right].
\end{align*}
It is easy to see that the negative term in the above expression goes to zero as $e_l \rightarrow \infty$.
Further, since
 \[ \left|\sqrt{b_l}\bar{U}^r_{m}\sqrt{b_l} \bar{V}^{(s,t)}_{m}\right|
 \leq \frac{1}{2}\left[b_l(\bar{U}^r_{m})^2\right] +\frac{1}{2}\left[ b_l(\bar{V}^{(s,t)}_{m})^2\right],
 \]
 we have
\[
\left|\frac{1}{e_l - 1} \sum_{m=0}^{e_l - 1}
  \sqrt{b_l}\bar{U}^r_{m}\sqrt{b_l} \bar{V}^{(s,t)}_{m} \right|
\leq \frac{1}{2}\frac{1}{e_l - 1}\sum_{m=0}^{e_l - 1}
  \left[b_l(\bar{U}^r_{m})^2\right] +\frac{1}{2}\frac{1}{e_l - 1}\sum_{m=0}^{e_l - 1}\left[ b_l(\bar{V}^{(s,t)}_{m})^2\right]
 \cas \frac{1}{2}\sigma^2_{r}+\frac{1}{2}\tau^2_{s,t},
 \]
 where the last step above is due to strong consistency of the BM estimators for the asymptotic variances of
 the sequences $\{p_r(X_j^{(l)}, \bzeta), j=1,\cdots,n_l\}$ and $\{ \partial
 p_s(X_j^{(l)}, \bzeta)/\partial \bzeta_t , j=1,\cdots,n_l\}$ respectively. Similarly,
we have
\[
\left|\frac{1}{e_l - 1} \sum_{m=0}^{e_l - 1}
  \sqrt{b_l} \bar{V}^{(r,t)}_{m} \sqrt{b_l}\bar{U}^s_{m} \right|
\leq \frac{1}{2}\frac{1}{e_l - 1}\sum_{m=0}^{e_l - 1}\left[ b_l(\bar{V}^{(r,t)}_{m})^2 + \frac{1}{2}\frac{1}{e_l - 1}\sum_{m=0}^{e_l - 1}
  \left[b_l(\bar{U}^s_{m})^2\right] \right]
 \cas \frac{1}{2}\tau^2_{r,t} + \frac{1}{2}\sigma^2_{s}.
 \]
 Note that the terms $U_m^r V_{m}^{(r,t)}, \sigma^2_{r},
 \tau^2_{r,t}$, etc, above actually depends on $\bzeta$, and we are
 indeed concerned with the case where $\bzeta$ takes on the value
 $\bzeta^*$, lying between $\hatbzeta$ and $\bzeta_0$. Since,
 $\hatbzeta \cas \bzeta_0$, $\bzeta^* \cas \bzeta_0$ as $n_l
 \rightarrow \infty$. Let $\|u \|$ denotes the $L_1$ norm of a vector
 $u \in \mathbb{R}^k$. So from \eqref{eq:Sigmv}, and the fact that $\partial  \widehat{\Sigma}_{rs}^{(l)} (\bzeta)/\partial \bzeta_t$ is bounded with probability one, we
 have
  \[
    |\widehat{\Sigma}_{rs}^{(l)} - \widehat{\Sigma}_{rs}^{(l)} (\bzeta_0)| \leq \underset{1 \le t \le k}{\max}\left\{\left|\frac{\partial  \widehat{\Sigma}_{rs}^{(l)} (\bzeta^*) }{\partial \bzeta_t}\right|\right\} \|\hatbzeta - \bzeta_0\| \cas 0
    \;\;\text{ as}\; n\rightarrow \infty.
 \]

 Since $\widehat{\Sigma}^{(l)} \cas \Sigma^{(l)},$ for $l =
 1,\dots,k$, it follows that $\widehat{\Sigma} \cas \Sigma$ where
 $\widehat{\Sigma}$ is defined in \eqref{eq:Sigmahat} and $\Sigma$ is
 the corresponding $k^2 \times k^2$ covariance matrix, that
 is, $\Sigma$ is a block diagonal matrix as $\widehat{\Sigma}$ with
 $\Sigma^{(l)}$ substituted for $\widehat{\Sigma}^{(l)}, l=1,\dots,k$.
Since $n_l/n
\rightarrow s_l$ for $l= 1,\dots,k$, we have $A_n \rightarrow A_s$ as $n \rightarrow
\infty$ where $A_n$ is defined in \eqref{eq:defA} and
\[
A_s = \left( - \sqrt{\frac{1}{s_1}} a_1 I_k \quad - \sqrt{\frac{1}{s_2}} a_2 I_k \quad \dots \quad - \sqrt{\frac{1}{s_k}} a_k I_k \right).
\]
Finally from \eqref{eq:Omega} and \eqref{eq:Sigmars} we see that
$\Omega = A_s \Sigma A_s^{T}$. So from \eqref{eq:Omegahat} we have
$\widehat{\Omega} \equiv A_n \widehat{\Sigma} A_n^{T} \cas A_s \Sigma
A_s^{T} = \Omega$ as $n \rightarrow \infty$.

\section{Proof of Theorem~\ref{thm:elnc}}
\label{sec:appthm2}
As in \citet{buta:doss:2011} we write
\begin{equation}
  \label{eq:huminu}
  \sqrt{n}(\hat{u}(\pi, \pi_1; \ba, \hatbd) - u(\pi, \pi_1)) = \sqrt{n}(\hat{u}(\pi, \pi_1; \ba, \hatbd) - \hat{u}(\pi, \pi_1; \ba, \bd)) + \sqrt{n}(\hat{u}(\pi, \pi_1; \ba, \bd) - u(\pi, \pi_1)).
\end{equation}
First, consider the 2nd term, which involves randomness only from the
2nd stage. From \eqref{eq:nucon} note that $\sum_{l=1}^k a_l E_{\pi_l}
u(X; \ba, \bd) = u(\pi, \pi_1)$. Then from \eqref{eq:uvhat} we have
 \[
 \sqrt{n}(\hat{u}(\pi, \pi_1; \ba, \bd) - u(\pi, \pi_1)) =
  \sum_{l=1}^k a_l \sqrt{\frac{n}{n_l}} \frac{\sum_{i=1}^{n_l} (u(X_i^{(l)}; \ba, \bd) - E_{\pi_l}
u(X; \ba, \bd))}{\sqrt{n_l}}.
\]
Since $\Phi_l$ is polynomially ergodic of order $m$ and $E_{\pi_l} |u(X; \ba,
\bd)|^{2+\delta}$ is finite where $m > 1 + 2/\delta$, it follows that $\sum_{i=1}^{n_l}
(u(X_i^{(l)}; \ba, \bd) - E_{\pi_l} u(X; \ba, \bd))/\sqrt{n_l} \cd
N(0, \tau^2_l(\pi ; \ba, \bd))$ where $\tau^2_l(\pi ; \ba, \bd)$ is
defined in \eqref{eq:taul}. As $n_l/n \rightarrow s_l$ and the Markov
chains $\Phi_l$'s are independent, it follows that $
\sqrt{n}(\hat{u}(\pi, \pi_1; \ba, \bd) - u(\pi, \pi_1)) \cd N(0,
\tau^2(\pi ; \ba, \bd))$.

Now we consider the 1st term in the right hand side of \eqref{eq:huminu}. Letting $F(\bz) = \hat{u}(\pi, \pi_1; \ba, \bz)$, by Taylor series
expansion of $F$ about $\bd$ we have
\begin{equation}
  \label{eq:Ftayl}
  \sqrt{n}(F(\hatbd) - F(\bd)) = \sqrt{n} \nabla F(\bd)^{\top} (\hatbd - \bd) + \frac{\sqrt{n}}{2} (\hatbd - \bd)^{\top} \nabla^2 F(\bd^*) (\hatbd - \bd),
\end{equation}
where $\bd^*$ is between $\bd$ and $\hatbd$.  Simple calculations show that
\begin{equation}
  \label{eq:delfconv}
  [\nabla F(\bd)]_{j-1} = \sum_{l=1}^k \frac{a_l}{n_l} \sum_{i=1}^{n_l} \frac{a_j
  \nu_{j} (X_i^{(l)}) \nu(X_i^{(l)})}{(\sum_{s = 1}^k a_s \nu_{s} (X_i^{(l)})/d_s)^2 d_j^2} \cas [c(\pi; \ba, \bd)]_{j-1}
\end{equation}
where $[c(\pi; \ba, \bd)]_{j-1}$ is defined in \eqref{eq:cdef}. We
know that $n/N \rightarrow q$. Using
similar arguments as in \citet{buta:doss:2011}, it follows that
$\nabla^2 F(\bd^*)$ is bounded in probability. Thus from \eqref{eq:Ftayl} we have
\begin{align*}
  \sqrt{n}(F(\hatbd) - F(\bd)) &= \sqrt{\frac{n}{N}}  \nabla F(\bd)^{\top} \sqrt{N} (\hatbd - \bd) +\frac{1}{2\sqrt{N}} \sqrt{\frac{n}{N}}[\sqrt{N}(\hatbd - \bd)]^{\top} \nabla^2 F(\bd^*) [\sqrt{N}(\hatbd - \bd)]\\
&=\sqrt{q} c(\pi; \ba, \bd)^{\top} \sqrt{N} (\hatbd - \bd) + o_p(1) .
\end{align*}
Then Theorem~\ref{thm:elnc} (1) follow from \eqref{eq:huminu} and the independence of the
two stages of Markov chain sampling.

Next to prove Theorem~\ref{thm:elnc} (2), note that, we already have a
consistent BM estimator $\widehat{V}$ of $V$. From
\eqref{eq:delfconv}, we have $[\hat{c}(\pi; \ba, \bd)]_{j-1} = [\nabla
F(\bd)]_{j-1} \cas [c(\pi; \ba, \bd)]_{j-1}$. Applying mean value
theorem on $[\nabla F(\bd)]_{j-1}$ and the fact that $\nabla^2
F(\bd^*)$ is bounded in probability, it follows that $[\hat{c}(\pi;
\ba, \hatbd)]_{j-1} - [\hat{c}(\pi; \ba, \bd)]_{j-1} \cas 0$. Writing
$c(\pi; \ba, \bd)^{\top} V c(\pi; \ba, \bd)$ as
$\sum_{i=1}^{k-1}\sum_{j=1}^{k-1}c_i V_{ij} c_j$, it then follows that
$\hat{c}(\pi;
\ba, \hatbd)^{\top} \widehat{V} \hat{c}(\pi;
\ba, \hatbd) \cas c(\pi; \ba, \bd)^{\top} V c(\pi; \ba, \bd)$.

We now show $\hat{\tau}_l^2 (\pi ; \ba, \hatbd)$ is a
consistent estimator of $\tau_l^2 (\pi ; \ba, \bd)$ where $\tau^2_l$
and $\hat{\tau}^2_l$ are defined in \eqref{eq:taul} and
\eqref{eq:tauldef}, respectively. Since the Markov chains
$\{X_i^{(l)}\}_{i=1}^{n_l}$ are independent, it
then follows that $\tau^2 (\pi ; \ba, \bd)$ is consistently estimated
by $\hat{\tau}^2 (\pi ; \ba, \hatbd)$ completing the proof of
Theorem~\ref{thm:elnc} (2).

If $\bd$ is known from the assumptions of Theorem~\ref{thm:elnc} (2)
and the results in \citet{vats:fleg:jone:2015:output}, %(see also \citet{jone:hara:caff:neat:2006, bedn:latu:2007}), 
we know that $\tau_l^2 (\pi ; \ba, \bd)$ is
consistently estimated by its BM estimator $\hat{\tau}_l^2
(\pi ; \ba, \bd)$. Note that, $\hat{\tau}_l^2 (\pi ; \ba, \bd)$ is
defined in terms of the quantities $u(X_i^{(l)}; \ba, \bd)$'s. We now
show that $\hat{\tau}_l^2 (\pi ; \ba, \hatbd) - \hat{\tau}_l^2 (\pi ;
\ba, \bd) \cas 0.$

Denoting $\hat{\tau}_l^2 (\pi ; \ba, \bz)$ by $G(\bz)$, by the mean value theorem (in multiple variables), there exists $\bd^*=t\hatbd+(1-t)\bd$ for some $t\in(0,1)$, such that $G(\hatbd)-G(\bd) = \nabla G(\bd^*) \cdot (\hatbd-\bd)$.  For any $j \in \{2,\cdots, k\}$, and $\bz \in {R^+}^{k-1}$,
\begin{equation}\label{eq:grad}
\frac{\partial G(\bz)}{\partial z_j}
=\frac{b_l}{e_l - 1}\left[ \sum_{m=0}^{e_l - 1} 2 ( \bar{u}_{m}(\ba, \bz) - \bar{\bar{u}}(\ba,
  \bz))  \left(\frac{\partial \bar{u}_{m}(\ba, \bz)}{\partial z_j}-\frac{\partial \bar{\bar{u}}(\ba, \bz)}{\partial z_j} \right) \right]
\end{equation}
Let $\bar{W}_{m} := \bar{u}_{m}(\ba, \bz) - E_{\pi_l} (u(X;\ba, \bz))$
and $\bar{\bar{W}}:= \bar{\bar{u}}(\ba, \bz) - E_{\pi_l} (u(X;\ba, \bz))$.
Note that, there exists, $\sigma^2 <\infty$ such that
$\sqrt{b_l}\bar{W}_{m}\cd \text{N}(0,\sigma^2) $, and
$\sqrt{n_l}\bar{\bar{W}}\cd \text{N}(0,\sigma^2)$.
Simple calculations show that
\[
  \frac{\partial \bar{u}_{m}(\ba, \bz)}{\partial z_j}
  =\frac{a_j}{ z_j^2} \frac{1}{b_l}\sum_{i=m b_l +1}^{(m+1)b_l} \left[\frac{\nu(X_i^{(l)}) \nu_{j}(X_i^{(l)})}{\left(\sum_s a_s \nu_{s}(X_i^{(l)})/z_s \right)^2}.  \right]
\]
Hence, letting $\alpha_j=E_{\pi_l} [\nu(X) \nu_{j}(X)/\left(\sum_s a_s
  \nu_{s}(X)/z_s \right)^2]$, we write
\[ \frac{\partial \bar{u}_{m}(\ba, \bz)}{\partial z_j}-\frac{\partial \bar{\bar{u}}(\ba, \bz)}{\partial z_j}
  \equiv \frac{a_j}{ z_j^2} \left\{\bar{Z}_{m,j} \right\}
  -\frac{a_j}{ z_j^2}\left\{\bar{\bar{Z}}_{j}\right\},
\]
where $\bar{Z}_{1,j} = (1/b_l)\sum_{i=1}^{b_l} [\nu(X_i^{(l)})
\nu_{j}(X_i^{(l)})/\{\sum_s a_s \nu_{s}(X_i^{(l)})/z_s \}^2]
-\alpha_j$ and $\bar{\bar{Z}}_{j}$ is similarly defined.
Note that, there exists $\tau^2_j < \infty$, such that
 $\sqrt{b_l}
\bar{Z}_{m,j}\cd \text{N}(0,\tau_{j}^2)$,
 and
$\sqrt{n_l}\bar{\bar{Z}}_{j} \cd \text{N}(0,\tau_{j}^2)$.  From \eqref{eq:grad} we have
\[\begin{split}
\frac{\partial G(\bz)}{\partial z_j}= &\frac{a_j}{ z_j^2}\frac{2}{e_l - 1}\sum_{m=0}^{e_l - 1} \left[ \sqrt{b_l}(\bar{W}_{m}-\bar{\bar{W}})  \sqrt{b_l}\left( \bar{Z}_{m,j} -\bar{\bar{Z}}_{j} \right) \right]\\
=&\frac{a_j}{ z_j^2}\frac{2}{e_l - 1}\sum_{m=0}^{e_l - 1}
 \left[ \sqrt{b_l}\bar{W}_{m}\sqrt{b_l} \bar{Z}_{m,j}\right] \\
 - &  \frac{a_j}{ z_j^2} 2b_l\left[ \bar{\bar{Z}}_{j} \frac{1}{e_l - 1}\sum_{m=0}^{e_l - 1}\bar{W}_{m}
 +  \bar{\bar{W}}\frac{1}{e_l - 1}\sum_{m=0}^{e_l - 1}\bar{Z}_{m,j}
 - \frac{e_l}{e_l-1}\bar{\bar{W}} \bar{\bar{Z}}_{j}\right] \\
=&\frac{a_j}{ z_j^2}\frac{2}{e_l - 1}\sum_{m=0}^{e_l - 1}
 \left[ \sqrt{b_l}\bar{W}_{m}\sqrt{b_l} \bar{Z}_{m,j}\right]
 - \frac{a_j}{ z_j^2}\frac{2}{e_l - 1} \left[ \sqrt{n_l}\bar{\bar{W}} \sqrt{n_l} \bar{\bar{Z}}_{j}\right].
\end{split}\]
Then using similar arguments as in the proof of Theorem~\ref{thm:CLT}, it can be shown that
$\partial G(\bz)/\partial z_j$ is bounded with probability one. Then it follows that
 \begin{equation*}
 |G(\hatbd)-G(\bd)| \leq  \underset{1 \le j \le k-1}{\max} \left\{\left| \frac{\partial G(\bd^*)}{\partial z_j}\right|\right\} \|\hatbd-\bd\| \cas 0.
 \end{equation*}

\section{Proof of Theorem~\ref{thm:elex}}
\label{sec:appthm3}
As in the proof of Theorem~\ref{thm:elnc} we write
\begin{equation}
  \label{eq:hvminv}
  \sqrt{n}(\hat{\eta}^{[f]}(\pi; \ba, \hatbd) - E_\pi f) = \sqrt{n}(\hat{\eta}^{[f]}(\pi; \ba, \hatbd) - \hat{\eta}^{[f]}(\pi; \ba, \bd)) + \sqrt{n}(\hat{\eta}^{[f]}(\pi; \ba, \bd) - E_\pi f).
\end{equation}
First, consider the 2nd term, which involves randomness only from the
2nd stage. Since
\begin{equation*}
    \hat{v} \cas \sum_{l=1}^k a_l E_{\pi_l} v^{[f]}(X; \ba, \bd) =  \int_{\sX} \frac{f(x) \sum_{l=1}^k a_l
              \nu_l(x) / m_l} {\sum_{s=1}^k a_s \nu_s(x) / (m_s /
              m_1)} \nu(x) \, \mu(dx) = \frac{m}{m_1} E_\pi f,
\end{equation*}
we have
$\sum_{l=1}^k a_l E_{\pi_l}
v^{[f]}(X; \ba, \bd) = E_\pi f u(\pi, \pi_1)$. Then from \eqref{eq:uvhat} we have
\begin{equation}
  \label{eq:multclt}
  \sqrt{n} \left(\begin{array}{c}
    \hat{v}^{[f]}(\pi; \ba, \bd) - E_\pi f u(\pi, \pi_1)\\
    \hat{u}(\pi, \pi_1; \ba, \bd) - u(\pi, \pi_1)\\
\end{array}
\right) =
  \sum_{l=1}^k a_l \sqrt{\frac{n}{n_l}} \frac{1}{\sqrt{n_l}} \sum_{i=1}^{n_l} \left(\begin{array}{c}
     v^{[f]}(X_i^{(l)}; \ba, \bd) - E_{\pi_l}
v^{[f]}(X; \ba, \bd) \\
    u(X_i^{(l)}; \ba, \bd) - E_{\pi_l}
u(X; \ba, \bd)\\
\end{array}
\right) .
\end{equation}
From the conditions of Theorem~\ref{thm:elex} and the fact that the
Markov chains $\Phi_l, l=1,\dots,k$ are independent, it follows that
the above vector \eqref{eq:multclt} converges in distribution to the
bivariate normal distribution with mean $0$ and covariance matrix
$\Gamma(\pi; \ba, \bd)$ defined in \eqref{eq:defgam}. Then applying
the Delta method to the function $g(x, y) = x/y$ we have a CLT for the
ratio estimator $\hat{\eta}^{[f]}(\pi; \ba, \bd)$, that is, we have
$\sqrt{n}(\hat{\eta}^{[f]}(\pi; \ba, \bd) - E_\pi f) \cd N(0,
\rho(\pi; \ba, \bd))$ where $\rho(\pi; \ba, \bd))$ is defined in \eqref{eq:defrho}.

Next letting $L(\bz) = \hat{\eta}^{[f]}(\pi; \ba, \bz)$, by Taylor series
expansion of $L$ about $\bd$ we have
\begin{equation}
  \label{eq:Ltayl}
  \sqrt{n}(L(\hatbd) - L(\bd)) = \sqrt{n} \nabla L(\bd)^{\top} (\hatbd - \bd) + \frac{\sqrt{n}}{2} (\hatbd - \bd)^{\top} \nabla^2 L(\bd^*) (\hatbd - \bd),
\end{equation}
where $\bd^*$ is between $\bd$ and $\hatbd$.  Simple calculations show that
\begin{equation}
  \label{eq:dellconv}
  [\nabla L(\bd)]_{j-1} = [\hat{e}(\pi; \ba, \bd)]_{j-1} \cas [e(\pi; \ba, \bd)]_{j-1}
\end{equation}
where $[e(\pi; \ba, \bd)]_{j-1}$ and $[\hat{e}(\pi; \ba, \bd)]_{j-1}$ are defined in \eqref{eq:defe} and
\eqref{eq:defhate} respectively. It can be shown that
$\nabla^2 L(\bd^*)$ is bounded in probability. Thus from \eqref{eq:Ltayl} we have $\sqrt{n}(L(\hatbd) - L(\bd)) = \sqrt{q} e(\pi; \ba, \bd)^{\top} \sqrt{N} (\hatbd - \bd) + o_p(1)$.  Then Theorem~\ref{thm:elex} (1) follow from \eqref{eq:hvminv} and the independence of the
two stages of Markov chain sampling.

Next to prove Theorem~\ref{thm:elex} (2), note that, we already know that $\widehat{V}$ is a
consistent BM estimator of $V$. From
\eqref{eq:dellconv}, we have $[\hat{e}(\pi; \ba, \bd)]_{j-1} \cas [e(\pi; \ba, \bd)]_{j-1}$. Applying mean value
theorem on $[\nabla L(\bd)]_{j-1}$ and the fact that $\nabla^2
L(\bd^*)$ is bounded in probability, it follows that $[\hat{e}(\pi;
\ba, \hatbd)]_{j-1} - [\hat{e}(\pi; \ba, \bd)]_{j-1} \cas 0$.

From \eqref{eq:elncclt} we know that $\hat{u}(\pi, \pi_1; \ba, \hatbd)
\cas u(\pi, \pi_1)$.  From \eqref{eq:elexclt} we know
$\hat{\eta}^{[f]}(\pi; \ba, \hatbd) \cas E_\pi f$. Since $
\hat{v}^{[f]} (\pi, \pi_1; \ba, \bd) = \hat{\eta}^{[f]} (\pi; \ba,
\bd) \hat{u}(\pi, \pi_1; \ba, \bd)$, it follows that $\hat{v}^{[f]}
(\pi, \pi_1; \ba, \hatbd) \cas E_\pi f u(\pi, \pi_1)$. Thus $\nabla
h(\hat{v}^{[f]}
(\pi, \pi_1; \ba, \hatbd), \hat{u}(\pi, \pi_1; \ba, \hatbd)) \cas \nabla h(E_\pi f u(\pi,
\pi_1), u(\pi, \pi_1))$. Thus to prove Theorem~\ref{thm:elex} (2), we only need to show that
$\widehat{\Gamma}_l(\pi; \ba, \hatbd) \cas \Gamma_l(\pi; \ba, \bd)$.

If $\bd$ is known from the assumptions of Theorem~\ref{thm:elex} (2)
and the results in \citet{vats:fleg:jone:2015:output}, we know that
$\Gamma_l (\pi ; \ba, \bd)$ is consistently estimated by its BM estimator $\widehat{\Gamma}_l (\pi ; \ba, \bd)$.  We now show
that $\widehat{\Gamma}_l (\pi ; \ba, \hatbd) - \widehat{\Gamma}_l
(\pi ; \ba, \bd) \cas 0.$

From Theorem~\ref{thm:elnc} (2), we know that $\hat{\gamma}_l^{22} (\pi
; \ba, \hatbd) - \hat{\gamma}_l^{22} (\pi ; \ba, \bd) \cas 0$. We now
show $\hat{\gamma}_l^{11} (\pi ; \ba, \hatbd) - \hat{\gamma}_l^{11}
(\pi ; \ba, \bd)\cas 0$.  Letting $\hat{\gamma}_l^{11} (\pi ; \ba, \bz)$ by $H(\bz)$,
 by the mean value theorem, there exists
$\bd^*=t\hatbd+(1-t)\bd$ for some $t\in(0,1)$, such that $H(\hatbd)-H(\bd) = \nabla H(\bd^*) \cdot (\hatbd-\bd)$.  For any $j \in \{2,\cdots, k\}$, and $\bz \in {R^+}^{k-1}$,
\begin{equation*}
\frac{\partial H(\bz)}{\partial z_j}
=\frac{b_l}{e_l - 1}\left[ \sum_{m=0}^{e_l - 1} 2 ( \bar{v}^{[f]}_{m}(\ba, \bz) - \bar{\bar{v}}^{[f]}(\ba,
  \bz))  \left(\frac{\partial \bar{v}_{m}^{[f]}(\ba, \bz)}{\partial z_j}-\frac{\partial \bar{\bar{v}}^{[f]}(\ba, \bz)}{\partial z_j} \right) \right].
\end{equation*}
Let $\bar{W}^{[f]}_{m} := \bar{v}^{[f]}_{m}(\ba, \bz) - E_{\pi_l} (v^{[f]}(X;\ba, \bz))$
and $\bar{\bar{W}}^{[f]}:= \bar{\bar{v}}^{[f]}(\ba, \bz) - E_{\pi_l} (v^{[f]}(X;\ba, \bz))$.
Note that, there exists, $\sigma^2_f <\infty$ such that
$\sqrt{b_l}\bar{W}^{[f]}_{m}\cd \text{N}(0,\sigma^2_f) $, and
$\sqrt{n_l}\bar{\bar{W}}^{[f]} \cd \text{N}(0,\sigma^2_f)$.
Simple calculations show that
\[
  \frac{\partial \bar{v}^{[f]}_{m}(\ba, \bz)}{\partial z_j}
  =\frac{a_j}{ z_j^2} \frac{1}{b_l}\sum_{i=m b_l +1}^{(m+1)b_l} \left[\frac{f(X_i^{(l)}) \nu(X_i^{(l)}) \nu_{j}(X_i^{(l)})}{\left(\sum_s a_s \nu_{s}(X_i^{(l)})/z_s \right)^2}.  \right]
\]
Hence, letting $\alpha^{[f]}_j=E_{\pi_l} [f(X) \nu(X) \nu_{j}(X)/\left(\sum_s a_s
  \nu_{s}(X)/z_s \right)^2]$, we write
\[ \frac{\partial \bar{v}^{[f]}_{m}(\ba, \bz)}{\partial z_j}-\frac{\partial \bar{\bar{u}}^{[f]}(\ba, \bz)}{\partial z_j}
  \equiv \frac{a_j}{ z_j^2} \left\{\bar{Z}^{[f]}_{m,j} \right\}
  -\frac{a_j}{ z_j^2}\left\{\bar{\bar{Z}}^{[f]}_{j}\right\},
\]
where $\bar{Z}^{[f]}_{1,j} = (1/b_l)\sum_{i=1}^{b_l} [ f(X_i^{(l)}) \nu(X_i^{(l)})
\nu_{j}(X_i^{(l)})/\{\sum_s a_s \nu_{s}(X_i^{(l)})/z_s \}^2]
-\alpha_j^{[f]}$ and $\bar{\bar{Z}}^{[f]}_{j}$ is similarly defined.
Note that, there exists $\tau^2_{j,f} < \infty$, such that
 $\sqrt{b_l}
\bar{Z}_{m,j}\cd \text{N}(0,\tau_{j, f}^2)$,
 and
$\sqrt{n_l}\bar{\bar{Z}}_{j} \cd \text{N}(0,\tau_{j, f}^2)$.  The rest of the proof is analogous to Theorem~\ref{thm:elnc}, in that we have
\[
\frac{\partial H(\bz)}{\partial z_j} = \frac{a_j}{ z_j^2}\frac{2}{e_l - 1}\sum_{m=0}^{e_l - 1}
 \left[ \sqrt{b_l}\bar{W}^{[f]}_{m}\sqrt{b_l} \bar{Z}^{[f]}_{m,j}\right]
 - \frac{a_j}{ z_j^2}\frac{2}{e_l - 1} \left[ \sqrt{n_l}\bar{\bar{W}}^{[f]} \sqrt{n_l} \bar{\bar{Z}}^{[f]}_{j}\right].
\]
Then it can be shown $\hat{\gamma}_l^{11} (\pi ; \ba, \hatbd) - \hat{\gamma}_l^{11} (\pi ; \ba, \bd)\cas 0$ and finally $\hat{\gamma}_l^{12} (\pi ; \ba, \hatbd) - \hat{\gamma}_l^{12} (\pi ; \ba, \bd)\cas 0$.

\section{Regeneration with general weights}
\label{sec:regen}
\citet{tan:doss:hobe:2015} provide a regeneration based central limit
theorem (CLT) for the estimators $\hat{\eta}$ and $\hu$ defined in 1.3
and 3.1 respectively in the main text. In the case
when $\bd$ is unknown, they allow only a special choice for the weight
vector, namely $\ba =(1, \hatbd)$ for their results to hold, where
$\hatbd$ is the estimator of $\bd$ based on the Stage 1 chains
discussed in Section~2 of the main text. In this section, we establish a regeneration based CLT for
$\hat{\eta}$ and $\hu$ with any choice of the weight vector $\ba$.

We will refer to the following conditions.
  \begin{enumerate}
  \item[A1] For each $l = 1, \ldots, k$, the Markov chain $\Phi_l = \{ X_0^{(l)}, X_1^{(l)}, \ldots
    \}$ is geometrically ergodic and has $\pi_l$ as its invariant
    density.
  \item[A2] Let $k_l : \sX \times \sX \rightarrow [0, \infty )$ be the Markov transition density for
    $\Phi_l$, so that for any measurable set $A$ we have $P\bigl(
    X_{n+1}^{(l)} \in A \mid X_n^{(l)} = x \bigr)=
    \int_A k_l(y |x) \mu (dy)$. Suppose that for each $l= 1, \ldots, k$, $k_l$
    satisfies the following {\it minorization condition}: 
\begin{equation}
  \label{eq:mino}
  k_l (y|x) \geq s_l(x) \, q_l(y) \qquad \text{for all } x, y \in \sX,
\end{equation}
where the function
    $s_l \colon \sX \rightarrow [0, 1)$ with $E_{\pi_l} s_l > 0$, and
     $q_l$ is a probability density function on $\sX$.
  \item[A3] Recall the functions $u(X; \ba,
    \bd)$ and $v^{[f]}(X; \ba,
    \bd)$ defined in (3.2) of our paper. There exists $\epsilon > 0$ such that $E_{\pi_l} |v^{[f]}(X; \ba,
    \bd)|^{2+\epsilon}$ and $E_{\pi_l} |u(X; \ba,
    \bd)|^{2+\epsilon}$ are finite.
  \item[A4] Suppose $\Phi_l$ is simulated for $R_l$ regenerative tours
for $l = 1, \ldots, k$. Assume $R_l / R_1 \rightarrow b_l \in(0, \infty)$ as $R_1
    \rightarrow \infty$.
  \end{enumerate}

%\newpage
%\lhead[\footnotesize\thepage\fancyplain{}\leftmark]{}\rhead[]{\fancyplain{}\rightmark\footnotesize\thepage}%Put this line in Page 2

  Following \citet{tan:doss:hobe:2015}, let the \textsl{regeneration
    times} for the $l^{\text{th}}$ Markov chain be $\tau_0^{(l)} = 0, \tau_1^{(l)}, \tau_2^{(l)}, \ldots$.
  Accordingly, the chain $\Phi_l$ is broken up into ``tours'' $\bigl\{
  \bigl( X_{\tau_{t-1}^{(l)}}, \ldots, X_{\tau_t^{(l)}-1} \bigr), \, t
  = 1, 2, \ldots \bigr\}$ that are independent stochastic replicas of
  each other. Suppose we simulate $R_l$ tours of the $l^{\text{th}}$ Markov chain
for $l = 1, \ldots, k$, so the length of the $l^{\text{th}}$ chain is $n_l = \tau_{R_l}^{(l)}$.
Also as in \citet{tan:doss:hobe:2015}, for
$t = 1, 2, \ldots, R_l$ define
\begin{equation}
  \label{eq:VU}
  V_t^{(l)} = \sum_{i=\tau_{t-1}^{(l)}}^{\tau_t^{(l)}-1}
  v^{[f]}(X_i^{(l)}; \ba, \bd), \quad U_t^{(l)} =
  \sum_{i=\tau_{t-1}^{(l)}}^{\tau_t^{(l)}-1} u(X_i^{(l)}; \ba, \bd),
  \quad \text{and} \quad T_t^{(l)} =
  \sum_{i=\tau_{t-1}^{(l)}}^{\tau_t^{(l)}-1} 1 = \tau_t^{(l)} -
  \tau_{t-1}^{(l)},
\end{equation}
where the sums range over the values of $i$ that constitute the
$t^{\text{th}}$ tour.

 Recall from Remark 4 in Section 3 of our paper, when $\bd$ is unknown, we set $\ba = \bw *(1, \hatbd)$ where $*$ denotes component-wise multiplication.  That is, $(a_1, \ldots, a_k) = (w_1, w_2, \ldots, w_k)*(1, \hat{d}_2, \ldots, \hat{d}_k) $ for any pre-determined weight $\bw$. With this choice, the expressions for $u$ and $v^{[f]}$ in (3.2) become
\begin{equation}
  \label{eq:uv-new}
  u\bigl( x; \bw*(1, \hatbd), \hatbd \, \bigr) = \frac{\nu(x)}
  {\sum_{l=1}^k w_l\nu_l(x)} \quad \text{and} \quad v^{[f]}\bigl( x; \bw*(1, \hatbd), \hatbd \, \bigr) = \frac{f(x) \nu(x)} {\sum_{l=1}^k
  w_l\nu_l(x)}.
\end{equation}
The above quantities do not involve $\hatbd$, and consequently for each $l$, the
triples $\bigl( V_t^{(l)}, U_t^{(l)}, T_t^{(l)} \bigr), \, t = 0, 1,
2, \ldots$ defined in \eqref{eq:VU} are iid, and we have
independence across $l$'s.  The estimator for $\eta$ reduces to
\begin{align}
  \label{eq:Ad}
  \hat{\eta} = \hat{\eta}_{N,n} \bigl( \bw*(1, \hatbd), \hatbd \,
  \bigr) & = \sum_{l=1}^k \frac{
  w_l \hat{d}_l }{n_l} \sum_{i=1}^{n_l}
             \frac{ f(X_i^{(l)}) \nu(X_i^{(l)}) } { \sum_{s=1}^k
             w_l\nu_s(X_i^{(l)}) } \bigg/ \sum_{l=1}^k \frac{w_l \hat{d}_l}
             {n_l} \sum_{i=1}^{n_l} \frac{ \nu(X_i^{(l)})  } {
             \sum_{s=1}^k w_s \nu_s(X_i^{(l)}) } \\[2mm]
         & = \sum_{l=1}^k \frac{w_l\hat{d}_l} {n_l}
             \sum_{t=1}^{R_l}V_t^{(l)} \bigg/ \sum_{l=1}^k
             \frac{w_l\hat{d}_l} {n_l} \sum_{t=1}^{R_l} U_t^{(l)}
             \nonumber \\[2mm]
         & = \sum_{l=1}^k w_l\hat{d_l} \frac{\bar{V}^{(l)}}
             {\bar{T}^{(l)}} \bigg/ \sum_{l=1}^k w_l\hat{d_l}
             \frac{\bar{U}^{(l)}} {\bar{T}^{(l)}}, \nonumber
\end{align}
where
\begin{equation*}
  U_t^{(l)} = \sum_{i=\tau_{t-1}^{(l)}}^{\tau_t^{(l)}-1}
  \frac{\nu(X_i^{(l)})} {\sum_{s=1}^k w_s\nu_s(X_i^{(l)})} \quad
  \text{and} \quad V_t^{(l)} =
  \sum_{i=\tau_{t-1}^{(l)}}^{\tau_t^{(l)}-1} \frac{f(X_i^{(l)})
  \nu(X_i^{(l)})} {\sum_{s=1}^k w_s\nu_s(X_i^{(l)})},
\end{equation*}
$\bar{T}^{(l)} = R_l^{-1} \sum_{t=1}^{R_l} T_t^{(l)}$ be the average
tour length and, analogously, $\bar{V}^{(l)} = R_l^{-1}
\sum_{t=1}^{R_l} V_t^{(l)}$ and $\bar{U}^{(l)} = R_l^{-1}
\sum_{t=1}^{R_l} U_t^{(l)}$.
Similarly, the estimator for $m/m_1$ reduces to
\begin{equation}
  \label{eq:Bd}
\begin{split}
  \hat{u} & = \hat{u}_{N,n} \bigl( \bw*(1, \hatbd), \hatbd \, \bigr) =
  \sum_{l=1}^k \frac{w_l\hat{d}_l}{n_l} \sum_{i=1}^{n_l}
  \frac{\nu(X_i^{(l)})} {\sum_{s=1}^k w_s\nu_s(X_i^{(l)})} \\
  & = \sum_{l=1}^k \frac{w_l\hat{d}_l}{n_l} \sum_{t=1}^{R_l} U_t^{(l)} =
  \sum_{l=1}^k w_l\hat{d_l} \frac{\bar{U}^{(l)}} {\bar{T}^{(l)}}.
\end{split}
\end{equation}

Theorem 4 below gives the asymptotic distributions of $\hat{\eta}$ and $\hat{u}$. It extends \pcite{tan:doss:hobe:2015} Theorem 2 to the general choice of weight vector $\ba$. To state the theorem, we first need to define some notation.  Let $\tilde{M}$ and $\tilde{L}$ be the vectors of length $k - 1$ for which the $(j - 1)^{\text{th}}$ coordinates are, for $j = 2, \dots, k$,
\begin{equation}
  \label{eq:L-new}
\begin{split}
  \tilde{M}_{j-1} & = w_j E_{\pi_j} u \text{ and } \\
  \tilde{L}_{j-1} & =
  \frac{w_j E_{\pi_j} v^{[f]}} {\sum_{l=1}^k w_l d_l E_{\pi_l} u} - \frac{\bigl(
  \sum_{l=1}^k w_l d_l E_{\pi_l} v^{[f]} \bigr) \bigl( w_j E_{\pi_j} u \bigr)}
  {\bigl( \sum_{l=1}^k w_l d_l E_{\pi_l} u \bigr)^2}.
\end{split}
\end{equation}
As in \citet{tan:doss:hobe:2015}, assume that in
Stage~$1$, for $l = 1, \ldots, k$, the $l^{\text{th}}$ chain has been run for
$\rho_l$ regenerations.  So the length of the $l^{\text{th}}$ chain,
$N_l = T_1^{(l)} + \ldots + T_{\rho_l}^{(l)}$, is random.  We assume
that $\rho_1, \ldots, \rho_k \rightarrow \infty$ in such a way that
$\rho_l / \rho_1 \rightarrow c_l \in (0, \infty)$, for $l = 1,
\ldots, k$.

\medskip

\noindent {\bf Theorem 4}
  %\label{thm:d-unknown}
  {\it Suppose that for the Stage~$1$ chains, conditions~A1 and~A2
  hold, and that for the Stage~$2$ chains,
  conditions~A1--A4 hold.  If $\rho_1 \rightarrow \infty$ and $R_1
  \rightarrow \infty$ in such a way that $R_1 / \rho_1 \rightarrow q
  \in [0, \infty)$, then
  \begin{equation*}
    R_1^{1/2} \bigl( \hat{u} - m/m_1 \bigr) \cd {\cal N} \bigl( 0,
    q \tilde{M}^{\top} W \tilde{M} + \kappa^2 \bigr)
  \end{equation*}
  and
  \begin{equation*}
    R_1^{1/2} \bigl( \hat{\eta} - \eta \bigr) \cd {\cal N} \bigl(
    0, q \tilde{L}^{\top} W \tilde{L} + \tau^2 \bigr),
  \end{equation*}
  
  \noindent where $\tilde{M}$, $\tilde{L}$ are given in
  equations~\eqref{eq:L-new}, $W$, $\kappa^2$ and $\tau^2$ are given in
  equations (2.15), (2.8),
  and (2.10) of \citet{tan:doss:hobe:2015}, respectively.  In their (2.8)
  and (2.10), $\ba$ is taken to be $\ba = \bw * (1, \bd)$.
  Furthermore, we can form strongly consistent estimates of the
  asymptotic variances if we use $\widehat{W}$, $\hat{\kappa}^2$,
  and $\hat{\tau}^2$ defined in (2.16)
  and (2.11) of \citet{tan:doss:hobe:2015}, respectively, and use the standard
  empirical estimates of $\tilde{M}$ and $\tilde{L}$.}
%\end{theorem}

\subsection{Proof of Theorem 4}

We first prove the CLT for $\hat{\eta}$.  Note that
\begin{equation}
  \label{eq:etad}
  R_1^{1/2} \bigl[ \hat{\eta} \bigl( \bw * (1, \hatbd), \hatbd \, \bigr) -
  \eta \bigr] = R_1^{1/2} \bigl[ \hat{\eta} \bigl( \bw * (1, \hatbd),
  \hatbd \, \bigr) - \hat{\eta} \bigl( \bw * (1, \bd), \bd \bigr) \bigr] +
  R_1^{1/2} \bigl[ \hat{\eta} \bigl( \bw * (1, \bd), \bd \bigr) - \eta
  \bigr].
\end{equation}

The second term on the right side of~\eqref{eq:etad} involves
randomness coming only from Stage~$2$ sampling, and its distribution
is given by Theorem 1 of \citet{tan:doss:hobe:2015}: it is
asymptotically normal with mean $0$ and variance
$\tau^2$.  The first term involves randomness from both Stage~$1$ and
Stage~$2$ sampling.  However, as in the proofs of Theorem 2 and 3, we
can show that for this term, the randomness from Stage~$2$ is
asymptotically negligible, so that only Stage~$1$ sampling contributes
to its asymptotic distribution. Finally, the asymptotic normality of the left side
of~\eqref{eq:etad} follows since the two stages of
sampling are independent. We now provide the details of the proof.

Consider the first term on the right side of~\eqref{eq:etad}.
Recall that if $\ba = \bw * (1, \bd)$, then
\begin{equation*}
  v^{[f]}(x) := v^{[f]}(x; \ba, \bd) = \frac{f(x) \nu(x)} {\sum_{l=1}^k
  w_l \nu_l(x)} \quad \text{and} \quad u(x) := u(x; \ba, \bd) =
  \frac{\nu(x)} {\sum_{l=1}^k w_l \nu_l(x)}.
\end{equation*}
With~\eqref{eq:Ad} and~\eqref{eq:Bd} in mind, define the function
\begin{equation*}
 A(\bz) = \hat{\eta} \bigl( \bw*(1, \bz), \bz \bigr) = \sum_{l=1}^k
  \frac{w_l z_l}{n_l} \sum_{i=1}^{n_l} v^{[f]}(X_i^{(l)}) \bigg/ \sum_{l=1}^k
  \frac{w_l z_l}{n_l} \sum_{i=1}^{n_l} u(X_i^{(l)})
\end{equation*}
for $\bz = (z_2, \ldots, z_k)^{\top}$, with $z_l > 0$ for $l = 2,
\ldots, k$, and $z_1 = 1$.  Note that setting $\bz = \bd$ gives
$A(\bd) = \hat{\eta}(\bw*(1, \bd), \bd)$, and setting $\bz = \hatbd$
gives $A(\hatbd) = \hat{\eta}(\bw*(1, \hatbd), \hatbd)$.

 By a Taylor series expansion of $A$ about $\bd$ we get
\begin{multline*}
  R_1^{1/2} \bigl[ \hat{\eta} \bigl( \bw*(1, \hatbd), \hatbd \, \bigr) -
  \hat{\eta} \bigl( \bw*(1, \bd), \bd \bigr) \bigr] = R_1^{1/2} \nabla
  A({\bd})^{\top} ({\hatbd} - \bd) + \frac{R_1^{1/2}} {2}
  ({\hatbd} - \bd)^{\top} \nabla^2 A(\bd^*) (\hatbd - \bd) \\
  = R_1^{1/2} \nabla A({\bd})^{\top} ({\hatbd} - \bd) +
  \frac{R_1^{1/2}} {2\rho_1} \bigl( \rho_1^{1/2} ({\hatbd} - \bd)
  \bigr)^{\top} \nabla^2 A(\bd^*) \bigl( \rho_1^{1/2} (\hatbd - \bd)
  \bigr),
\end{multline*}
where $\bd^{*}$ is between ${\bd}$ and ${\hatbd}$.  As $R_1
\rightarrow \infty$, $n_l \rightarrow \infty $ for each $l$.  We
first show that the gradient $\nabla A({\bd})$ converges almost
surely to a finite constant vector by proving that each one of its
components, $[A(\bd)]_{j-1}, \, j = 2, \ldots, k$, converges almost
surely as $R_1 \rightarrow \infty$.  As $n_l \rightarrow \infty$ for
$l = 1, \ldots, k$, for $j = 2, \ldots, k$, we have
\begin{align*}
  [\nabla
  A(\bd)]_{j-1} & = \frac{ (w_j/n_j) \sum_{i=1}^{n_j} v^{[f]}(X_i^{(j)}) }
                    {\sum_{l=1}^k (w_l d_l / n_l) \sum_{i=1}^{n_l}
                    u(X_i^{(l)})} \\
                    & - \frac{\bigl( \sum_{l=1}^k (w_ld_l /
                    n_l) \sum_{i=1}^{n_l} v^{[f]}(X_i^{(l)}) \bigr) \bigl(
                    (w_j/n_j) \sum_{i=1}^{n_j} u(X_i^{(j)}) \bigr)}
                    {\bigl( \sum_{l=1}^k (w_l d_l / n_l)
                    \sum_{i=1}^{n_l} u(X_i^{(l)}) \bigr)^2}
                    \notag \\
                & \cas \frac{w_j E_{\pi_j} v^{[f]}}
                    {\sum_{l=1}^k w_l d_l E_{\pi_l} u} - \frac{\bigl(
                    \sum_{l=1}^k w_l d_l E_{\pi_l} v^{[f]} \bigr) \bigl(
                    w_j E_{\pi_j} u \bigr)} {\bigl( \sum_{l=1}^k w_l d_l
                    E_{\pi_l} u \bigr)^2}. \notag
\end{align*}
The expression above corresponds to $\tilde{L}_{j-1}$, which
is defined in~\eqref{eq:L-new}, and it is finite by assumption~$A3$.
Next, we show that the random Hessian matrix $\nabla^{2} A(\bd^*)$ is
bounded in probability, i.e., each element of this matrix is $O_p(1)$.
As $n_l \rightarrow \infty$ for $l = 1, \ldots, k$, for any $j, t \in
\{ 2, \ldots, k \}, j \neq t$, we have
\begin{align*}
  & [\nabla ^2 A (\bd^*)]_{t-1,j-1} = -\frac{\bigl( \frac{w_j}{n_j}
            \sum_{i=1}^{n_j} v^{[f]}(X_i^{(j)}) \bigr) \bigl(
            \frac{w_t}{n_t} \sum_{i=1}^{n_t} u(X_i^{(t)}) \bigr)}
            {\bigl( \sum_{l=1}^k \frac{w_l d_l^*}{n_l} \sum_{i=1}^{n_l}
            u(X_i^{(l)}) \bigr)^2 } \\[2mm]
          & - \biggl( \frac{w_j}{n_j} \sum_{i=1}^{n_j} u(X_i^{(j)})
            \biggr) \Biggl[ \frac{\frac{w_t}{n_t} \sum_{i=1}^{n_t}
            v^{[f]}(X_i^{(t)})} {\bigl( \sum_{l=1}^k \frac{w_l d_l^*}{n_l}
            \sum_{i=1}^{n_l} u(X_i^{(l)}) \bigr)^2} - 2 \frac{\bigl(
            \sum_{l=1}^k \frac{w_l d_l^*}{n_l} \sum_{i=1}^{n_l}
            v^{[f]}(X_i^{(l)}) \bigr) \bigl( \frac{w_t}{n_t}
            \sum_{i=1}^{n_t} u(X_i^{(t)}) \bigr)} {\bigl(
            \sum_{l=1}^k \frac{w_l d_l^*}{n_l} \sum_{i=1}^{n_l}
            u(X_i^{(l)}) \bigr)^3} \Biggr] \\[2mm]
          & \hspace{-2mm} \cas - \frac{(w_j E_{\pi_j} v^{[f]}) (w_t E_{\pi_t} u)}
            {\bigl( \sum_{l=1}^k w_l d_l E_{\pi_l} u \bigr)^2} -
            (w_j E_{\pi_j} u) \Biggl[ \frac{w_t E_{\pi_t} v^{[f]}} {\bigl(
            \sum_{l=1}^k w_l d_l E_{\pi_l} u \bigr)^2} - 2 \frac{\bigl(
            \sum_{l=1}^k w_l d_l E_{\pi_l} v^{[f]} \bigr) (w_t E_{\pi_t} u )}
            {\bigl( \sum_{l=1}^k w_l d_l E_{\pi_l} u \bigr)^3} \Biggr],
\end{align*}
where the limits are also finite.

Now, we can rewrite~\eqref{eq:etad} as
\begin{align*}
  R_1^{1/2} \bigl[ \hat{\eta} \bigl( \bw * (1, \hatbd), \hatbd \, \bigr) -
  \eta \bigr] & = (R_1 / \rho_1)^{1/2} \nabla A(\bd)^{\top}
                  \rho_1^{1/2} (\hatbd - \bd) + R_1^{1/2} \bigl[
                  \hat{\eta} \bigl(\bw * (1, \bd), \bd \bigr) - \eta
                  \bigr] \\
              & \hspace{8mm} + \frac{1}{2 \rho_1^{1/2}} (R_1 /
                  \rho_1)^{1/2} \bigl[ \rho_1^{1/2} (\hatbd - \bd)
                  \bigr]^{\top} \nabla^{2} A(\bd^*) \bigl[
                  \rho_1^{1/2} (\hatbd - \bd) \bigr] \\
              & = q^{1/2} [\nabla A(\bd)]^{\top} \rho_1^{1/2}
                  (\hatbd - \bd) + R_1^{1/2} \bigl[ \hat{\eta}
                  \bigl( \bw * (1, \bd), \bd \bigr) - \eta \bigr] +
                  o_p(1).
\end{align*}
Since from \citet{tan:doss:hobe:2015} we have $\rho_1^{1/2} (\hatbd -
\bd) \cd {\cal N} (0, W)$ and the two sampling stages are assumed to
be independent, we conclude that
\begin{equation*}
  R_1^{1/2} \bigl[ \hat{\eta} \bigl( \bw * (1, \hatbd), \hatbd \, \bigr) -
  \eta \bigr] \cd {\cal N} \bigl( 0, q \tilde{L}^{\top} W
  \tilde{L} + \tau^2 \bigr).
\end{equation*}

The proof of the CLT for $\hat{u}$ is similar.  As in~\eqref{eq:etad}, we have
\begin{equation}
  \label{eq:ud}
\begin{split}
  R_1^{1/2} \bigl[ \hat{u} \bigl( \bw * (1, \hatbd), \hatbd \, \bigr) -
  m/m_1 \bigr] & = R_1^{1/2} \bigl[ \hat{u} \bigl( \bw * (1, \hatbd), \hatbd
  \, \bigr) - \hat{u} \bigl( \bw * (1, \bd), \bd \bigr) \bigr] \\
   & + R_1^{1/2}
  \bigl[ \hat{u} \bigl( \bw * (1, \bd), \bd \bigr) - m/m_1 \bigr].
\end{split}
\end{equation}
The asymptotic distribution of the second term in~\eqref{eq:ud} is
given in \pcite{tan:doss:hobe:2015} Theorem 1.  The first term is linear in
$\hatbd - \bd$:
\begin{equation}
  \label{eq:uhat-diff}
  \hat{u} \bigl( \bw*(1, \hatbd), \hatbd \, \bigr) - \hat{u} \bigl(  \bw*(1,
  \bd), \bd \bigr) = \sum_{j=2}^k w_j \biggl( \frac{1}{n_j}
  \sum_{i=1}^{n_j} u(X_i^{(j)}) \biggr) (\hat{d}_j - d_j).
\end{equation}
For $j = 2, \ldots, k$, the coefficient of $(\hat{d}_{j} - d_{j})$
in~\eqref{eq:uhat-diff} converges almost surely to $w_j E_{\pi_{j}} u$,
which is the term $\tilde{M}_{j-1}$ defined in~\eqref{eq:L-new}. 

  Finally,
from the independence of the two terms in~\eqref{eq:ud} we conclude
that
\begin{equation*}
  R_1^{1/2} \bigl[ \hat{u} \bigl( (1, \hatbd), \hatbd \, \bigr) -
  m/m_1 \bigr] \cd {\cal N} \bigl( 0, q \tilde{M}^{\top} W \tilde{M} + \kappa^2
  \bigr).
\end{equation*}

\par

%\begin{equation}
%\mbox {The 1st display equation of section 1.}
%\end{equation}

%\begin{equation}
%\mbox {The 2nd display equation of section 1.}
%\end{equation}

%\section{Title of section 2}
%\setcounter{equation}{0}

\section{Toy example}
\label{sec:toyapp}
\setcounter{equation}{0}

In this section, we follow up on the simulation studies that involve t distributions from Section~4 of the main text to verify Theorems~1-3. We also discuss different weights in forming generalized IS estimators and their effects on estimates of expectations and ratios of normalizing constants.

Let $t_{r,\mu}$ denote the t-distribution with degree of freedom $r$ and central parameter $\mu$. We consider $\pi_1(\cdot)$ and $\pi_2(\cdot)$ as the density functions for a $t_{5, \mu_1=1}$ and $t_{5, \mu_2=0}$, respectively. For simplicity, let $\nu_i(\cdot)=\pi_i(\cdot)$ for $i=1, 2$.  Our plan is to first estimate the ratio between the two normalizing constants, $d=m_2/m_1$.  Then we will study a sea of $t$-distributions $\Pi=\{t_{5,\mu}: \mu \in M\}$ where $M$ is a fine grid over $[0,1]$, say $M=\{0,.01,\cdots,.99, 1\}$. For each $\mu \in M$, we assume that $\nu_\mu(\cdot)=\pi_\mu(\cdot)$ and we estimate the ratio between its normalizing constant and $m_1$, denoted by $d_\mu:=\frac{m_\mu}{m_1}$. We also estimate the expectation of each distribution in $\Pi$, denoted $\E_{t_{5,\mu}}X$ or $\E_\mu X$ for short.  Clearly, the exact answers are $d=d_\mu=1$ and $\E_{\mu}X=\mu$ for any $\mu\in M$. Nevertheless, we follow the two-stage procedure from Sections~2 and ~3 to generate Markov chains from $\pi_1$ and $\pi_2$ and build MCMC estimators from Theorems~1-3. The primary goal is to compare the performance of BM and RS estimators.  

We draw iid samples from $\pi_1$ and Markov chain samples from $\pi_2$ using the so called independent Metropolis Hastings algorithm with proposal density $t_{5,1}$. For RS, we follow the idea of \citet[Section 4.1]{mykl:tier:yu:1995} on constructing minorization conditions to identify regeneration times. Based on a carefully tuned minoration condition, the Markov chain for $\pi_2$ regenerates about every 3 iterations on average. In contrast, for users of the BM method proposed in this paper, no such theoretical development is needed. For $i=1, 2$ we draw $N_i$ observations from $\pi_i$ in stage 1 and $n_i$ observations from $\pi_i$ in stage 2. We set $N_1 = N_2$ and $n_1 = n_2 = N_1/10 = N_2/10$.  Recall the reason for smaller stage~2 sample sizes is due to computing cost.  For completeness, note generating Markov chain samples using RS results in a random chain length so these chains were run in such a way that $N_1 \sim N_2$ and $n_1 \sim n_2$.

% \begin{figure}[h!]
%   \begin{center}
%     \includegraphics[width=.45\linewidth]{logtrace-av-a1=05-qmiu1-nrep1000sub100.pdf}
%         \includegraphics[width=.45\linewidth]{logtrace-av-a1=082-qmiu1-nrep1000sub100.pdf}
%   \end{center}
%        \caption{Plots of BM and RS estimates of the asymptotic variance of $\hat{d}$ in stage~1 for 100 randomly chosen replications. The left panel is based on the naive weight, $\baone=(0.5, 0.5)$ and the right panel is based on a close-to-optimal weight, $\baone=(0.82, 0.18)$. Horizontal lines represent the empirical asymptotic variance of $\hat{d}$ obtained over all replications.}
% \label{fig:a1}
% \end{figure}

For estimators based on stage~1 samples, Theorem~1 allows any choice of weight, $\baone$.  For estimators based on stage~2 samples, Theorem~2 and 3 allow any choice of weight, $\batwo$, in constructing consistent BM estimators of the asymptotic variances.  RS based estimators in stage~2 are calculated using Theorems stated in \citet{DossTan:2014} and \citet{tan:doss:hobe:2015} with a general weight choice noted in Remark~4.  This is an important generalization in that now any non-negative numerical weight vector can be used.  We discuss the choice of weights and their impact on the estimators later in this section. 

\begin{figure}[h!]
  \begin{center}
      \includegraphics[width=.32\linewidth]{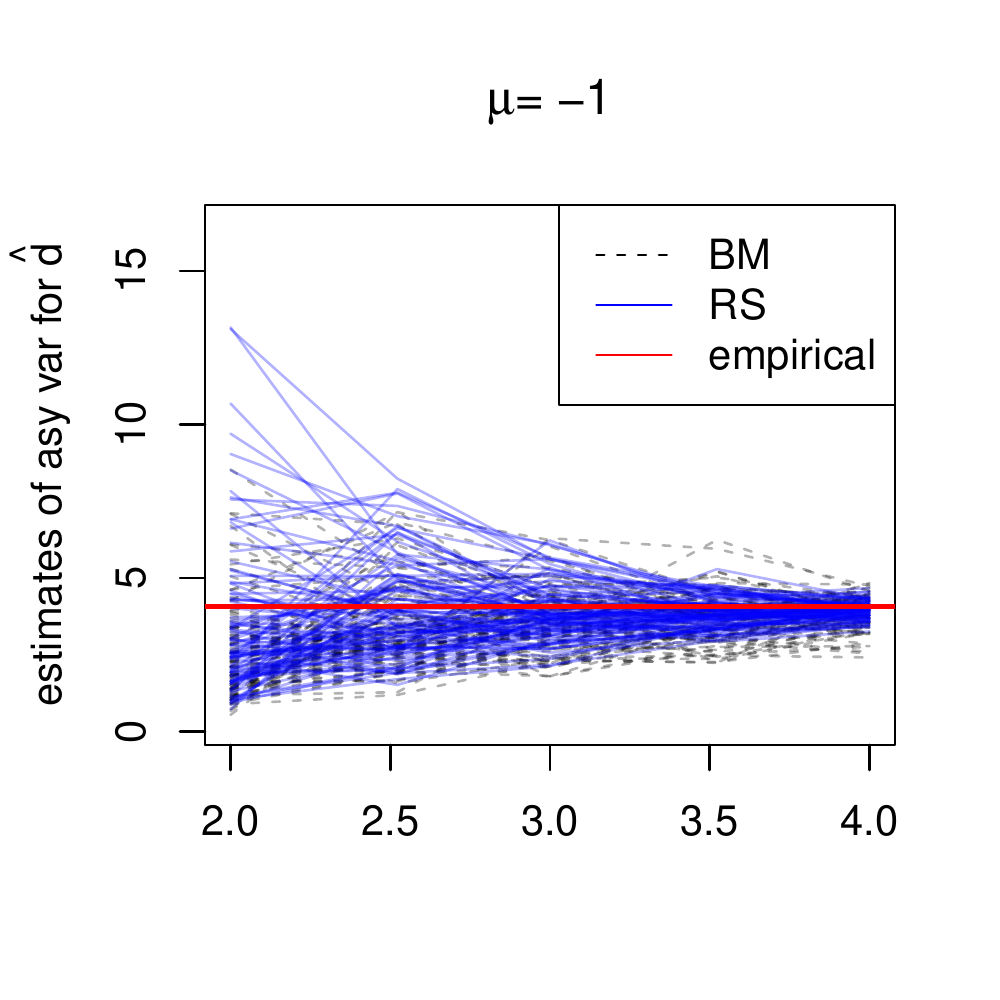}
            \includegraphics[width=.32\linewidth]{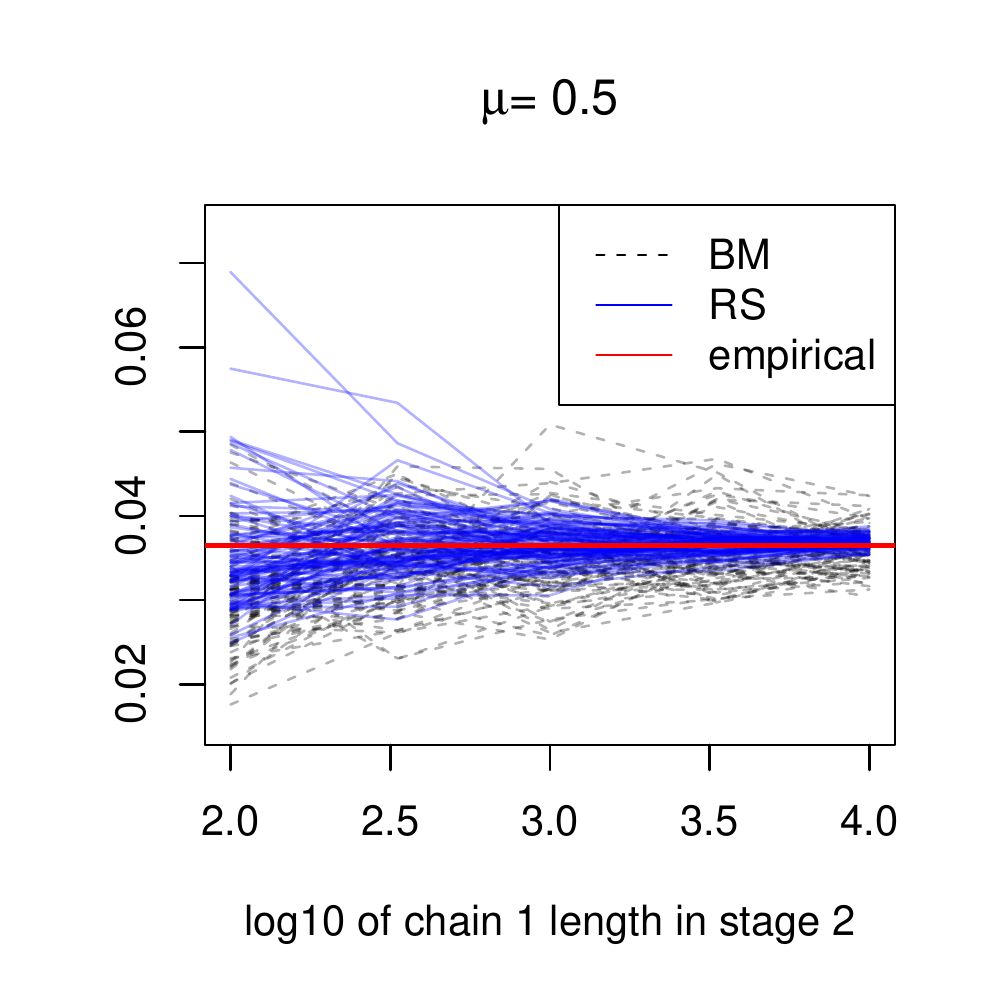}
               \includegraphics[width=.32\linewidth]{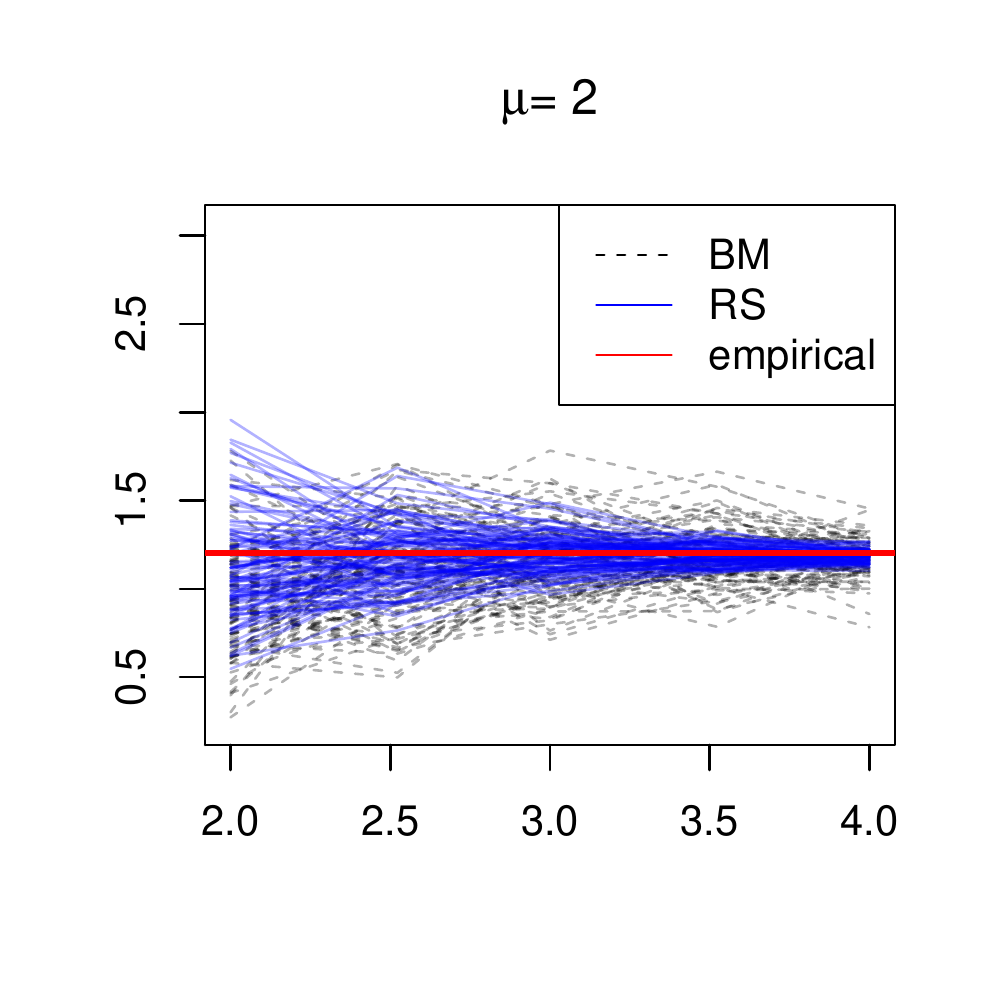}
 
                   \includegraphics[width=.32\linewidth]{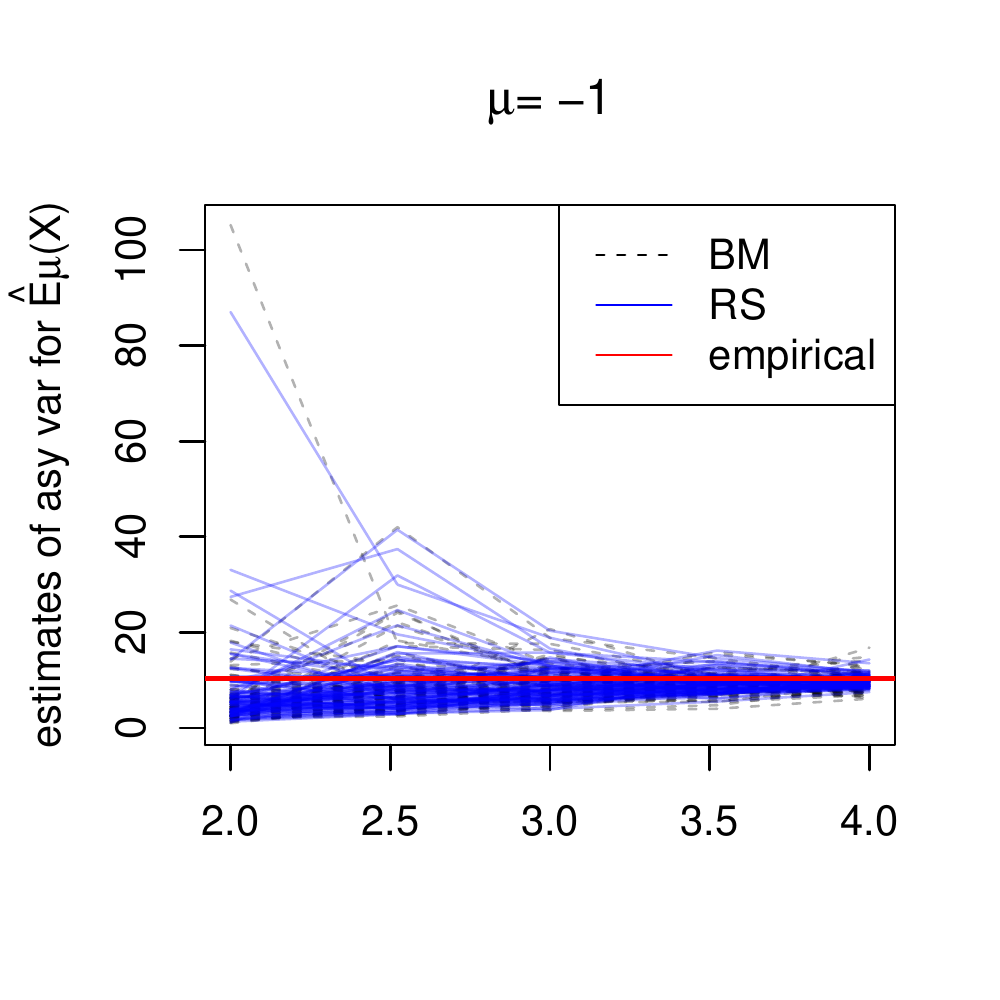}
        \includegraphics[width=.32\linewidth]{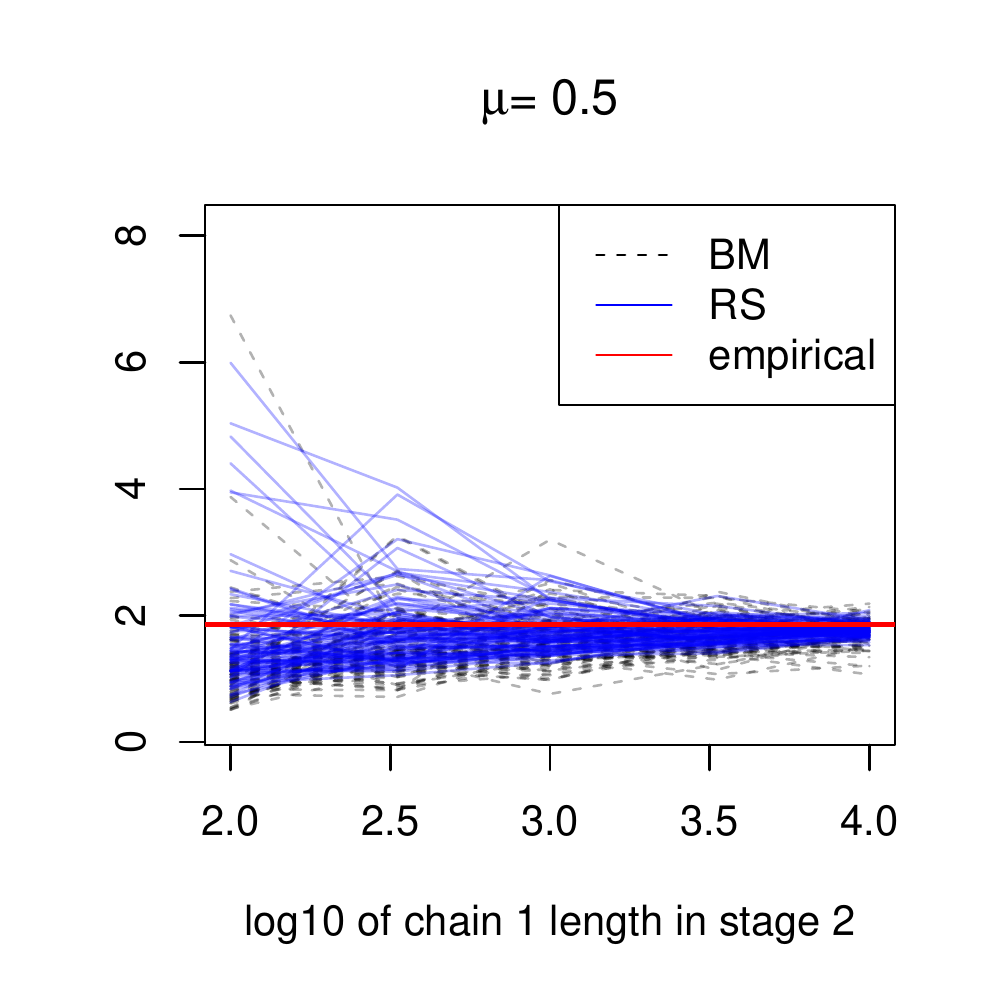}
        \includegraphics[width=.32\linewidth]{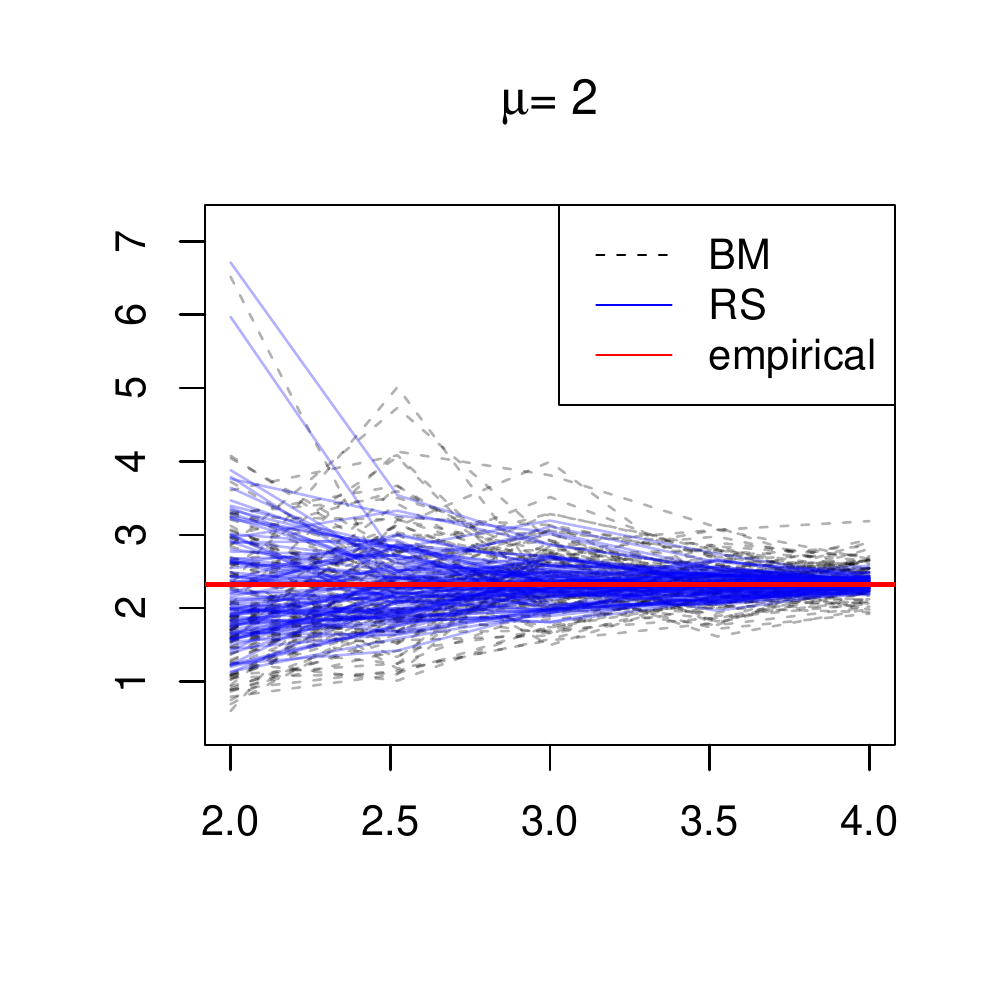}
  
  \end{center}
\caption{Estimates of the asymptotic variance of $\hat{d}_\mu$ (upper panels) and $\hat{\E}_\mu(X)$ (lower panels) in stage~2, with naive weight $\batwo=(.5, .5)$.}
\label{fig:dhat-Ehat-1}
\end{figure}

The following details the simulation study presented in the main text.  We consider increasing sample sizes from $N_1=10^3$ to $10^5$ in order to examine trace plots for BM and RS estimates.  The two stage procedure is repeated 1000 times independently. The unknown true value of the asymptotic variance of $\hat{d}$ is estimated by its empirical asymptotic variance over the $1000$ replications at $N_1=10^5$. %Figure~\ref{fig:a1} displays traces of the BM and the RS estimates of the asymptotic variance of $\hat{d}$ in stage~1, in dashed and solid lines, respectively. For clarity, the traces are only plotted for a randomly chosen subset of 100 replications. The left panel is based on the naive weight, $\baone=(0.5, 0.5)$, that is proportional to the sample sizes; and the right panel is based on $\baone=(0.82, 0.18)$ that weighs the iid sample more than the Markov chain sample.  For either choice of weight, we see that both the BM and the RS estimates approach the empirical asymptotic variance (shown as the horizontal line across the picture) as the sample size increases, suggesting their consistency. As expected due to the frequency of regenerations, BM estimates are more variable than RS estimates. 
We consider the naive weight, $\baone=(0.5, 0.5)$, that is proportional to the sample sizes, and an alternative $\baone=(0.82, 0.18)$ that weighs the iid sample more than the Markov chain sample.  As illustrated in Figure 1 of the main text, both the BM and the RS estimates approach the empirical asymptotic variance as the sample size increases suggesting consistency. Similarly for stage~2, Figure~\ref{fig:dhat-Ehat-1} shows convergence of the BM and the RS estimates to the corresponding empirical asymptotic variances of $\hat{d}_\mu$ and $\hat{\E}_\mu(X)$.  Plots for other $\mu \in M$ show similar results, but are not included here.

Overall, the simulation study suggests BM and RS methods provide consistent estimators for the true asymptotic variance. RS estimators enjoy smaller mean squared error in most cases. Nevertheless, when the number of regenerations is not great, BM estimators could be the more stable estimator. For example, in the top left panel of Figure~\ref{fig:dhat-Ehat-1}, at stage~2 sample size $n_2 \approx n_1 = 100$, or about $35$ regenerations for chain~2, the RS method substantially over-estimated the target in about $5\%$ of the replications. Further, in the cases where regeneration is unavailable or the number of regenerations is extremely small, then BM would be the more viable estimator.

\subsection{Choice of stage 1 weights}
For stage~1, we recommend obtaining a close-to-optimal weight $\baoneopt$ using a pilot study described in \citet{DossTan:2014}. In short, one can generate samples of small size from $\pi_1$ and $\pi_2$, estimate $\hat{d}$ and its asymptotic variance based on Theorem~1 for a grid of weights, and then identify the weight that minimizes the estimated variance.  With a small pilot study based on samples of size $1000$ from both distributions, we obtained $\baoneopt=(0.82, 0.18)$. As depicted by the horizontal lines accross the pictures in Figure~1 of the main text, the asymptotic variance of the estimator $\hat{d}$ based on $\baoneopt$ is approximately $0.07$, which is more than $30\%$ smaller than of the estimator based on the naive choice $\baone=(.5, .5)$. Note that the naive weight is proportional to the sample sizes from $\pi_1$ and $\pi_2$, which is asymptotically optimal if both samples were independent.  However, since sample 2 is from a Markov chain sample, using a weight that appropriately favors the independent sample has lead to smaller error in the estimator. The gain in efficiency using a close-to-optimal weight will be more pronounced if the difference in the mixing rates of the two samples is larger. 

\subsection{Choice of stage 2 weights}
In stage~2, for each $\mu \in M$, the asymptotic variance of $\hat{d}_\mu$ and $\hat{E}_{\mu}(X)$ are minimized at different weights. Instead of searching for each of the $2|M|$ optimal weights in a pilot study, it is more practical to set sub-optimal weights using less costly strategies.  Below, we perform a simulation study to examine three simple weighting strategies:
\begin{enumerate}
\item naive: $\batwo \propto (n_1, n_2)$,
\item inverse distance (inv-dist): $\batwo(\mu) \propto \left(\frac{n_1}{|\mu-\mu_1|}, \frac{n_2}{|\mu-\mu_2|}\right)$, 
\item effective sample size (ess) by inverse distance (inv-dist): $\batwo(\mu) \propto  \left(\frac{\text{ess}_1}{|\mu-\mu_1|}, \frac{\text{ess}_2}{|\mu-\mu_2|}\right)$.
\end{enumerate}

\begin{figure}[h!]
  \begin{center}
    \includegraphics[width=.45\linewidth]{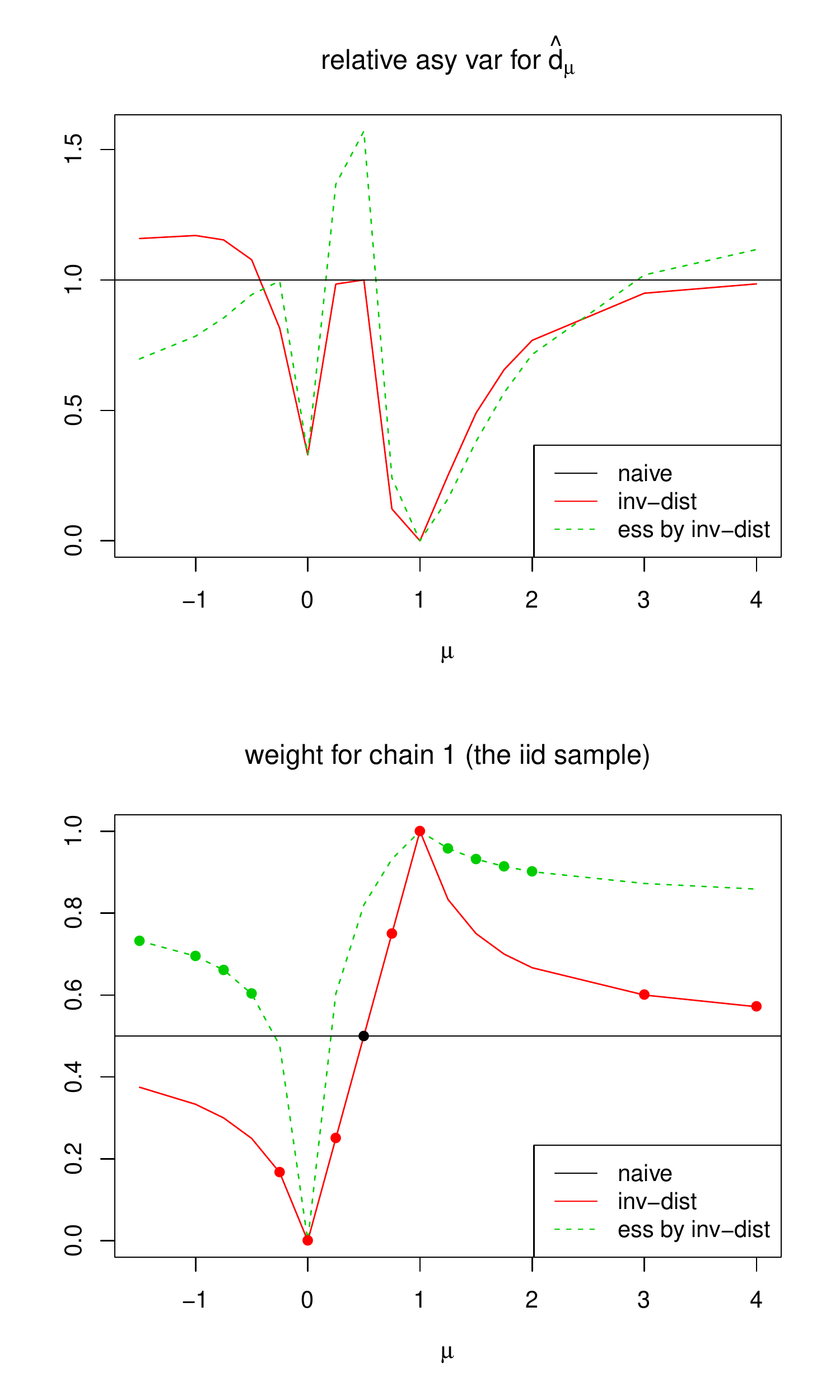}
        \includegraphics[width=.45\linewidth]{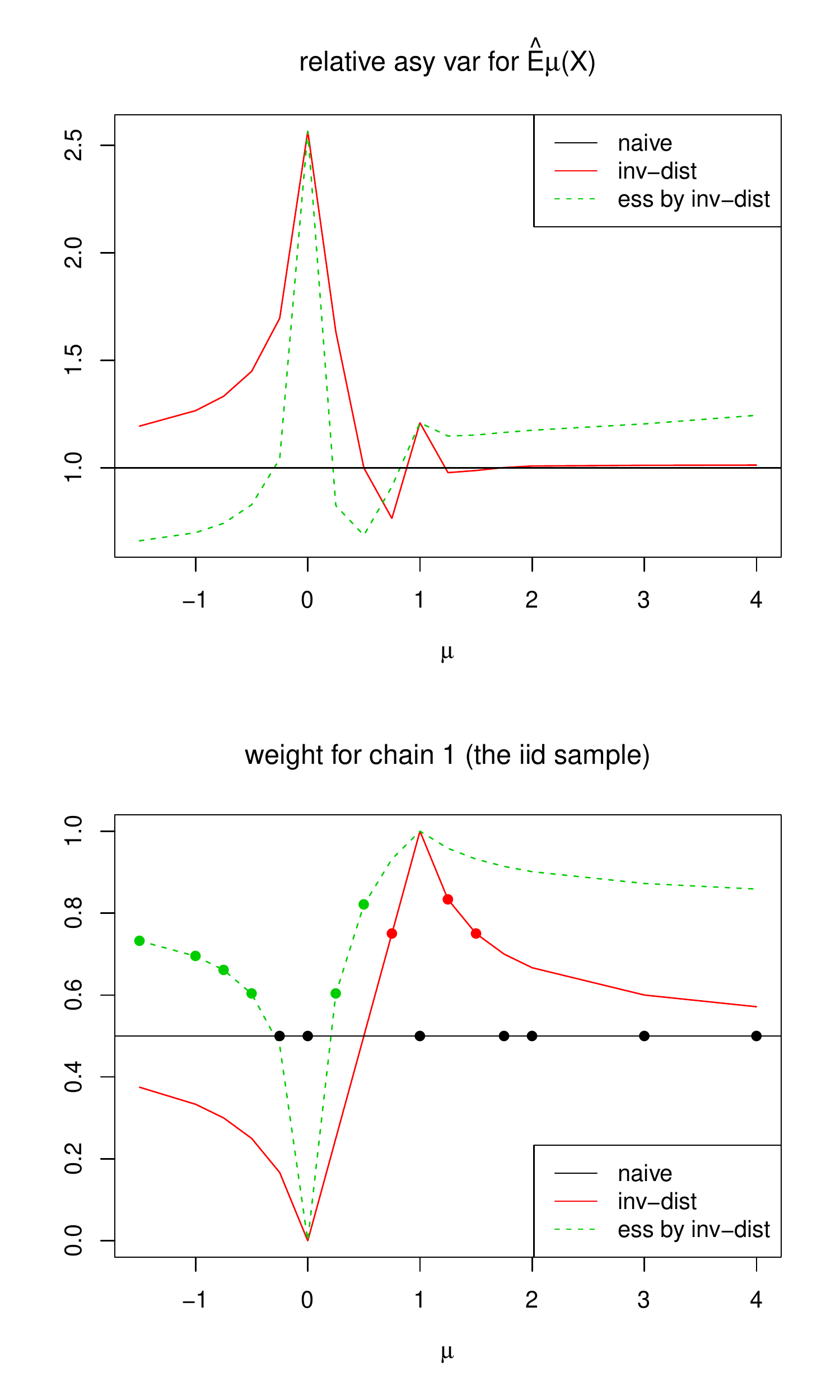}
        \caption{Comparisons of three weight strategies in terms of the asymptotic variance of the
          corresponding estimators $\hat{d}_\mu$ and $\hat{\E}_\mu(X)$. The solid dots show which strategy achieves the smallest
          asymptotic variance among the three at any given $\mu$ (ties awarded to the more basic strategy).}
    \label{fig:a2}
  \end{center}
\end{figure}

\begin{figure}[h!]
  \begin{center}
      \includegraphics[width=.32\linewidth]{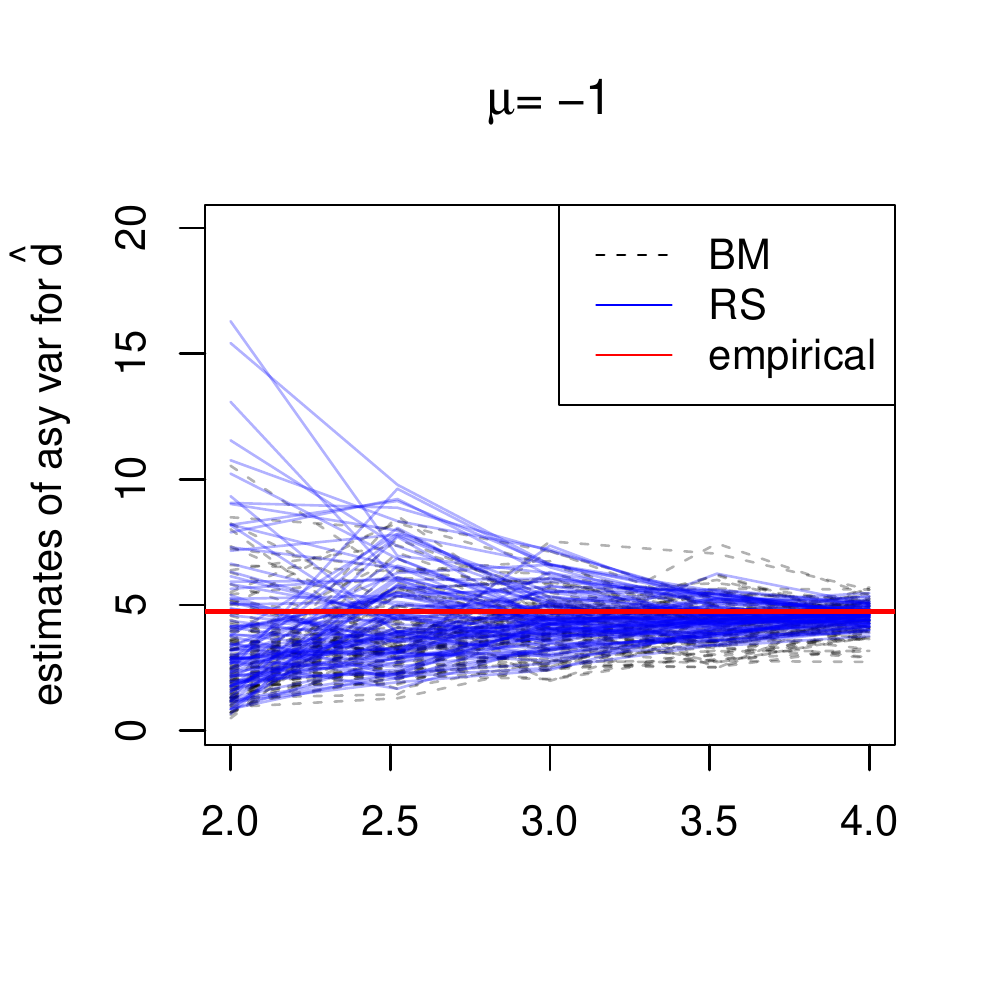}
            \includegraphics[width=.32\linewidth]{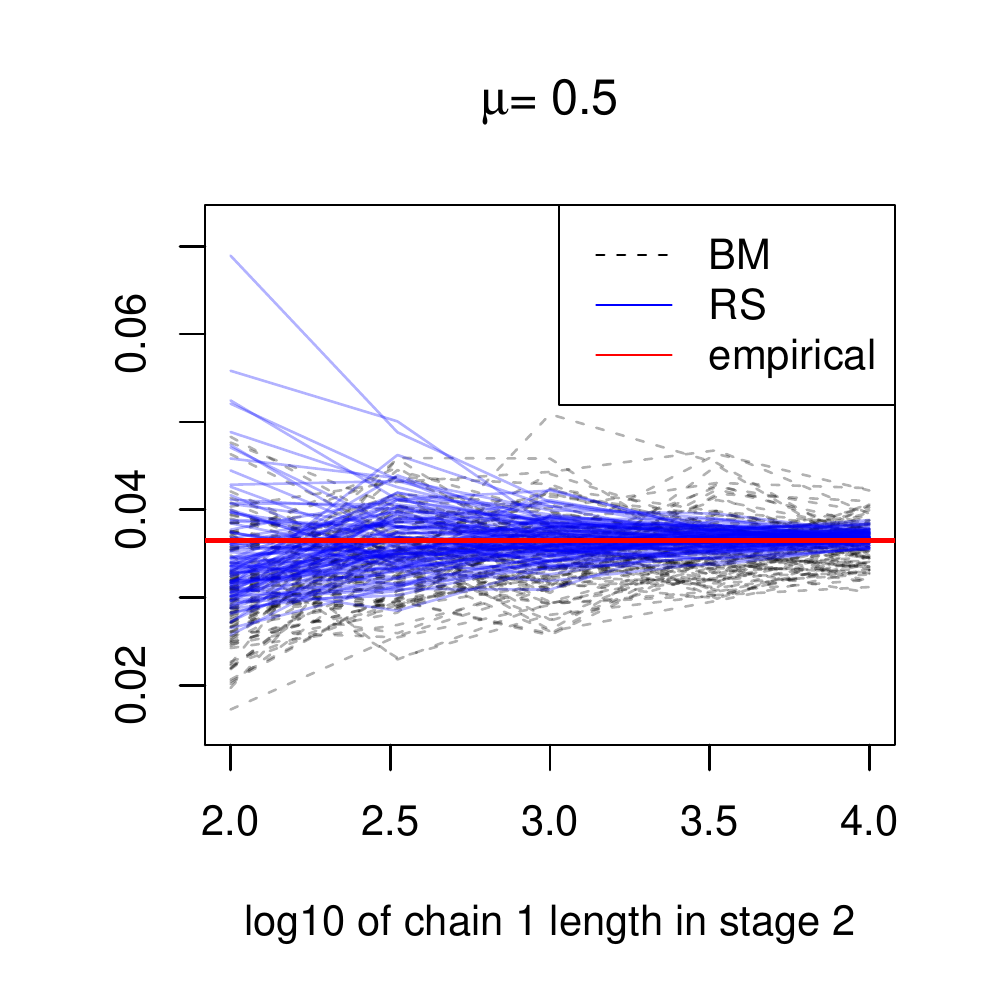}
               \includegraphics[width=.32\linewidth]{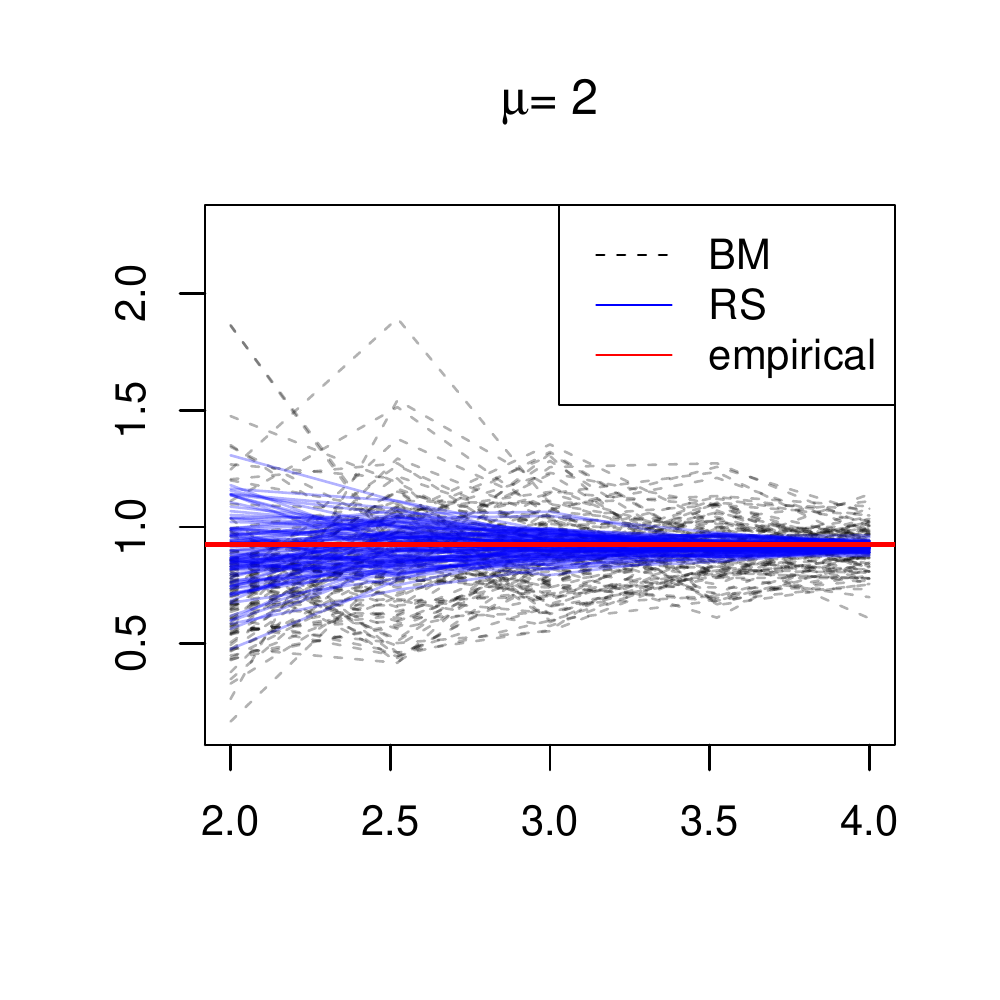}
 
                   \includegraphics[width=.32\linewidth]{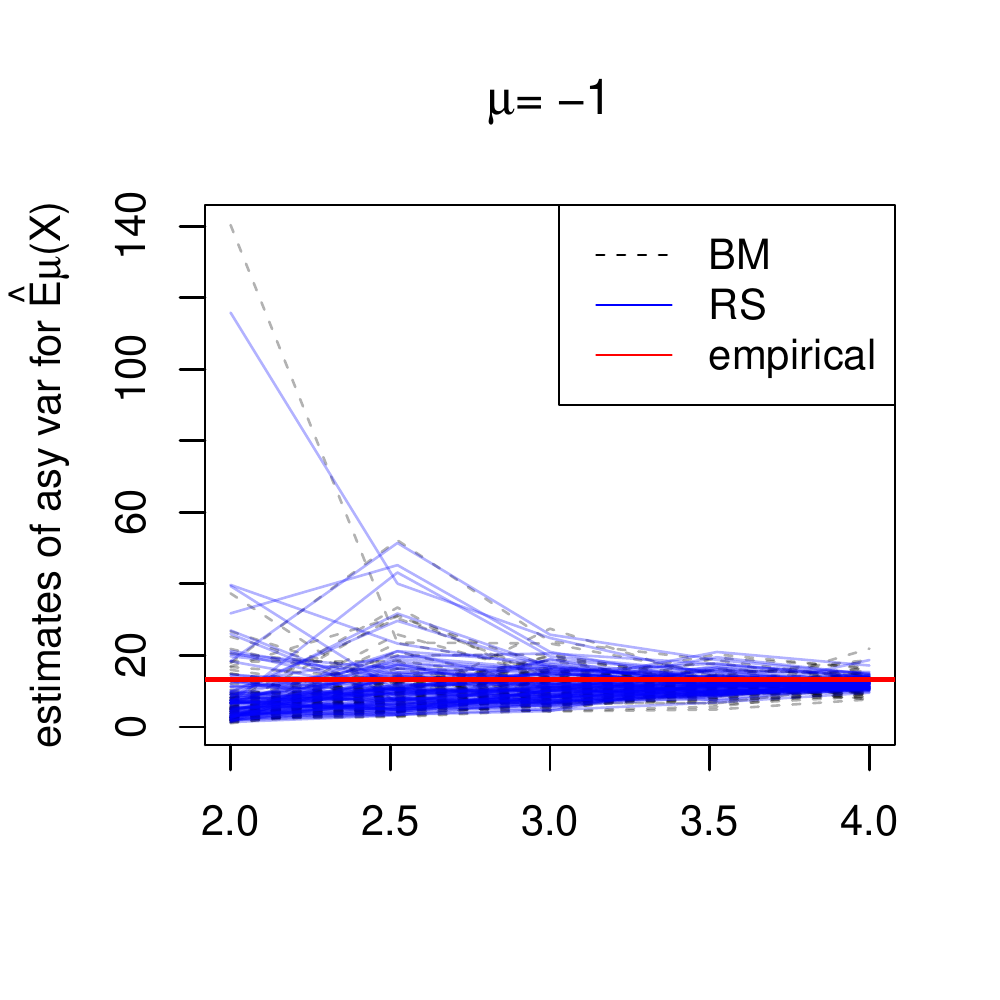}
        \includegraphics[width=.32\linewidth]{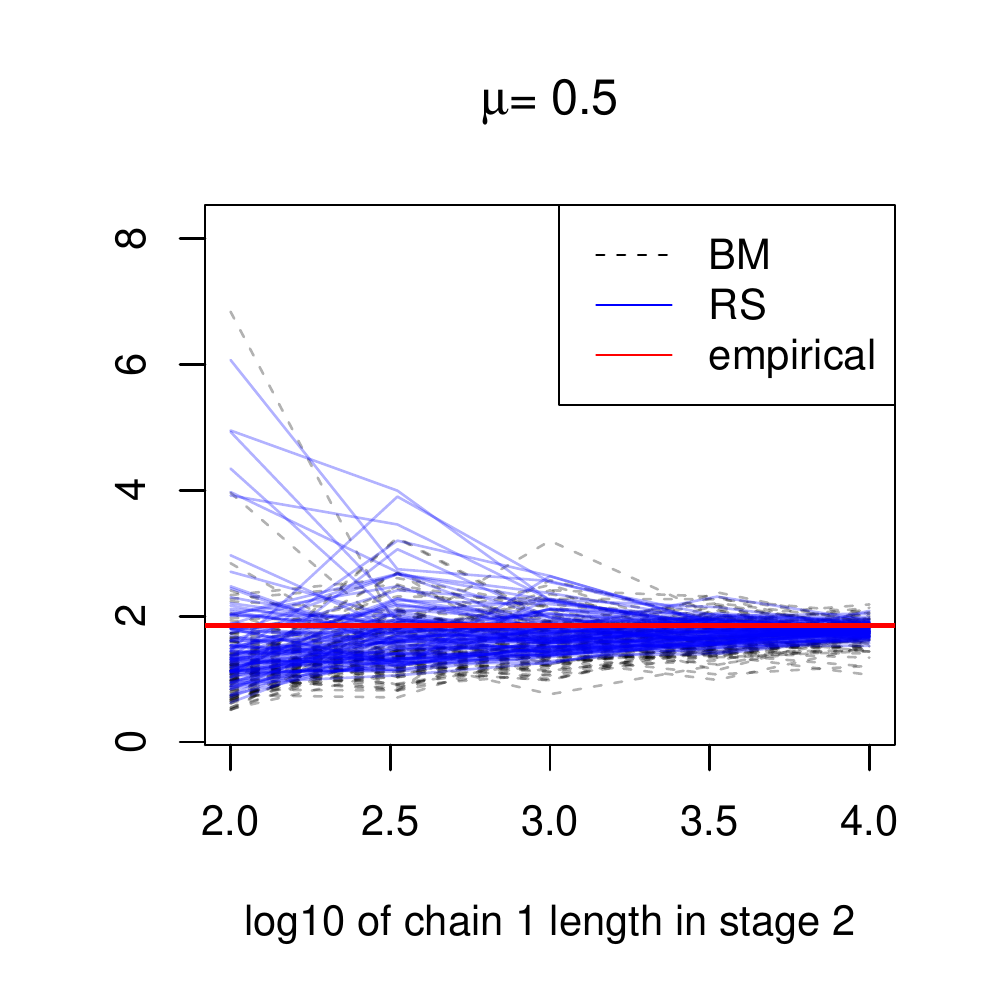}
        \includegraphics[width=.32\linewidth]{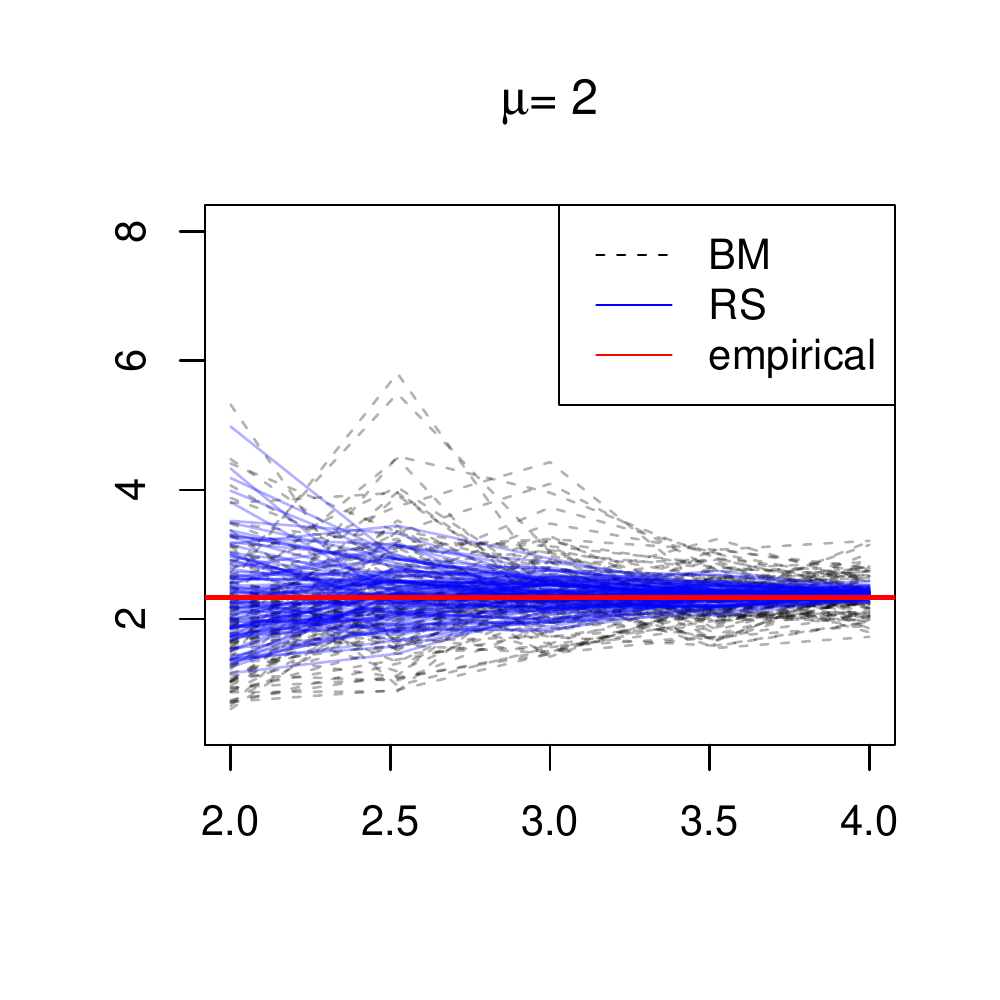}
    \end{center}
\caption{Estimates of the asymptotic variance of $\hat{d}_\mu$ (upper panels) and $\hat{\E}_\mu(X)$ (lower panels) in stage~2, with weight $\batwo(\mu)$ chosen by strategy~2.}
\label{fig:dhat-Ehat-2}
\end{figure}

\begin{figure}[h!]
  \begin{center}
      \includegraphics[width=.32\linewidth]{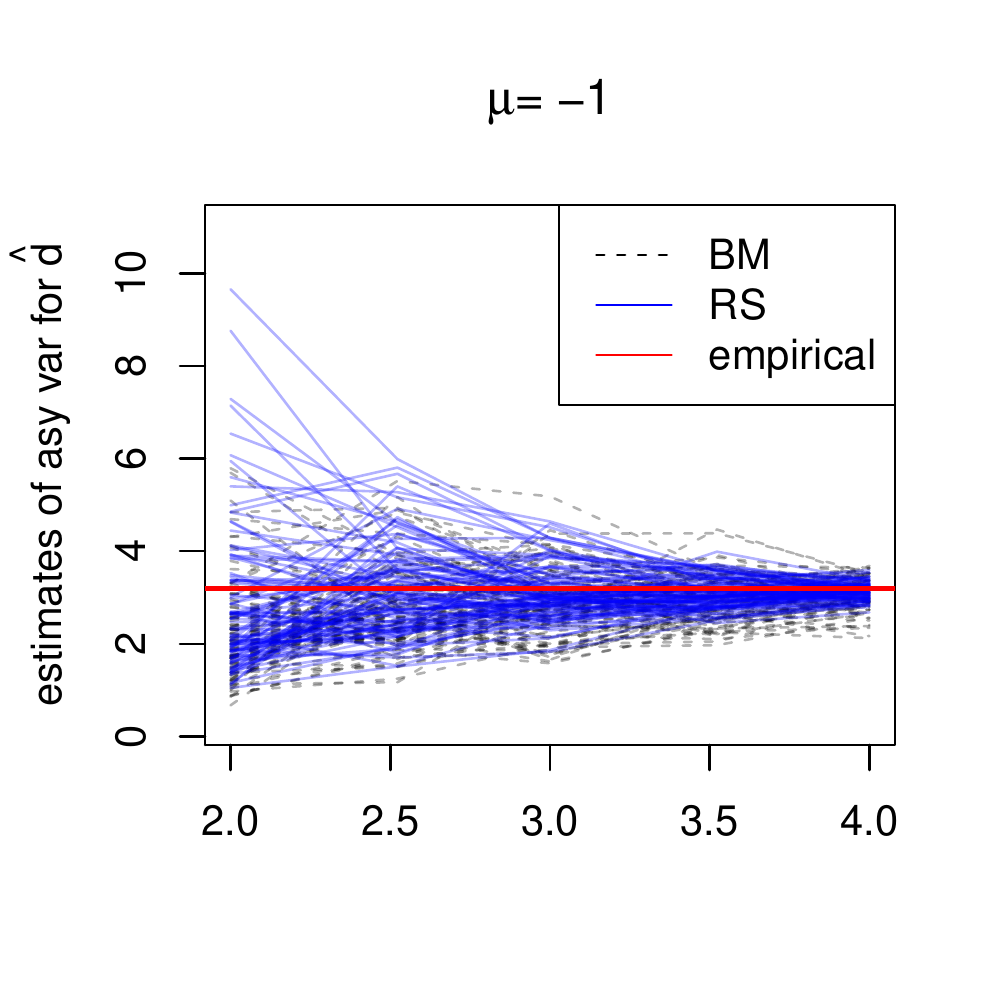}
            \includegraphics[width=.32\linewidth]{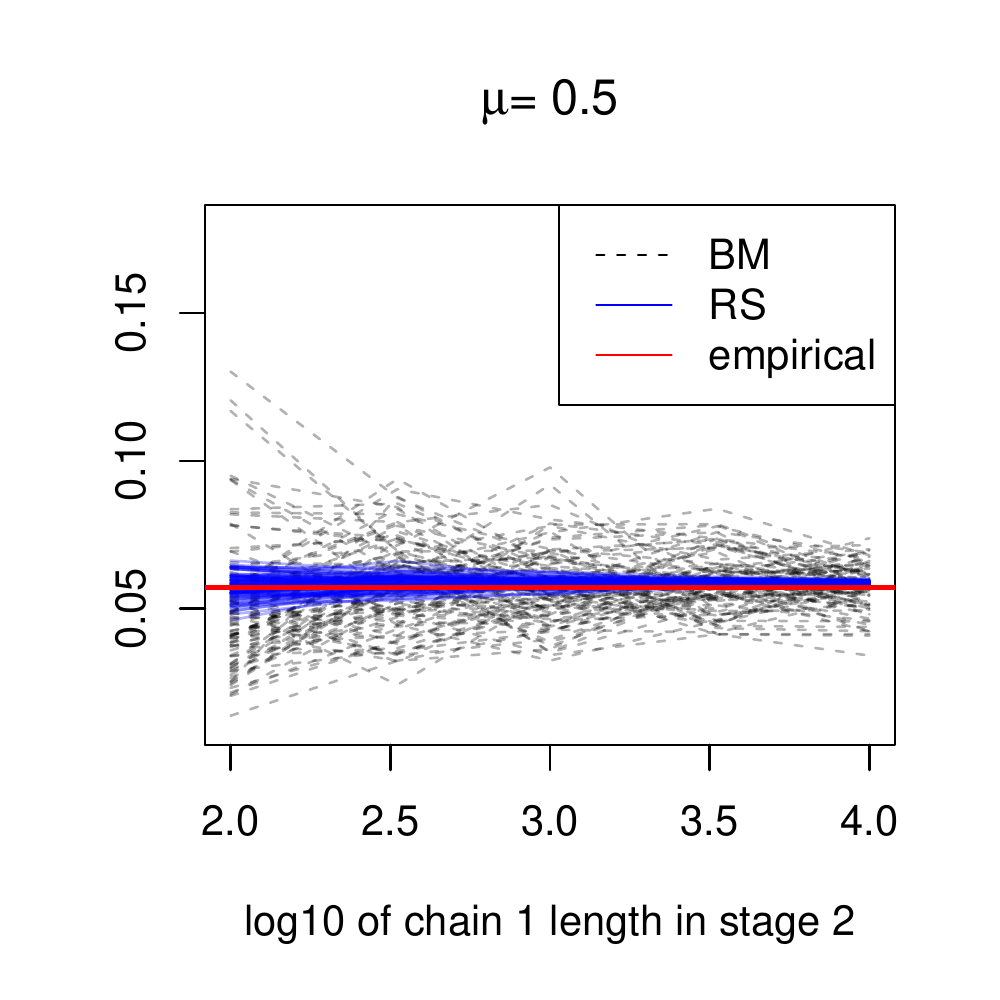}
               \includegraphics[width=.32\linewidth]{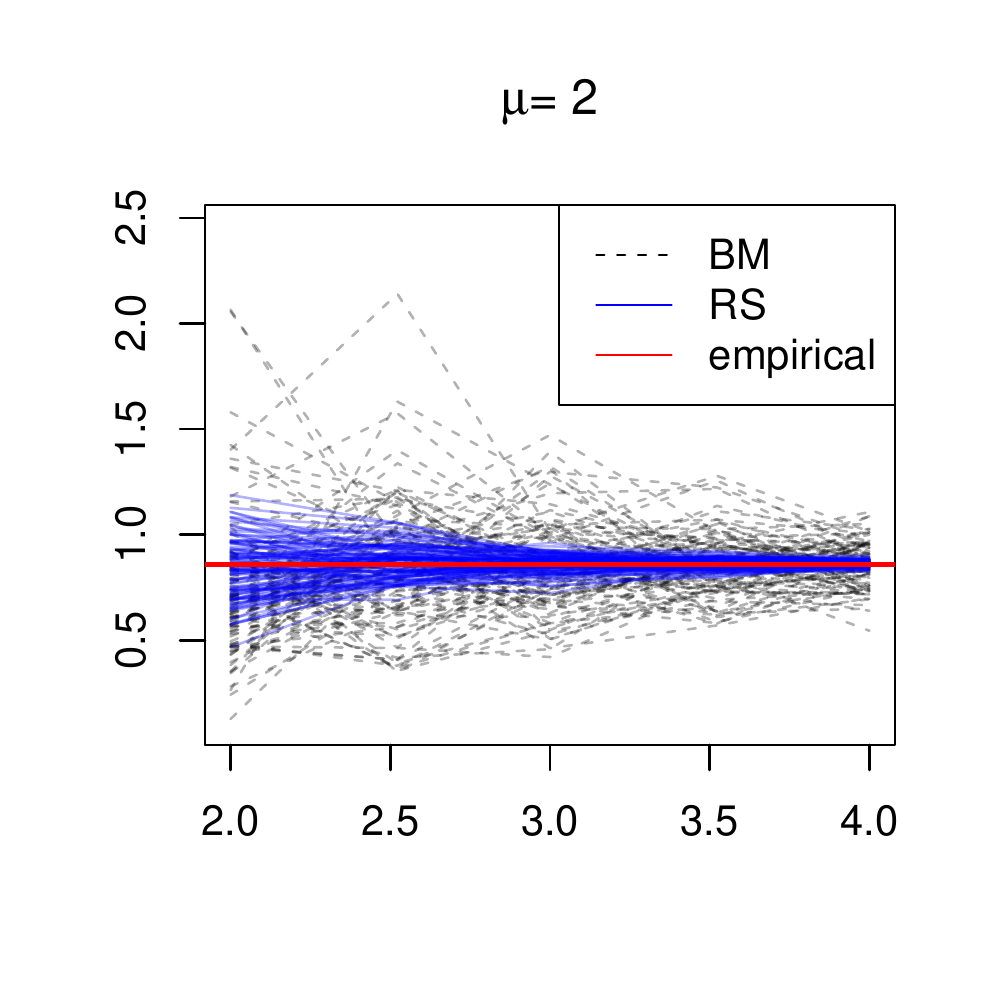}
 
                   \includegraphics[width=.32\linewidth]{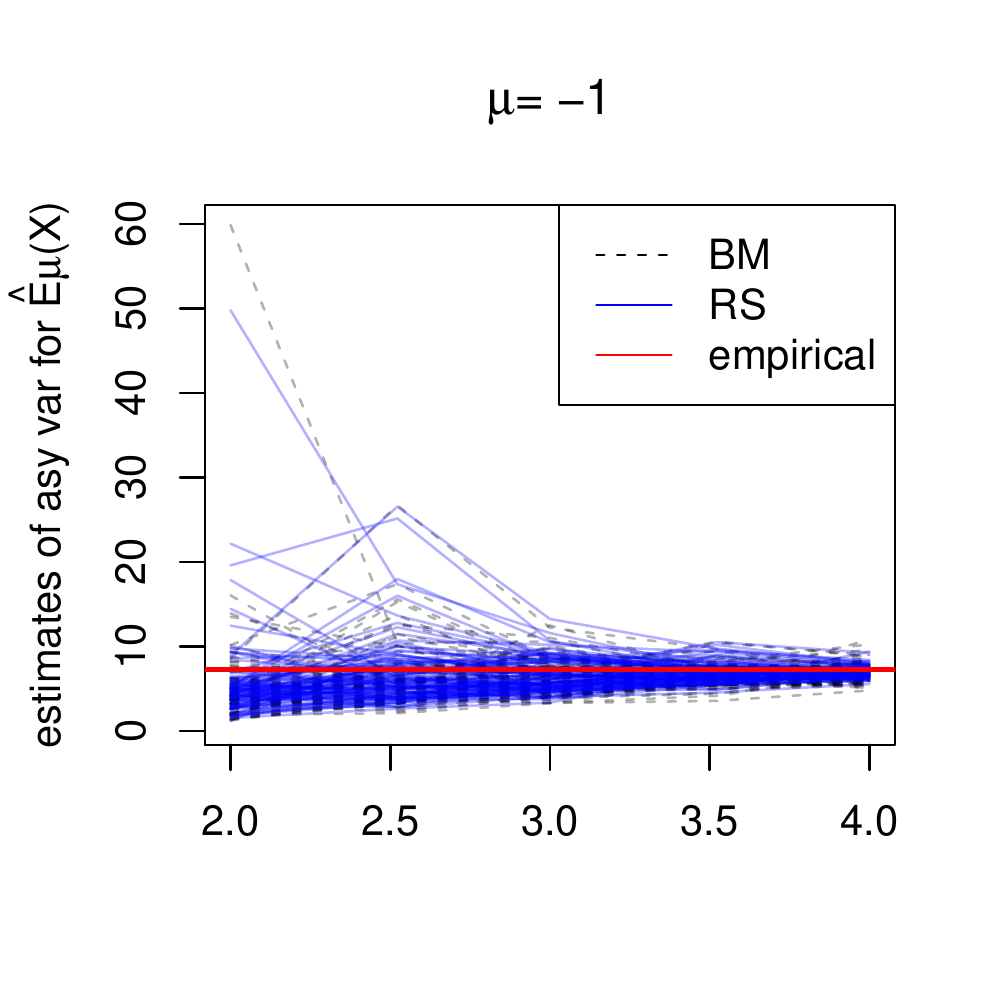}
        \includegraphics[width=.32\linewidth]{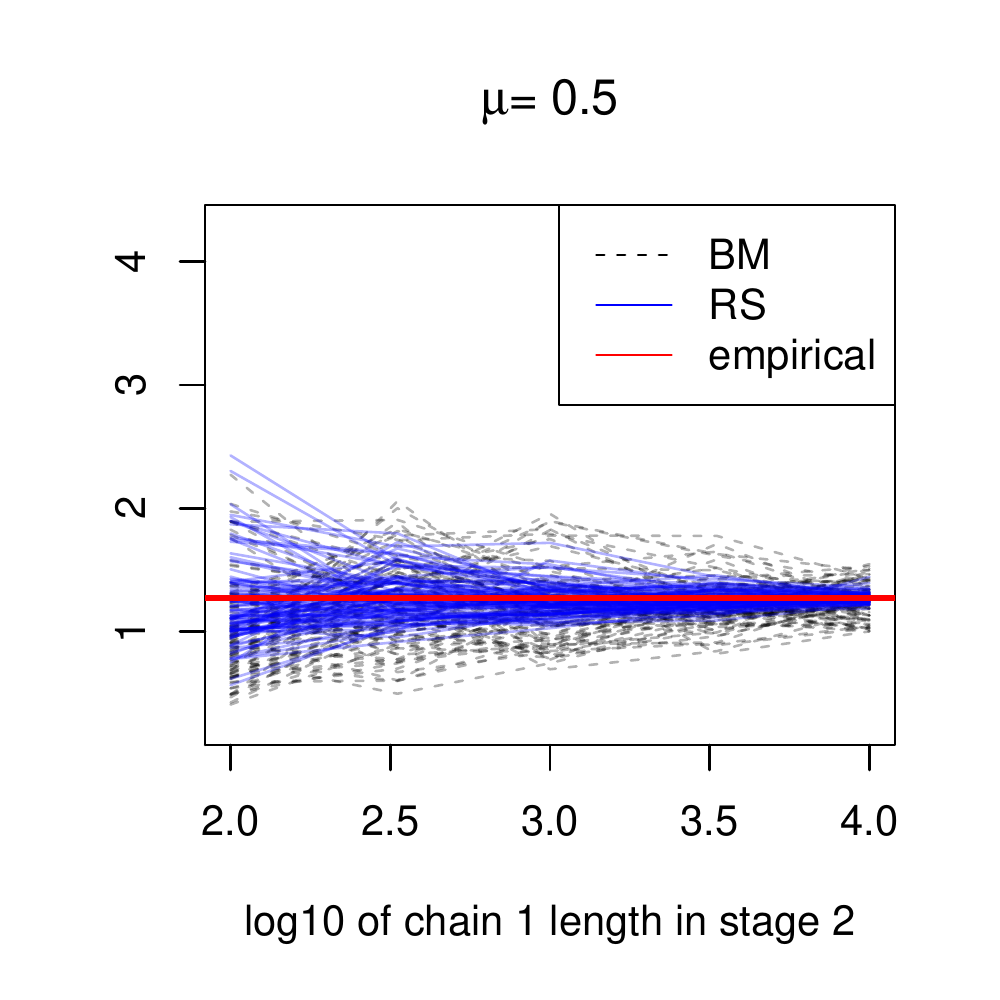}
        \includegraphics[width=.32\linewidth]{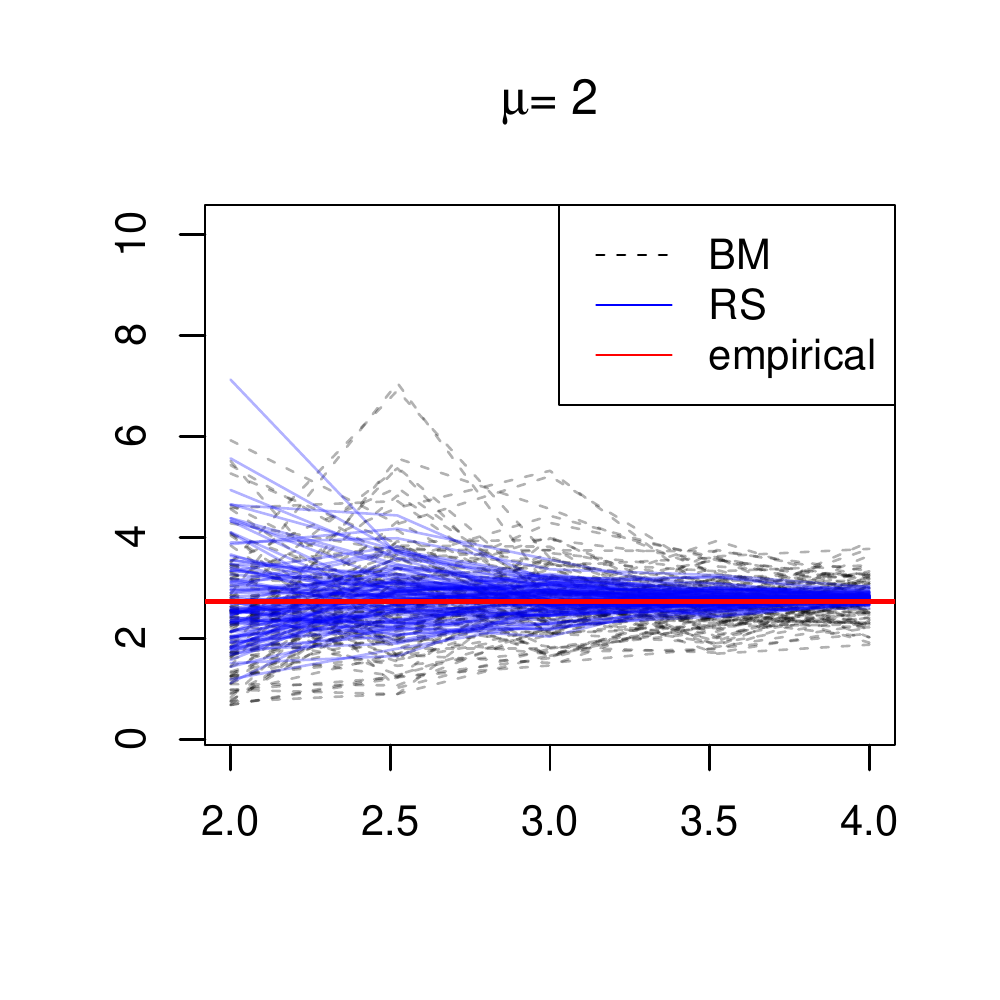}
    \end{center}
\caption{Estimates of the asymptotic variance of $\hat{d}_\mu$ (upper panels) and $\hat{\E}_\mu(X)$ (lower panels) in stage~2, with naive weight $\batwo(\mu)$ chosen by strategy~3.}
\label{fig:dhat-Ehat-3}
\end{figure}

\noindent Using each of the three strategies, we construct generalized IS estimators for $d_\mu$ and $\E_\mu(X)$ for a grid of $\mu$ values between $-1.5$ and $4$. Note that samples are drawn from two reference distributions indexed by $\mu=1$ and $\mu=0$ respectively. Hence our simulation study concerns both interpolation and extrapolation. A summary of their performance is provided in Figure~\ref{fig:a2}, and detailed results for selected simulation setups are shown in Figures~\ref{fig:dhat-Ehat-1},~\ref{fig:dhat-Ehat-2}, and~~\ref{fig:dhat-Ehat-3} for strategies~1, 2, and 3, respectively. Figure~\ref{fig:a2} suggests that none of the three strategies is uniformly better than the others. In particular, we observe the following.

\begin{enumerate}
\item \textbf{ For estimating $d_\mu$}
\begin{enumerate}
\item For $\mu \in (0, 1)$, strategy 2 works the best.
\item For $\mu = 0$, strategies 2 and 3 work better than strategy
  1. Indeed, both of them simply set their stage~2 estimates
  $\hat{d}_0$ to be the stage~1 estimate, $\hat{d}$. This would be a
  better choice than strategy 1 because in a two-step procedure,
  stage~1 chains are often much longer than stage~2 chains, and hence
  $\hat{d}$ is already a very accurate estimate for $d_0=d$.
\item For $\mu \notin [0,1]$, strategies 2 and 3 generally lead to more stable
  estimates of $d_\mu$. However, all strategies lead to very large
  asymptotic variances for $\mu<0$. Hence, one needs to be mindful when doing extrapolation with IS estimators --- always obtain an estimate of the standard error, or reconsider the placement of the reference points. 
\end{enumerate}
\item \textbf{ For estimating $\E_\mu(X)$}
\begin{enumerate}
\item For $\mu \in (0, 1)$, strategy~2 works the best in
  general, while strategy~3 is very unstable.
\item For either $\mu = 0$ or $1$, strategy~2 and 3 are the same, and
  they only utilize the reference chain from $\mu$. This was a wise
  choice for estimating $d_\mu$ as explained before, but not so for other quantities of interest.
\item For $\mu \notin [0,1]$, all strategies lead to fairly large
  asymptotic variances, especially for $\mu<0$.
  \end{enumerate}
\end{enumerate}

Overall in stage~2, strategy 2 has an advantage when the estimands are ratios between normalizing constants. However, when estimating $\E_\mu(X)$, the situation is more complicated. Our impression is that assigning any extreme weight will lead to high variability in the estimator. So it is reasonable to simply use the naive weight, or other strategies that bound the weights away from $0$ and $1$.
%CONTENT OF THIS SECTION........\\
%CONTENT OF THIS SECTION.\\
%CONTENT OF THIS SECTION.\\

\par

\section{Bayesian variable selection models }
\label{sec:realdata}
Here, we consider a class of Bayesian variable selection (BVS) models for linear regression with independent normal priors on the regression coefficients. This model involves a 2-dimensional prior hyperparameter that influences inference, yet no default choice guarantees good performance in practice. Hence, displaying the effect of different hyperparameter values on the posterior distribution would greatly benefit users of the model. When the number of predictors, $q$, is large, the computing is challenging. Our solution is to obtain MCMC samples for a small number of models with different hyperparameter values, based on which generalized IS estimates can be obtained for BFs and other posterior expectations for a large number of models. Again, an important problem in practice is how long the Markov chains need to be run? In this context, the only affordable method that we are aware of is to estimate the SE of these IS estimators using the proposed BM method.

As introduced by \citet{MitchellBeauchamp:1988}, let $Y = (Y_1, \ldots, Y_m)^{\top}$ denote the vector of responses and $X_1, \ldots, X_q$ denote $q$ potential predictors, each a vector of length $m$. The predictors are standardized, so that for $j=1,\ldots,q$, $1_m^T X_j=0$ and $X_j^T X_j= m$, where $1_m$ is the vector of $m$ $1$'s. The BVS model is given by:
\begin{subequations}
  \label{eq:vsblm}
  \begin{alignat}{2}
    \label{eq:vsblm-a}
  & \text{given } \gamma, \sigma^2, \beta_0, \beta_{\gamma}, \quad & Y                   & \sim {\cal N}_{m}(1_m \beta_0 + X_{\gamma} \beta_{\gamma}, \sigma^2 I),               \\[1.5mm]
    \label{eq:vsblm-b}
  & \text{given } \gamma, \sigma^2, \beta_0, \quad                          & \beta_j     & \ind {\cal N} \bigl( 0, \frac{\gamma_j}{\lambda} \sigma^2  \bigr) \;\;\text{for $j=1,\ldots,q$}, \\[1.5mm]
    \label{eq:vsblm-c}
  &  \text{given } \gamma,                                                            & (\sigma^2, \beta_0) & \sim p(\beta_0, \sigma^2) \propto 1 / \sigma^2,                                       \\[1.5mm]
    \label{eq:vsblm-d}
  & \quad                                                        & \gamma              & \sim p(\gamma) = w^{q_{\gamma}} (1 - w)^{q - q_{\gamma}}.
  \end{alignat}
\end{subequations}
The binary vector $\gamma = (\gamma_1, \ldots, \gamma_q)^{\top} \in \{ 0, 1 \}^q$ identifies a subset of predictors, such that  $X_j$ is included in the model if and only if $\gamma_j = 1$, and $|\gamma|=\sum_{j=1}^q \gamma_j$ denotes the number of predictors included. So \eqref{eq:vsblm-a}  says that each $\gamma$ corresponds to a model given by $ Y = 1_m \beta_0 + X_{\gamma} \beta_{\gamma} + \epsilon\,,$ where $X_{\gamma}$ is an $n \times |\gamma|$ sub-matrix of $X$ that consists of predictors included by $\gamma$, $\beta_{\gamma}$ is the vector that contains corresponding coefficients, and $\epsilon \sim {\cal N}_m(0, \sigma^2 I)$. It is sometimes more convenient to use the notation, $Y = X_{0\gamma} \beta_{0\gamma} + \epsilon$, where $X_{0\gamma}$ has one more column of $1$'s than $X_{\gamma}$ and $\beta_{0\gamma}^T=(\beta_0, \beta_\gamma^T)$.  Unknown parameters are $\theta=(\gamma, \sigma, \beta_0, \beta_{\gamma})$ for which we set a hierarchical prior in \eqref{eq:vsblm-b} to \eqref{eq:vsblm-d}.  In \eqref{eq:vsblm-d}, an independent Bernoulli prior is set for $\gamma$, where $w \in (0,1)$ is a hyperparameter that reflects the prior inclusion probability of each predictor. In~\eqref{eq:vsblm-c}, a non-informative prior is set for $(\sigma^2, \beta_0)$. In \eqref{eq:vsblm-b}, an independent normal prior is assigned to $\beta_{\gamma}$, where $\lambda >0$ is a second hyperparameter, that controls the precision of the prior.  Overall, $\theta$ is given an improper prior due to \eqref{eq:vsblm-c} but the posterior of $\theta$ is indeed proper. 

One can actually integrate out $(\beta_\gamma, \beta_0,\sigma^2)$ %from \eqref{eq:vsblm-a} to \eqref{eq:vsblm-c} 
and arrive at the following model with parameter $\gamma$ only:
\begin{equation} \label{eq:vsb}
  \begin{split}
     Y| \gamma \sim \ell_h(\gamma; Y)  &= \int_{{\Real}_+} \int_\Real \int_{\Real^{|\gamma|}} f(Y|\gamma,  \sigma^2, \beta_0, \beta_{\gamma}) f(\beta_{\gamma}|\gamma,  \sigma^2, \beta_0) f(\sigma^2, \beta_0) d \beta_{\gamma} d\beta_0 d {\sigma^2}\\
&  = c_m \,\lambda^{\frac{|\gamma|}{2}} \big| A_{0\gamma} \big|^{-\frac{1}{2}}
 \bigl[  (Y-\overline{Y})^T (Y-\overline{Y}) -\widetilde{\beta}_\gamma^T A_{0\gamma} \widetilde{\beta}_\gamma \bigr]^{-(m-1)/2}\,,\\
   \gamma   \sim p_h(\gamma) \quad &= \quad w^{q_{\gamma}} (1 - w)^{q - q_{\gamma}}.
   \end{split}
\end{equation}
Here, $c_m$ is a constant depending only on the sample size
$m$. Further, $A_{0\gamma}=X_{0\gamma}^T X_{0\gamma}+ \Lambda_{0\gamma} $, where
$\Lambda_{0\gamma}$ is a diagonal matrix, the main diagonal of which
is the $(1+|\gamma|)$-dimensional vector $(0, \lambda,\cdots,
\lambda)$. %It's easy to see that $A_{0\gamma}$ is always invertible. 
Finally, $\widetilde{\beta}_\gamma=A_{0\gamma}^{-1} X_{0\gamma}^T Y$.  
% Note that our notations $\ell_h(\gamma; Y)$ and
% $p_h(\gamma)$ emphasize the fact that both the likelihood and the
% prior depend on the hyperparameter $h=(w, \lambda)$. Based on
% \eqref{eq:vsb}, the posterior distribution of $\gamma$ is given by
% \begin{equation}
%   \label{eq:pih}
%   \begin{split}
%   \pi_h(\gamma|Y) 
%   %&= \int \int \int \pi_h(\gamma,  \sigma, \beta_0, \beta_{\gamma}| y) d\beta_0, d\beta_{\gamma} d_{\sigma^2}\\
%   &= \frac{\ell_h(\gamma; Y) p_h(\gamma)}{\sum_{\gamma'} \ell_{h}(\gamma'; Y) p_{h}(\gamma')} =:  \frac{\ell_h(\gamma; Y) p_h(\gamma)}{m_h}\,.
% %  &= c'_m \,\lambda^{\frac{|\gamma|}{2}} \big| X_{\gamma}^T X_{\gamma}+ \lambda \I_{\gamma} \big|^{-\frac{1}{2}}
% % \bigl[  (y-\overline{y})^T (y-\overline{y}) -\widetilde{\beta}^T A_\gamma \widetilde{\beta} \bigr]^{-(m-1)/2}\ w^{q_{\gamma}} (1 -
% %  w)^{q-q_{\gamma}}\,,
%   \end{split}
% \end{equation}

Using the model at \eqref{eq:vsb} requires specification of the hyperparameter $h = (w, \lambda)$.  Smaller $w$ values assign high prior probabilities to models with fewer predictors, and priors with smaller $\lambda$ values allow selected predictors to have large coefficients. It is common to set $w=0.5$ (a uniform prior on the model space) and $\lambda=1$ (a unit information prior for uncorrelated predictors, see e.g.\ \citet{kass:raft:1995}). One can also choose $h$ adaptively, say according to the marginal likelihood $m_h=\sum_\gamma \ell_h(\gamma; Y) p_h(\gamma)$. A small value of $m_h$ indicates that the prior $p_h$ is not compatible with the observed data, while $h_{\text{EB}}=\arg\max m_h$ is defined to be the empirical Bayesian choice of $h$. The empirical Bayes idea has been successfully applied to various models with variable selection components (see e.g.\ \citet{geor:fost:2000, yuan:lin:2005}).  However, we have not seen this idea being carried out for the model in \eqref{eq:vsblm}, except where $n=p$ and the design matrices are orthogonal (\citet{john:silv:2005, clyd:geor:2000}). Due to the improper prior in \eqref{eq:vsblm-d}, $m_h$ is not uniquely defined. Nevertheless, the Bayes factor among any two models, say $m_h/m_{h'}$, is well-defined because the same improper prior is assigned to the shared parameters of the two models (see e.g.\ \citet[sec.5]{kass:raft:1995} and \citet[sec.2]{LiangEtal:2008}). 
% Note also that the EB method is not without potential problems --- \citet{scot:berg:2010} show for a BVS model similar to ours that uses the g-prior does not necessarily lead to a solution that approaches any fully Bayesian procedures asymptotically. Hence, the motivation is strong to perform  a sensitivity analysis over the choice of $h$. 

Here, we concentrate on two goals. The first is to evaluate $\{ m_h/m_{h_1}, h \in {\cal H} \}$, the marginal likelihood of model $h$ relative to a reference model $h_1$, which allows us to identify the empirical Bayesian choice of $h$. The second is to evaluate the posterior mean of the vector of coefficients $\beta$ for each $h \in {\cal H}$, which we denote by $\eb_h$. Predictions can then be made for new observations using $Y^\text{(new)}=(x^\text{(new)})^T \eb_h$. %, and their performance can be checked on a test dataset when available. 
%Our goal is to estimate $m_h / m_{h_1}$ and $\eb_h$ for a large set of $h \in {\cal H}$ using MCMC methods.  

%Note that exact evaluation can be done for $q<20$ via brute force, and is manageable for $q<25$ had gray code ordering been used (\citet{diac:holm:1994}). Larger $q$ values require MCMC methods. It turns out that approximating $m_h$ itself is very difficult (\citet{NewtonRaftery:1994}) whereas approximating $m_h / m_{h_1}$ for any fixed $h_1 \in {\cal H}$ is computationally easier, and still produces the BFs and EB estimates we desire. Therefore, our goal is to estimate $m_h / m_{h_1}$ and $\eb_h$ for a large set of $h \in {\cal H}$ using MCMC methods.  Our strategy is to generate Markov chains for $\gamma$ at several $h$ values that scatters in ${\cal H}$, from which we build generalized IS estimators. 

%Recall the formula of the generalized IS estimators in \eqref{eq:uvhat} and \eqref{eq:uv}. Calculating these reduces to evaluating $f(\gamma)=\widetilde{\beta}_\gamma$ and ratios of the following sort:
% \[
% \frac{\nu_{h'}(\gamma)}{\nu_h(\gamma)}=\frac{ \ell_{h'} (\gamma; Y) p_{h'}(\gamma) } { \ell_{h}(\gamma; Y) p_{h}(\gamma) }\,,
% \]
% where the likelihood and the prior were given in \eqref{eq:vsb}.

% many MCMC algorithms can be used to explore the
% distribution $\pi_h$. Well known examples include the simple Gibbs
% sampler (\citet{geor:mccu:1993}), various random walk Metropolis
% Hastings (MH) algorithms (\cite{nott:gree:2004}), and the ODA algorithm
% and its variations (\citet{ghos:clyd:2011, ghos:tan:2015}).  Here, 

For model~\eqref{eq:vsb} with a fixed $h$, a Metropolis Hastings random-swap algorithm (\citet{clyd:ghos:litt:2011}) can be used to generate Markov chains of $\gamma$ from its posterior distribution. In each iteration, with probability $\rho(\gamma)$, we propose flipping a random pair of $0$ and $1$ in $\gamma$, and with probability $1-\rho(\gamma)$, we propose changing $\gamma_j$ to $1-\gamma_j$ for a random $j$ while leaving other coordinates untouched. We set $\rho(\gamma)=0$ when $\gamma$ corresponds to the null model or the full model, and $\rho(\gamma)=.5$ otherwise. Finally, the proposal is accepted with an appropriate probability. Since this Markov chain lies on a finite state space, it is uniformly ergodic and hence polynomially ergodic as well. Further, moment conditions in Theorems 2 and 3 are satisfied because they reduce to summations of $2^q$ terms, a large but finite number. %After all, the multiple importance sampling estimators and our batch means method can be used to study the BVS models in this section.}%In words, the random-swap algorithm allows the $\gamma$-chain a mixed experience of birth, death, and flip, which greatly facilitates the chain to escape local modes in the posterior distribution (\citet{deni:mall:smit:1998}).
To achieve the goal, we generate Markov chains of $\gamma$ with respect to model~\eqref{eq:vsb} at several $h$ values that scatters in ${\cal H}$, from which we build generalized IS estimators, and estimate their standard errors.

% \begin{remark} We mention that there is a popular alternative prior
%   for $\beta_\gamma$, namely the Zellner's $g$-prior
%   (\citet{zell:1986}). Instead of \eqref{eq:vsblm-b}, this prior sets $
%   \beta_\gamma \sim {\cal N}_{|\gamma|} \bigl( 0, g \left(X_\gamma^T
%     X_\gamma\right)^{-1} \sigma^2 \bigr)$. A lot of research has been
%   done on how to specify the hyperparameter, $(w, g)$, in the BVS
%   model with the $g$-prior, partly because such a choice is closely
%   related to variable selection criteria such as AIC or BIC (see e.g.\ \citet{clyd:geor:2004, LiangEtal:2008}).  Also, a sensitivity analysis for its choice of hyperparameter based on
%   generalized IS estimator was performed in
%   \citet{tan:doss:hobe:2015}. 
%   On the other hand, there are limited
% discussions on the choice of the hyperparameter for the BVS model with
% the independent normal prior. 
%   The $g$-prior also does not apply directly to linear
%   regression problems with $n<q$, unless extra care is given to
%   restrict the $\gamma$ space to include models with $|\gamma|<q$
%   only.
% \end{remark}

\subsection{Cookie dough data}
We demonstrate the aforementioned sensitivity analysis using the biscuit dough dataset (\citet{osborne:1984, brown:2001}). The dataset, available in the {\tt R} package {\tt ppls} (\citet{krae:boul:2012}), contains a training set of $39$ observations and a test set of $31$ observations. These data were obtained from a near-infrared spectroscopy experiment that study the composition of biscuit dough pieces. For each biscuit, the reflectance spectrum is measured at $700$ evenly spaced wavelengths. % between $1100$nm to $2498$nm. 
We use these measurements as covariates to predict the response variable, the percentage of water in each dough.  % Since measurements at $700$ consecutive wavelengths are highly correlated, w
We follow previous studies (\citet{hans:2011}) and thin the spectral to $q=50$ evenly spaced wavelengths. 

\label{sec:ozone}
\begin{figure}[h!]
  \begin{center}
    \includegraphics[width=.8\linewidth]{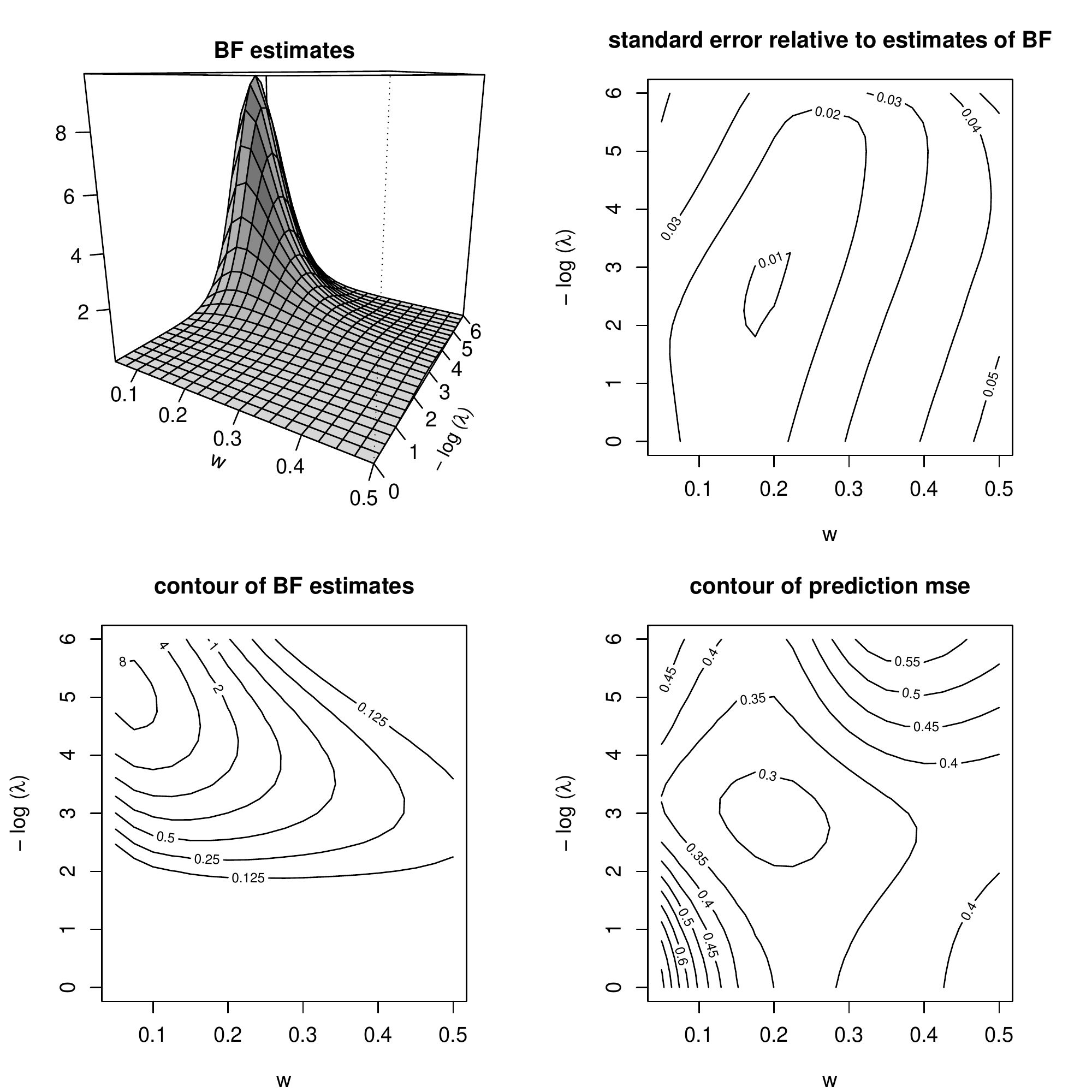} 
  \end{center}
  \caption{
  Left panels provide a surface plot and a contour plot for BF estimates. The upper-right panel displays standard errors with respect to the BF estimates.  The lower-right panel shows pmse over the test set.
}
      \label{fig:bf}
\end{figure}

\begin{figure}[h!]
  \begin{center}
    \includegraphics[width=.4\linewidth]{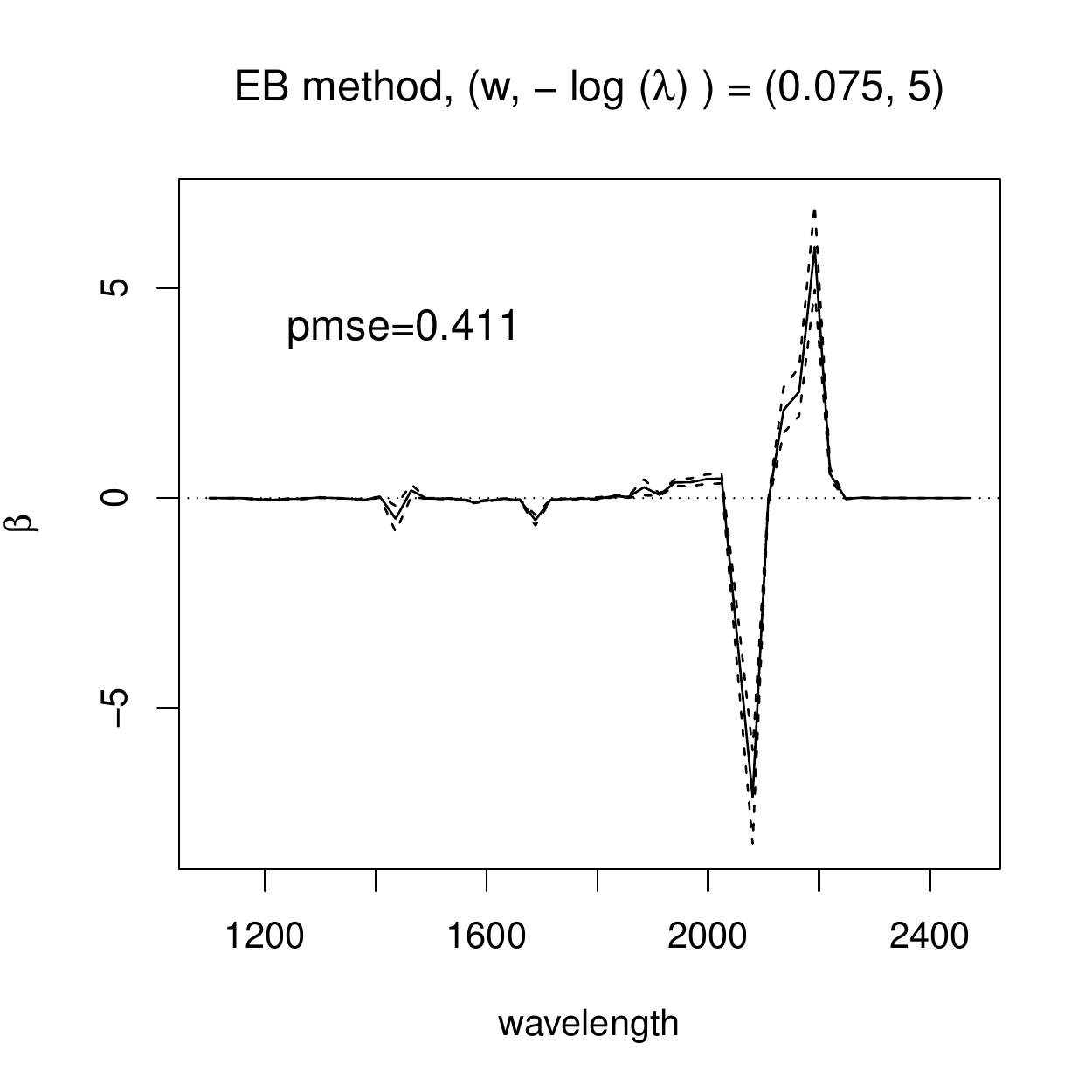} 
        \includegraphics[width=.4\linewidth]{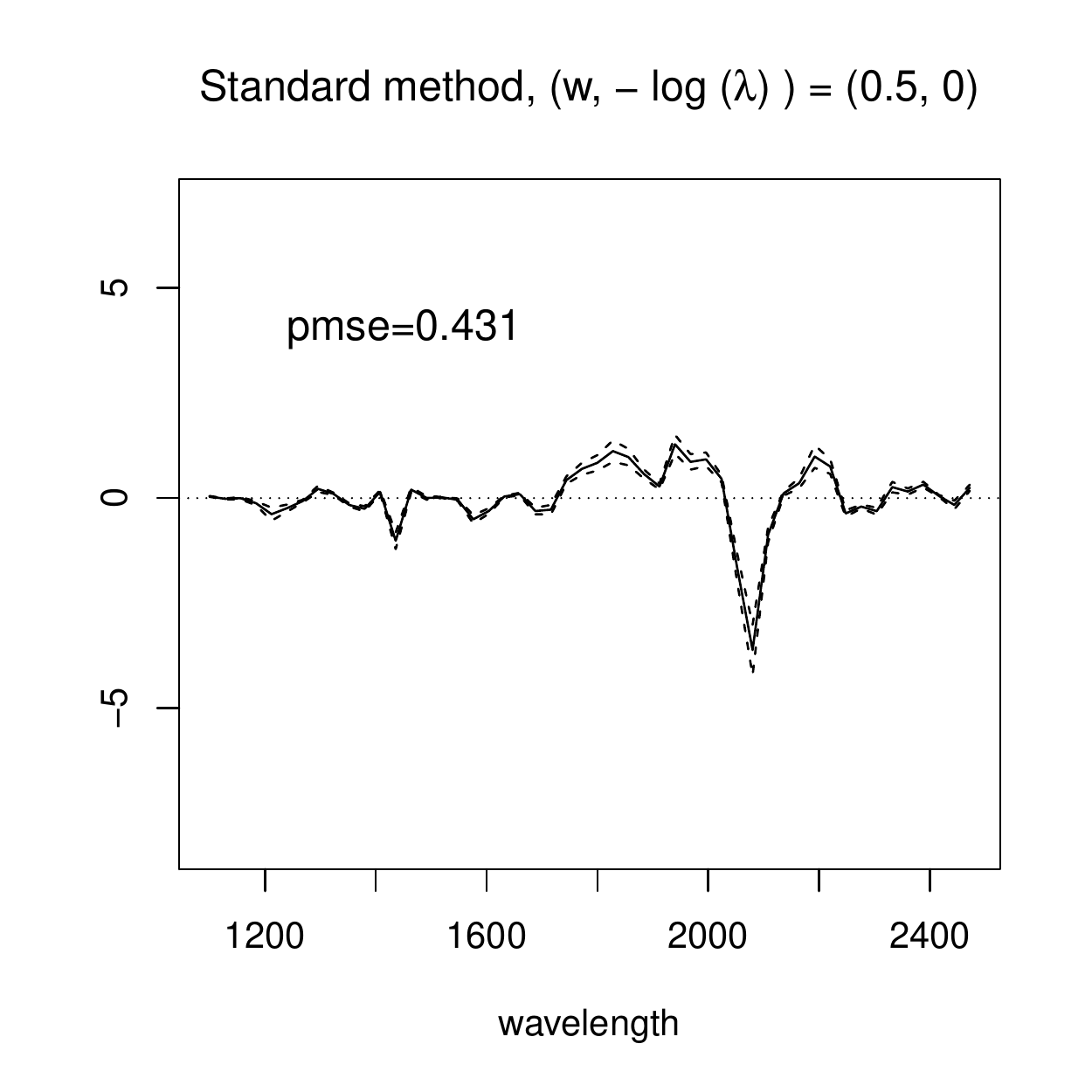} 
  \end{center}
       \caption{Estimated posterior mean of regression coefficients (with $95\%$ point-wise confidence intervals) for the empirical Bayesian and the standard choice of $h$, respectively.}
      \label{fig:beta}
\end{figure}

Figure~\ref{fig:bf} provides a general picture of the sensitivity
analysis. The left plots provide two ways to visualize
estimates of the BFs. To form the plot, we took the $12$ reference
values of $h=(w, \lambda)$ to be such that $\left(w, -\log(\lambda)\right) \in \{ 0.1, 0.2, 0.3, 0.4 \} \times \{ 1, 3, 5\}$.
In stage~$1$ we ran each of the $12$ Markov chains at the above values
of $h$ for $10^5$ iterations to obtain $\hatbd$.  In stage~$2$,
we ran the same $12$ Markov chains for $50,000$ iterations each, to
form the estimates $\hat{u}_n$ over a fine grid that consists of $475$
different $h$ values, with the $w$ component ranging from $0.05$ to
$0.5$ in increments of $0.025$ and the $-\log(\lambda)$ component
ranging from $0$ to $6$ in increments of $0.25$.  % We designed ${\cal
%   H}$ to contain the common choice $(w, \lambda) = (0.5, 1)$, as well
% as a set of $h$ values that have relatively high marginal
% likelihoods. The latter was identified roughly using a small pilot
% study. The reference points were chosen to be a spaced out grid
% inside ${\cal H}$. Finally, the baseline $h_1$ could be any point, and
% we selected $h_1 = (0.2, e^{-3})$ because it corresponds to a
% high marginal likelihood that allows the surface of BFs to
% be displayed easily.

How trustworthy are these BF estimates? Their estimated standard errors are obtained using the BM method, based on Theorem 2. We choose to display the relative SE with respect to the BF estimates, as shown in the upper right panel of Figure~\ref{fig:bf}. The relative SEs are smaller than or equal to $5\%$, and we believe the BF estimates are accurate enough. Finally, the lower-right panel of Figure~\ref{fig:bf}  shows the prediction mean squared error (pmse) over the test set for all $h$.

Based on our estimation, the BF attains the maximum value $9.75$ at
$h_{\text{EB}}=(0.075, e^{-5})$. Recall when comparing any two models
indexed by $h$ and $h'$ respectively, the BF between them is given by
$\text{BF}_{h,{h'}}=\frac{m_h/m_{h_1}}{m_{h'}/m_{h_1}}$. Also,
according to \citet{jeff:1998} and \citet{kass:raft:1995}, the
evidence for $h$ over $h'$ is considered to be strong only if
$\text{BF}_{h,h'}$ is greater than $10$ or $20$. Hence, all $h$ with
$\text{BF}$ over $1/10$ or $1/20$ times the maximum BF can be
considered as reasonably well supported by the data as that of the empirical Bayesian
choice. Comparing the lower two plots of Figure~\ref{fig:bf},
we see
that the set $A_c:=\{h \in {\cal H}: \text{BF}>c\}$ for $c=1$ and
$0.5$ do overlap with an area that corresponds to relatively small
pmse. Outside $A_{0.5}$, a region that consists of larger $w$ and
smaller $-\log(\lambda)$ also enjoys small pmse values, at around
$0.3$ to $0.4$. This region includes the common choice of $h_0=(0.5,
e^0)$. These suggest that $h_{\text{EB}}$ and its vicinity might not be the only area of $h$ that has good prediction performances. 

To better compare the effect of $h_{\text{EB}}=(0.075, e^{-5})$ and
the commonly used $h=(0.5, e^0)$, Figure~\ref{fig:beta} 
displays the
estimated posterior mean of regression coefficients at both choices of
$h$, together with the point-wise $95\%$ confidence intervals for the
posterior means. %These are obtained with the estimators we proposed in Sections~2 and 3. 
Due to the small size of $w_{\text{EB}}$ and
$\lambda_{\text{EB}}$, the empirical Bayesian method yields a model with a few
covariates that have big coefficients. In comparison, the common
choice has larger $w$ and $\lambda$ values, leading to a regression
model that combines more covariates each having smaller effects. It
turns out these two opposite strategy of modeling both predict the
test dataset well, with pmse being $0.411$ and $0.431$
respectively. For comparison, pmses were calculated for several
frequentist penalized linear regression methods with their respective penalty
parameters chosen by $10$-fold cross validation.  The resulting pmses for 
the ridge, the lasso and the elastic net method are $4.675$, $0.633$ and $0.536$,
respectively. %(see cookie-simple-p=50.html)

The BM method for estimating SE is carried out above without the need of further user input. Theoretically, its competitor RS can be developed too, if enough regeneration times can be identified for each Markov chain. 
%Using the technique introduced in ???, the best regeneration scheme we can come up with among all $12$ Markov chains under consideration regenerated about 10 times every $10^4$ iterations, and that the number of iterations between consecutive regeneration times are highly variable. 
Recall that with the random-swap algorithm, each Markov chain lives on the discrete state space $\Gamma$ of size $2^q$. A naive way to introduce regeneration is to specify a single point $\gamma_1$, then each visit of the Markov chain to $\gamma_1$ marks a regeneration time. Note that the chance of visiting $\gamma_1$ converges to $\pi(\gamma_1)$, the posterior probability of $\gamma_1$. In our BVS model with $2^{50}\approx 1.1\times 10^{15}$ states of $\gamma$, even $\max_{\gamma \in  \Gamma} \pi(\gamma)$ could be very small. Take for example the Markov chain for the BVS model with $h=(w, \log(g))=(.4, 1)$, the point with the highest frequency appeared only $8$ out of a run of $10^4$ iterations. And that the waiting times between consecutive regenerations are highly variable, which ranges from less than ten iterations to a few thousand iterations. 
%Even more troublesome is the fact that all $8$ appearances occurred within a span of 10 consecutive iterations. 
%not enough to yield a stable RS estimator. %Even for other choices of $h$ that a higher frequency, but 
%This is a common problem of identifying regeneration times of Metropolis Hastings algorithms this way, in that the random waiting time between regenerations has huge variability, due to short waiting times resulting from the chain getting stuck at $\gamma_1$, and very long waiting times once the chain departs $\gamma_1$. 
% \blue{Say, after a visit to $\gamma_1$ that counts as a regeneration time, a proposal to another state is rejected, then the Markov chain stays at $\gamma_1$, which constitutes another regeneration. And once the Chain leaves $\gamma_1$ eventually, it takes a long time for it to return. Many of the waiting times is $1$, as the next proposal is rejected, but then it can take a long time for the chain to revisit. }
To obtain alternative schemes of identifying regeneration times, one can take the general minorization condition approach. It could potentially increase the chance of regeneration and reduce variability of the waiting times. Specifically, for any $\alpha \in \{1, 2, \cdots, 2^q\}$, one could define $D_\ap$ to contain the $\ap$ points with the highest posterior probabilities, and find $\epsilon_\ap \in (0,1]$ and a probability mass function $k_\ap(\cdot)$ such that $p(\gamma' | \gamma) \geq \epsilon_\ap I_{D_\ap}(\gamma) k_\ap(\gamma')$ for all $\gamma' \in \Gamma$. Note that as $\ap$ increases, the chance of visiting $D_\ap$ improves, but $\epsilon_\ap$, the conditional rate of regeneration given the current state $\gamma$ is in $D_\ap$, would decrease sharply. Finding a good $\ap$ to maximize the overall chance of regeneration requires tuning that is specific for both the model specification $h$ and the dataset. Even if we can find the optimal $\ap$ for each Markov chain used in the example, it is unlikely that all of them would regenerate often enough for the RS estimator to be stable.

\end{appendix}

%%%%%%%%%%%%%%%%%%%%%%%%%%%%%%%%%%%%%%%%%%%%%%%%%%%%%%%%%%%%%%%%%%%%%%%%%%%%%%%%%%%%%%%%%%%%%%%%%%%%%%%%%%%%%%%%%%%%%%%%%%%%
%\vskip 14pt
%\noindent {\large\bf Supplementary Materials}
%
%The supplement to this paper contains a proof to the extension of the CLT based on regenerative simulation mentioned in Remark~\ref{rem:tdhdiff}. Also included is a simulation study that demonstrates consistency of the BM and the RS estimators in stage~2 of the generalized IS estimators, as well as a comparison among three different weighting strategies. Finally, we study a linear regression model and use the BM estimator to aid the process of empirical Bayes variable selection.  
%%Contain the brief description of the online supplementary materials.
%\par
%%%%%%%%%%%%%%%%%%%%%%%%%%%%%%%%%%%%%%%%%%%%%%%%%%%%%%%%%%%%%%%%%%%%%%%%%%%%%%%%%%%%%%%%%%%%%%%%%%%%%%%%%%%%%%%%%%%%%%%%%%%%
%\vskip 14pt
%\noindent {\large\bf Acknowledgements}

%Write the acknowledgements here.
%\par

% \bibitem[{Zellner(1986)}]{zell:1986}
% \text{Zellner, A.} (1986).
% \newblock {On assessing prior distributions and Bayesian regression analysis
%   with $g$-prior distributions}.
% \newblock \textit{Bayesian Inference and Decision Techniques: Essays in Honor
%   of Bruno de Finetti} \textbf{6} 233--243.
%\end{thebibliography}

\bibliographystyle{apalike}
\vskip .65cm
\noindent
{Department of Statistics, Iowa State University}
\vskip 2pt
\noindent
E-mail: vroy@iastate.edu
\vskip 2pt

\noindent
{Department of Statistics and Actuarial Science, University of Iowa}
\vskip 2pt
\noindent
E-mail: aixin-tan@uiowa.edu

\noindent
{Department of Statistics, University of California, Riverside}
\vskip 2pt
\noindent
E-mail: jflegal@ucr.edu 
\end{document}